\documentclass[leqno]{article}

\usepackage{makeidx}
\makeindex

\def\Int{\mathop{\rm Int}}

\newtheorem{theorem}{Theorem}
\newtheorem{lemma}[theorem]{Lemma}
\newtheorem{proposition}[theorem]{Proposition}
\newtheorem{definition}[theorem]{Definition}
\newtheorem{corollary}[theorem]{Corollary}

\newcommand{\begintheorem}{\addtocounter{equation}{1}\begin{theorem}}
\newcommand{\beginlemma}{\addtocounter{equation}{1}\begin{lemma}}
\newcommand{\beginproposition}{\addtocounter{equation}{1}\begin{proposition}}
\newcommand{\begindefinition}{\addtocounter{equation}{1}\begin{definition}}
\newcommand{\begincorollary}{\addtocounter{equation}{1}\begin{corollary}}

\begin{document}

\title{Some topics related to metrics and norms, including ultrametrics
and ultranorms, 3}

\author{Stephen Semmes \\
        Rice University}

\date{}

\maketitle

\begin{abstract}
Some basic geometric properties related to connectedness and
topological dimension $0$ are discussed, especially in connection
with the ultrametric version of the triangle inequality.
\end{abstract}

\tableofcontents

\part{Semimetrics and uniform structures}
\label{semimetrics, uniform structures}

\section{Semimetrics}
\label{semimetrics}
\setcounter{equation}{0}

        Let $X$ be a set.  A nonnegative real-valued function $d(x, y)$
defined for $x, y \in X$ is said to be a \emph{semimetric}\index{semimetrics}
on $X$ if it satisfies the following three conditions.  First,
\begin{equation}
\label{d(x, x) = 0}
        d(x, x) = 0
\end{equation}
for every $x \in X$.  Second, $d(x, y)$ should be symmetric in $x$ and
$y$, so that
\begin{equation}
\label{d(x, y) = d(y, x)}
        d(x, y) = d(y, x)
\end{equation}
for every $x, y \in X$.  Third, $d(\cdot, \cdot)$ should satisfy the
\emph{triangle inequality},\index{triangle inequality} which is to say
that
\begin{equation}
\label{d(x, z) le d(x, y) + d(y, z)}
        d(x, z) \le d(x, y) + d(y, z)
\end{equation}
for every $x, y z \in X$.  If $d(\cdot, \cdot)$ also has the property
that
\begin{equation}
\label{d(x, y) > 0}
        d(x, y) > 0
\end{equation}
for every $x, y \in X$ with $x \ne y$, then $d(\cdot, \cdot)$ is said
to be a \emph{metric}\index{metrics} on $X$.  The \emph{discrete
  metric}\index{discrete metric} on $X$ is defined by putting $d(x,
y)$ equal to $1$ when $x \ne y$ and to $0$ when $x = y$, and it is
easy to see that this defines a metric on $X$.

        Let $d(\cdot, \cdot)$ be any semimetric on $X$.  The \emph{open
ball}\index{open balls} in $X$ centered at a point $x \in X$ with
radius $r > 0$ associated to $d(\cdot, \cdot)$ is defined as usual by
\begin{equation}
\label{B(x, r) = B_d(x, r) = {y in X : d(x, y) < r}}
        B(x, r) = B_d(x, r) = \{y \in X : d(x, y) < r\}.
\end{equation}
Similarly, the \emph{closed ball}\index{closed balls} in $X$ centered
at $x \in X$ with radius $r \ge 0$ associated to $d(\cdot, \cdot)$
is defined by
\begin{equation}
\label{overline{B}(x, r) = overline{B}_d(x, r) = {y in X : d(x, y) le r}}
 \overline{B}(x, r) = \overline{B}_d(x, r) = \{y \in X : d(x, y) \le r\}.
\end{equation}
If $y \in B(x, r)$, then
\begin{equation}
\label{t = r - d(x, y)}
        t = r - d(x, y)
\end{equation}
is strictly positive, and one can check that
\begin{equation}
\label{B(y, t) subseteq B(x, r)}
        B(y, t) \subseteq B(x, r),
\end{equation}
using the triangle inequality.  In the same way, if $y \in
\overline{B}(x, r)$, then (\ref{t = r - d(x, y)}) is greater than or
equal to $0$, and
\begin{equation}
\label{overline{B}(y, t) subseteq overline{B}(x, r)}
        \overline{B}(y, t) \subseteq \overline{B}(x, r).
\end{equation}

        Let us say that $W \subseteq X$ is an \emph{open set}\index{open sets}
with respect to $d(\cdot, \cdot)$ if for each $x \in W$ there is an
$r > 0$ such that
\begin{equation}
\label{B(x, r) subseteq W}
        B(x, r) \subseteq W.
\end{equation}
The collection of open subsets of $X$ with respect to $d(\cdot,
\cdot)$ defines a topology on $X$, as usual.  Every open ball in $X$
with respect to $d(\cdot, \cdot)$ is an open set with respect to the
topology determined on $X$ by $d(\cdot, \cdot)$, by (\ref{B(y, t)
  subseteq B(x, r)}).  If $d(\cdot, \cdot)$ is a metric on $X$, then
$X$ is Hausdorff with respect to the topology determined by $d(\cdot,
\cdot)$.  If $d(\cdot, \cdot)$ is the discrete metric on $X$, then the
corresponding topology on $X$ is the same as the discrete topology, in
which every subset of $X$ is an open set.

        Let $d(\cdot, \cdot)$ be any semimetric on $X$ again, and put
\begin{equation}
\label{V(x, r) = X setminus overline{B}(x, r) = {y in X : d(x, y) > r}}
 V(x, r) = X \setminus \overline{B}(x, r) = \{y \in X : d(x, y) > r\}
\end{equation}
for every $x \in X$ and $r \ge 0$.  If $y \in V(x, r)$, then
\begin{equation}
\label{t' = d(x, y) - r > 0}
        t' = d(x, y) - r > 0,
\end{equation}
and one can verify that
\begin{equation}
\label{B(y, t') subseteq V(x, r)}
        B(y, t') \subseteq V(x, r),
\end{equation}
using the triangle inequality.  This implies that $V(x, r)$ is an open
set in $X$ with respect to the topology determined by $d(\cdot,
\cdot)$.  Equivalently, this means that $\overline{B}(x, r)$ is a
closed set in $X$ with respect to this topology.  Alternatively, one
can use the triangle inequality to check that $\overline{B}(x, r)$
contains all of its limit points with respect to this topology.

\section{Relations}
\label{relations}
\setcounter{equation}{0}

        Let $X$ be a set again, and let $X \times X$ be the Cartesian
product of $X$ with itself, which is the set of ordered pairs $(x, y)$
with $x, y \in X$.  Thus a subset of $X \times X$ is the same as a
(binary) \emph{relation}\index{relations} on $X$.  The \emph{identity
  relation}\index{identity relation} on $X$ corresponds to the
\emph{diagonal}\index{diagonal set} in $X \times X$, defined by
\begin{equation}
\label{Delta = Delta_X = {(x, x) : x in X}}
        \Delta = \Delta_X = \{(x, x) : x \in X\}.
\end{equation}
If $U \subseteq X \times X$ is any relation on $X$, then put
\begin{equation}
\label{widetilde{U} = {(x, y) : (y, x) in U}}
        \widetilde{U} = \{(x, y) : (y, x) \in U\}.
\end{equation}
This is sometimes denoted $U^{-1}$, and described as the inverse
relation associated to $U$.

        If $U, V \subseteq X \times X$ are relations on $X$, then
let $U * V$ be the relation on $X$ defined by
\begin{eqnarray}
\label{U * V = ...}
 U * V & = & \{(x, z) \in X \times X : \hbox{ there is a } y \in X
                  \hbox{ such that} \\
 & & \qquad\qquad\qquad\qquad (x, y) \in U \hbox{ and } (y, z) \in V\}.
                                                                \nonumber
\end{eqnarray}
This is the same as the composition $V \circ U$ of $U$ and $V$ in the
notation and terminology on p7 of \cite{jk}.  Note that
\begin{equation}
\label{U * Delta = Delta * U = U}
        U * \Delta = \Delta * U = U
\end{equation}
for every $U \subseteq X \times X$.  If $U, V, W \subseteq X \times X$
are relations on $X$, then it is easy to see that
\begin{equation}
\label{(U * V) * W = U * (V * W)}
        (U * V) * W = U * (V * W).
\end{equation}
More precisely, both sides of (\ref{(U * V) * W = U * (V * W)})
consist of the $(x, t) \in X \times X$ for which there are $y, z \in
X$ such that $(x, y) \in U$, $(y, z) \in V$, and $(z, t) \in W$.

        Suppose that $d(\cdot, \cdot)$ is a semimetric on $X$, and put
\begin{equation}
\label{U(r) = U_d(r) = {(x, y) in X times X : d(x, y) < r}}
        U(r) = U_d(r) = \{(x, y) \in X \times X : d(x, y) < r\}
\end{equation}
for every $r > 0$.  Observe that
\begin{equation}
\label{Delta subseteq U(r)}
        \Delta \subseteq U(r)
\end{equation}
for each $r > 0$, by (\ref{d(x, x) = 0}), and that
\begin{equation}
\label{widetilde{U(r)} = U(r)}
        \widetilde{U(r)} = U(r)
\end{equation}
for every $r > 0$, by (\ref{d(x, y) = d(y, x)}).  The triangle
inequality (\ref{d(x, z) le d(x, y) + d(y, z)}) implies that
\begin{equation}
\label{U(r_1) * U(r_2) subseteq U(r_1 + r_2)}
        U(r_1) * U(r_2) \subseteq U(r_1 + r_2)
\end{equation}
for every $r_1, r_2 > 0$.  By construction,
\begin{equation}
\label{U(r) subseteq U(t)}
        U(r) \subseteq U(t)
\end{equation}
when $r < t$, and
\begin{equation}
\label{bigcap_{r > 0} U(r) = {(x, y) in X times X : d(x, y) = 0}}
        \bigcap_{r > 0} U(r) = \{(x, y) \in X \times X : d(x, y) = 0\}.
\end{equation}
In particular,
\begin{equation}
\label{bigcap_{r > 0} U(r) = Delta}
        \bigcap_{r > 0} U(r) = \Delta
\end{equation}
exactly when $d(\cdot, \cdot)$ is a metric on $X$.

        If $U \subseteq X \times X$ is a relation on $X$ and $A$
is a subset of $X$, then we put
\begin{equation}
\label{U[A] = {y in X : there is an x in A such that (x, y) in U}}
        U[A] = \{y \in X : \hbox{ there is an } x \in A \hbox{ such that }
                                                           (x, y) \in U\}.
\end{equation}
If $V \subseteq X \times X$ is another relation on $X$, then it is
easy to see that
\begin{equation}
\label{(U * V)[A] = V[U[A]]}
        (U * V)[A] = V[U[A]].
\end{equation}
In this case, $U \cap V$ is a relation on $X$ too, and
\begin{equation}
\label{(U cap V)[A] subseteq (U[A]) cap (V[A])}
        (U \cap V)[A] \subseteq (U[A]) \cap (V[A]).
\end{equation}
If $x \in X$, then it will be convenient to put
\begin{equation}
\label{U[x] = U[{x}] = {y in X : (x, y) in U}}
        U[x] = U[\{x\}] = \{y \in X : (x, y) \in U\}.
\end{equation}
Using this notation, we have equivalently that
\begin{equation}
\label{U[A] = bigcup_{x in A} U[x]}
        U[A] = \bigcup_{x \in A} U[x]
\end{equation}
for each $A \subseteq X$.  Note that
\begin{equation}
\label{(U cap V)[x] = (U[x]) cap (V[x])}
        (U \cap V)[x] = (U[x]) \cap (V[x]),
\end{equation}
which is to say that equality holds in (\ref{(U cap V)[A] subseteq
  (U[A]) cap (V[A])}) when $A$ has only one element.  Let $d(\cdot,
\cdot)$ be a semimetric on $X$ again, and let $U(r)$ be as in
(\ref{U(r) = U_d(r) = {(x, y) in X times X : d(x, y) < r}}) for some
$r > 0$.  Using the notation (\ref{U[x] = U[{x}] = {y in X : (x, y) in
    U}}) with $U = U(r)$, we get that
\begin{equation}
\label{(U(r))[x] = B(x, r)}
        (U(r))[x] = B(x, r)
\end{equation}
for every $x \in X$, where $B(x, r)$ is as in (\ref{B(x, r) = B_d(x,
  r) = {y in X : d(x, y) < r}}).

\section{Uniform structures}
\label{uniform structures}
\setcounter{equation}{0}

        Let $X$ be a set again, and let $\mathcal{U}$ be a nonempty
collection of subsets of $X \times X$.  If $\mathcal{U}$ satisfies the
following five conditions, then $\mathcal{U}$ is said to define a
\emph{uniformity}\index{uniformities} on $X$, and $(X, \mathcal{U})$
is said to be a \emph{uniform space},\index{uniform spaces} as on p176
of \cite{jk}.  First, for each $U \in \mathcal{U}$, we should have
that
\begin{equation}
\label{Delta subseteq U}
        \Delta \subseteq U,
\end{equation}
where $\Delta$ is as in (\ref{Delta = Delta_X = {(x, x) : x in X}}).
Second, for each $U \in \mathcal{U}$, we ask that
\begin{equation}
\label{widetilde{U} in mathcal{U}}
        \widetilde{U} \in \mathcal{U}
\end{equation}
too, where $\widetilde{U}$ is as in (\ref{widetilde{U} = {(x, y) : (y,
    x) in U}}).  Third, for every $U \in \mathcal{U}$, there should be
a $V \in \mathcal{U}$ such that
\begin{equation}
\label{V * V subseteq U}
        V * V \subseteq U,
\end{equation}
where $V * V$ is as defined in (\ref{U * V = ...}).  Fourth, for any
two elements $U$, $V$ of $\mathcal{U}$, we ask that
\begin{equation}
\label{U cap V in mathcal{U}}
        U \cap V \in \mathcal{U}
\end{equation}
as well.  The fifth condition is that if $U \in \mathcal{U}$ and $U
\subseteq V \subseteq X \times X$, then we also have that
\begin{equation}
\label{V in mathcal{U}}
        V \in \mathcal{U}.
\end{equation}
This implies that $X \times X$ should be an element of $\mathcal{U}$,
since $\mathcal{U}$ is supposed to be nonempty.

        Let $d(\cdot, \cdot)$ be a semimetric on $X$, and let $U_d(r)
\subseteq X \times X$ be as in (\ref{U(r) = U_d(r) = {(x, y) in X times X : 
d(x, y) < r}}) for each $r > 0$.  Also let $\mathcal{U}_d$ be the collection
of subsets $U$ of $X \times X$ for which there is an $r > 0$ such that
\begin{equation}
\label{U_d(r) subseteq U}
        U_d(r) \subseteq U.
\end{equation}
It is easy to see that this defines a uniformity on $X$, using
(\ref{Delta subseteq U(r)}), (\ref{widetilde{U(r)} = U(r)}),
(\ref{U(r_1) * U(r_2) subseteq U(r_1 + r_2)}), and (\ref{U(r) subseteq
  U(t)}).  More precisely, one can check that $\mathcal{U}_d$
satisfies the first four requirements of a uniformity using these
properties of $U_d(r)$.  The fifth requirement of a uniformity is
built into the definition of $\mathcal{U}_d$.  If $d(\cdot, \cdot)$ is
the discrete metric on $X$, then
\begin{equation}
\label{U_d(r) = Delta}
        U_d(r) = \Delta
\end{equation}
when $0 < r \le 1$.  In this case, $\mathcal{U}_d$ consists of all
subsets of $X \times X$ that contain the diagonal $\Delta$ as a
subset.  If instead $d(x, y) = 0$ for every $x, y \in X$, then
\begin{equation}
\label{U_d(r) = X times X}
        U_d(r) = X \times X
\end{equation}
for every $r > 0$, and $X \times X$ is the only element of
$\mathcal{U}_d$.

        Let $\mathcal{U}$ be any uniformity on $X$.  A subcollection
$\mathcal{B}$ of $\mathcal{U}$ is said to be a \emph{base}\index{base for
a uniformity} for $\mathcal{U}$ if for each $U \in \mathcal{U}$ there
is a $V \in \mathcal{B}$ such that
\begin{equation}
\label{V subseteq U}
        V \subseteq U.
\end{equation}
In this case, $\mathcal{U}$ is exactly the same as the collection of
subsets $U$ of $X \times X$ for which there is a $V \in \mathcal{B}$
that satisfies (\ref{V subseteq U}), because of the fifth requirement
of a uniformity.  In particular, $\mathcal{U}$ is uniquely determined
by $\mathcal{B}$ under these conditions.  If $d(\cdot, \cdot)$ is a
semimetric on $X$, and if $\mathcal{B}_d$ is the collection of subsets
of $X \times X$ of the form $U_d(r)$ as in (\ref{U(r) = U_d(r) = {(x,
    y) in X times X : d(x, y) < r}}) for some $r > 0$, then
$\mathcal{B}_d$ is a base for the uniformity $\mathcal{U}_d$ on $X$
associated to $d$ as in the preceding paragraph.

        If $\mathcal{B}$ is a base for a uniformity $\mathcal{U}$
on $X$, then $\mathcal{B}$ is a nonempty collection of subsets of $X
\times X$, because $\mathcal{U} \ne \emptyset$.  Each $U \in
\mathcal{B}$ has to contain the diagonal $\Delta$ as a subset, because
of (\ref{Delta subseteq U}).  If $U \in \mathcal{B}$, then there
should be a $V_0 \in \mathcal{B}$ such that
\begin{equation}
\label{V_0 subseteq widetilde{U}}
        V_0 \subseteq \widetilde{U},
\end{equation}
by (\ref{widetilde{U} in mathcal{U}}).  There should also be a $V \in
\mathcal{B}$ that satisfies (\ref{V * V subseteq U}).  If $U, V \in
\mathcal{B}$, then there should be a $W \in \mathcal{B}$ such that
\begin{equation}
\label{W subseteq U cap V}
        W \subseteq U \cap V,
\end{equation}
because of (\ref{U cap V in mathcal{U}}).  Conversely, if
$\mathcal{B}$ is a nonempty collection of subsets of $X \times X$
that satisfies these four conditions, then one can check that
\begin{equation}
\label{mathcal{U}(mathcal{B}) = ...}
 \mathcal{U}(\mathcal{B}) = \{U \subseteq X \times X : \hbox{ there is a }
                           V \in \mathcal{B} \hbox{ such that } V \subseteq U\}
\end{equation}
is a uniformity on $X$.  Note that $\mathcal{B}$ is automatically a
base for (\ref{mathcal{U}(mathcal{B}) = ...}), by construction.  If
$\mathcal{B}$ is a base for any uniformity $\mathcal{U}$ on $X$, then
$\mathcal{U}$ has to be the same as (\ref{mathcal{U}(mathcal{B}) =
  ...}), as in the previous paragraph.

\section{The associated topology}
\label{associated topology}
\setcounter{equation}{0}

        Let $(X, \mathcal{U})$ be a uniform space.  Let us say that
$A \subseteq X$ is an \emph{open set}\index{open sets} with respect
to the uniformity $\mathcal{U}$ if for each $x \in A$ there is a
$U \in \mathcal{U}$ such that
\begin{equation}
\label{U[x] subseteq A}
        U[x] \subseteq A,
\end{equation}
where $U[x]$ is as defined in (\ref{U[x] = U[{x}] = {y in X : (x, y)
    in U}}).  Note that the empty set and $X$ itself are automatically
open sets with respect to $\mathcal{U}$.  The union of any collection
of open subsets of $X$ with respect to $\mathcal{U}$ is clearly open
with respect to $\mathcal{U}$ too.  Suppose that $A, B \subseteq X$ are
open sets with respect to $\mathcal{U}$, and let us check that $A \cap B$
is open with respect to $\mathcal{U}$ as well.  Let $x \in A \cap B$
be given, and let $U$, $V$ be elements of $\mathcal{U}$ such that
(\ref{U[x] subseteq A}) holds and similarly
\begin{equation}
\label{V[x] subseteq B}
        V[x] \subseteq B.
\end{equation}
The fourth condition (\ref{U cap V in mathcal{U}}) in the definition
of a uniformity implies that $U \cap V$ is an element of $\mathcal{U}$,
and we have that
\begin{equation}
\label{(U cap V)[x] = (U[x]) cap (V[x]) subseteq A cap B}
        (U \cap V)[x] = (U[x]) \cap (V[x]) \subseteq A \cap B,
\end{equation}
using (\ref{(U cap V)[x] = (U[x]) cap (V[x])}) in the first step, and
(\ref{U[x] subseteq A}), (\ref{V[x] subseteq B}) in the second step.
This shows that $A \cap B$ is an open set in $X$ with respect to
$\mathcal{U}$, as desired.  It follows that the collection of open
subsets of $X$ with respect to $\mathcal{U}$ defines a topology on
$X$.

        Let $d(\cdot, \cdot)$ be a semimetric on a set $X$, and let
$\mathcal{U}_d$ be the corresponding uniformity on $X$, as in the
previous section.  In this case, the open subsets of $X$ with respect
to $\mathcal{U}_d$ are the same as the open subsets of $X$ with
respect to the topology determined by $d(\cdot, \cdot)$ in the usual
way, as in Section \ref{semimetrics}.  This can be verified using
(\ref{(U(r))[x] = B(x, r)}) and the relevant definitions.  In particular,
if $d(\cdot, \cdot)$ is the discrete metric on $X$, then $\mathcal{U}_d$
consists of all subsets of $X \times X$ that contain the diagonal $\Delta$,
and the corresponding topology on $X$ is the discrete topology.  If
instead $d(x, y) = 0$ for every $x, y \in X$, then $X \times X$
is the only element of $\mathcal{U}_d$, and the corresponding topology
on $X$ is the discrete topology.

        Let $(X, \mathcal{U})$ be any uniform space again, and let
$A$ be a subset of $X$.  Put
\begin{equation}
\label{A_0 = {x in X : there is a U in mathcal{U} such that U[x] subseteq A}}
        A_0 = \{x \in X : \hbox{ there is a } U \in \mathcal{U}
                             \hbox{ such that } U[x] \subseteq A\}.
\end{equation}
Note that
\begin{equation}
\label{x in U[x]}
        x \in U[x]
\end{equation}
for every $x \in X$ and $U \in \mathcal{U}$, because of (\ref{Delta
  subseteq U}).  Thus
\begin{equation}
\label{A_0 subseteq A}
        A_0 \subseteq A
\end{equation}
automatically.  As in Theorem 4 on p178 of \cite{jk},
\begin{equation}
\label{A_0 = Int(A)}
        A_0 = \Int(A),
\end{equation}
where $\Int(A)$ denotes the interior of $A$ with respect to the
topology on $X$ associated to $\mathcal{U}$.  More precisely, it is
easy to see that
\begin{equation}
\label{Int(A) subseteq A_0}
        \Int(A) \subseteq A_0,
\end{equation}
directly from the definitions.  In order to show that (\ref{A_0 =
  Int(A)}) holds, it suffices to verify that $A_0$ is an open set in
$X$ with respect to the topology associated to $\mathcal{U}$.  Let $x
\in A_0$ be given, and let $U$ be an element of $\mathcal{U}$ such
that (\ref{U[x] subseteq A}) holds, as in (\ref{A_0 = {x in X : there
    is a U in mathcal{U} such that U[x] subseteq A}}).  By definition
of a uniformity, there is a $V \in \mathcal{U}$ that satisfies (\ref{V
  * V subseteq U}).  It follows that
\begin{equation}
\label{(V * V)[x] subseteq U[x] subseteq A}
        (V * V)[x] \subseteq U[x] \subseteq A.
\end{equation}
This implies that
\begin{equation}
\label{V[V[x]] = (V * V)[x] subseteq A}
        V[V[x]] = (V * V)[x] \subseteq A,
\end{equation}
using (\ref{(U * V)[A] = V[U[A]]}) in the first step.  Equivalently,
this means that
\begin{equation}
\label{V[y] subseteq A}
        V[y] \subseteq A
\end{equation}
for every $y \in V[x]$, as in (\ref{U[A] = bigcup_{x in A} U[x]}).
Thus $y \in A_0$ for each $y \in V[x]$, so that
\begin{equation}
\label{V[x] subseteq A_0}
        V[x] \subseteq A_0.
\end{equation}
This shows that $A_0$ is an open set in $X$ with respect to the
topology associated to $\mathcal{U}$, as desired, because $x \in A_0$
is arbitrary.

        In particular, if $x$ is any element of $X$ and $U$ is any
element of $\mathcal{U}$, then
\begin{equation}
\label{x in Int(U[x])}
        x \in \Int(U[x]).
\end{equation}
This follows by applying (\ref{A_0 = Int(A)}) to $A = U[x]$, in which
case $x \in A_0$ automatically.  Now let $A \subseteq X$ and $U \in
\mathcal{U}$ be given, and put
\begin{equation}
\label{B = bigcup_{x in A} Int(U[x])}
        B = \bigcup_{x \in A} \Int(U[x]).
\end{equation}
Thus $B$ is automatically an open set in $X$ with respect to the
topology associated to $\mathcal{U}$, since it is a union of open
sets.  We also have that
\begin{equation}
\label{A subseteq B subseteq U[A]}
        A \subseteq B \subseteq U[A],
\end{equation}
where $U[A]$ is defined in (\ref{U[A] = {y in X : there is an x in A
    such that (x, y) in U}}).  More precisely, the first inclusion in
(\ref{A subseteq B subseteq U[A]}) uses (\ref{x in Int(U[x])}).  The
second inclusion in (\ref{A subseteq B subseteq U[A]}) uses the
equivalent expression for $U[A]$ in (\ref{U[A] = bigcup_{x in A}
  U[x]}), and the fact that the interior of $U[x]$ is automatically
contained in $U[x]$.

\section{Closure and regularity}
\label{closure, regularity}
\setcounter{equation}{0}

        Let $(X, \mathcal{U})$ be a uniform space, let $A$ be a subset
of $X$, and let $\overline{A}$ be the closure of $A$ in $X$ with respect
to the topology associated to $\mathcal{U}$.  Let us check that
\begin{equation}
\label{overline{A} = ...}
 \overline{A} = \{x \in X : U[x] \cap A \ne \emptyset \hbox{ for every }
                                             U \in \mathcal{U}\},
\end{equation}
where $U[x]$ is as in (\ref{U[x] = U[{x}] = {y in X : (x, y) in U}}).
Of course, $\overline{A}$ can normally be defined as the set of $x \in
X$ such that every open subset of $X$ that contains $x$ also
intersects $A$.  The fact that the right side of (\ref{overline{A} =
  ...})  is contained in $\overline{A}$ thus follows directly from the
definition of the topology on $X$ associated to $\mathcal{U}$.  The
opposite inclusion can be obtained similarly, using (\ref{x in
  Int(U[x])}).

        Now let us use (\ref{overline{A} = ...}) to verify that
\begin{equation}
\label{overline{A} = bigcap_{V in mathcal{U}} V[A]}
        \overline{A} = \bigcap_{V \in \mathcal{U}} V[A],
\end{equation}
where $V[A]$ is as in (\ref{U[A] = {y in X : there is an x in A such
    that (x, y) in U}}).  This is the same as saying that
\begin{equation}
\label{overline{A} = bigcap_{U in mathcal{U}} widetilde{U}[A]}
        \overline{A} = \bigcap_{U \in \mathcal{U}} \widetilde{U}[A],
\end{equation}
where $\widetilde{U}$ is as defined in (\ref{widetilde{U} = {(x, y) :
    (y, x) in U}}).  The equivalence of (\ref{overline{A} = bigcap_{V
    in mathcal{U}} V[A]}) and (\ref{overline{A} = bigcap_{U in
    mathcal{U}} widetilde{U}[A]}) uses the second condition
(\ref{widetilde{U} in mathcal{U}}) in the definition of a uniformity,
and the simple fact that
\begin{equation}
\label{widetilde{(widetilde{U})} = U}
        \widetilde{(\widetilde{U})} = U
\end{equation}
for every $U \subseteq X \times X$.  It is easy to see that
\begin{equation}
\label{widetilde{U}[A] = {x in X : U[x] cap A ne emptyset}}
        \widetilde{U}[A] = \{x \in X : U[x] \cap A \ne \emptyset\}
\end{equation}
for every $U \subseteq X \times X$, directly from the definitions.
Thus (\ref{overline{A} = bigcap_{U in mathcal{U}} widetilde{U}[A]}) is
a reformulation of (\ref{overline{A} = ...}), as desired.

        Remember that a topological space is said to be
\emph{regular}\index{regular topological spaces} in the strict sense
if for each point $x$ and closed set $E$ in the space with $x \not\in E$
there are disjoint open sets that contain $x$ and $E$.  This is equivalent
to saying that for each point $x$ in the space and open set $W$ that
contains $x$ there is an open set $W_1$ that contains $x$ and whose
closure is contained in $W$.  If the space also satisfies the first
or even $0$th separation condition, then it follows that the space is
Hausdorff.  Sometimes one may include the first or $0$th separation
condition in the definition of regularity of a topological space,
and so we refer to regularity in the strict sense to indicate that the
first or $0$th separation condition is not necessarily included.
One may also say that a topological space satisfies the \emph{third
separation condition}\index{separation conditions} when it is regular
in the strict sense and satisfies the first or $0$th separation condition,
and hence is Hausdorff.

        Let $(X, \mathcal{U})$ be a uniform space again, and let us
check that $X$ is regular in the strict sense with respect to the
topology associated to $\mathcal{U}$.  Let $x \in X$ be given, and
suppose that $A \subseteq X$ is an open set with $x \in A$.  Thus
there is a $U \in \mathcal{U}$ such that
\begin{equation}
\label{U[x] subseteq A, 2}
        U[x] \subseteq A,
\end{equation}
as in the previous section.  The third condition in the definition of
a uniformity implies that there is a $V \in \mathcal{U}$ that satisfies
(\ref{V * V subseteq U}), so that
\begin{equation}
\label{V[V[x]] = (V * V)[x] subseteq U[x] subseteq A}
        V[V[x]] = (V * V)[x] \subseteq U[x] \subseteq A.
\end{equation}
It follows that
\begin{equation}
\label{overline{V[x]} subseteq V[V[x]] subseteq A}
        \overline{V[x]} \subseteq V[V[x]] \subseteq A,
\end{equation}
using (\ref{overline{A} = bigcap_{V in mathcal{U}} V[A]}) applied to
$V[x]$ in the first step.  In particular,
\begin{equation}
\label{overline{Int(V[x])} subseteq A}
        \overline{\Int(V[x])} \subseteq A,
\end{equation}
where $\Int(V[x])$ is the interior of $V[x]$ in $X$, as before.  Of
course, $\Int(V[x])$ is an open set in $X$ that contains $x$, by
(\ref{x in Int(U[x])}).  Hence (\ref{overline{Int(V[x])} subseteq A})
implies that $X$ is regular in the strict sense, as desired.

        Here is a slightly different version of the same type of argument.
Let $x \in X$ and a closed set $E \subseteq X$ be given, with $x
\not\in E$.  Thus $A = X \setminus E$ is an open set that contains
$x$, and so there is a $U \in \mathcal{U}$ that satisfies (\ref{U[x]
  subseteq A, 2}) as before.  This leads to a $V \in \mathcal{U}$ that
satisfies (\ref{V[V[x]] = (V * V)[x] subseteq U[x] subseteq A}),
which can be reformulated as saying that
\begin{equation}
\label{V[x] cap (widetilde{V}[E]) = emptyset}
        V[x] \cap (\widetilde{V}[E]) = \emptyset.
\end{equation}
Remember that $\widetilde{V} \in \mathcal{U}$ too, by the second
condition (\ref{widetilde{U} in mathcal{U}}) in the definition of a
uniformity.  It follows that there is an open subset of $X$ that
contains $E$ and is contained in $\widetilde{V}[E]$, as in (\ref{B =
  bigcup_{x in A} Int(U[x])}) and (\ref{A subseteq B subseteq U[A]}).
We have also seen that $x \in \Int(V[x])$, so that (\ref{V[x] cap
  (widetilde{V}[E]) = emptyset}) implies that $x$ and $E$ are
contained in disjoint open subsets of $X$.

        Now let $x, y \in X$ be given, with $x \ne y$, and suppose that
there is a $U \in \mathcal{U}$ such that
\begin{equation}
\label{y not in U[x]}
        y \not\in U[x].
\end{equation}
Let $V$ be an element of $\mathcal{U}$ that satisfies (\ref{V * V
  subseteq U}), so that
\begin{equation}
\label{V[V[x]] = (V * V)[x] subseteq U[x]}
        V[V[x]] = (V * V)[x] \subseteq U[x],
\end{equation}
using (\ref{(U * V)[A] = V[U[A]]}) in the first step.  Thus (\ref{y
  not in U[x]}) implies that $y \not\in V[V[x]]$, which is the
same as saying that
\begin{equation}
\label{V[x] cap (widetilde{V}[y]) = emptyset}
        V[x] \cap (\widetilde{V}[y]) = \emptyset,
\end{equation}
where $\widetilde{V}$ is as in (\ref{widetilde{U} = {(x, y) : (y, x)
    in U}}).  In particular, (\ref{V[x] cap (widetilde{V}[y]) =
  emptyset}) implies that the interiors of $V[x]$ and
$\widetilde{V}[x]$ are disjoint.  Of course, we already know that $x$
and $y$ are contained in the interiors of $V[x]$ and
$\widetilde{V}[y]$, respectively, as in (\ref{x in Int(U[x])}).  This
is basically a more direct version of the proof that $X$ is Hausdorff
when $X$ satisfies the first or $0$th separation condition, since
we already know that $X$ is regular in the strict sense.  This can
also be formulated as saying that $X$ is Hausdorff when
\begin{equation}
\label{bigcap_{U in mathcal{U}} U = Delta}
        \bigcap_{U \in \mathcal{U}} U = \Delta,
\end{equation}
where $\Delta$ is the diagonal (\ref{Delta = Delta_X = {(x, x) : x in
    X}}) in $X \times X$.  Conversely, if $X$ satisfies the first
separation condition, then it is easy to see that (\ref{bigcap_{U in
    mathcal{U}} U = Delta}) holds.  This also works with the $0$th
separation condition, using the second condition (\ref{widetilde{U} in
  mathcal{U}}) in the definition of a uniformity.

        This discussion applies in particular to the case where the
topology on $X$ is determined by a semimetric $d(\cdot, \cdot)$, as in
Section \ref{semimetrics}.  In this case, regularity in the strict
sense can be derived from the fact that closed balls are closed sets,
or equivalently that subsets of $X$ as in (\ref{V(x, r) = X setminus
  overline{B}(x, r) = {y in X : d(x, y) > r}}) are open.  Of course,
one can check directly that $X$ is Hausdorff when $d(\cdot, \cdot)$
is a metric on $X$.

\section{Symmetry and interior}
\label{symmetry, interior}
\setcounter{equation}{0}

        Let $X$ be a set.  As usual, a set $U \subseteq X \times X$
is said to be \emph{symmetric}\index{symmetric sets} if
\begin{equation}
\label{widetilde{U} = U}
        \widetilde{U} = U,
\end{equation}
where $\widetilde{U}$ is as in (\ref{widetilde{U} = {(x, y) : (y, x)
    in U}}).  If $U$ is any subset of $X \times X$, then
\begin{equation}
\label{V = U cap widetilde{U}}
        V = U \cap \widetilde{U}
\end{equation}
is automatically symmetric.  If $\mathcal{U}$ is a uniformity on $X$
and $U \in \mathcal{U}$, then (\ref{V = U cap widetilde{U}}) is an
element of $\mathcal{U}$ as well.  It follows that the symmetric
elements of $\mathcal{U}$ form a base for $\mathcal{U}$.

        If $U$ is any subset of $X \times X$ and $x \in X$, then
$U[x]$ was defined in (\ref{U[x] = U[{x}] = {y in X : (x, y) in U}}).
Thus, for each $y \in X$, we have that
\begin{equation}
\label{widetilde{U}[y] = ... = {x in X : (x, y) in U}}
        \widetilde{U}[y] = \{x \in X : (y, x) \in \widetilde{U}\}
                         = \{x \in X : (x, y) \in U\}.
\end{equation}
If $V$ is another subset of $X \times X$, then $U * V$ is
defined in (\ref{U * V = ...}), and can be given equivalently as
\begin{equation}
\label{U * V = bigcup_{y in X} widetilde{U}[y] times V[y]}
        U * V = \bigcup_{y \in X} \widetilde{U}[y] \times V[y].
\end{equation}
If $W \subseteq X \times X$ too, then
\begin{equation}
\label{(U * V) * W = U * (V * W) = ...}
 (U * V) * W = U * (V * W) = \bigcup_{(y, z) \in V} \widetilde{U}[y] \times W[z],
\end{equation}
as in (\ref{(U * V) * W = U * (V * W)}).  This is a reformulation of
Lemma 1 on p176 of \cite{jk}, and one may drop the parentheses on the
left side, because of associativity.

        Let $\mathcal{U}$ be a uniformity on $X$ again, let
$U \in \mathcal{U}$ be given, and let $V$ be an element of $\mathcal{U}$
that satisfies (\ref{V * V subseteq U}).  This implies that
\begin{equation}
\label{bigcup_{y in X} widetilde{V}[y] times V[y] subseteq U}
        \bigcup_{y \in X} \widetilde{V}[y] \times V[y] \subseteq U,
\end{equation}
as in (\ref{U * V = bigcup_{y in X} widetilde{U}[y] times V[y]}).
Remember that $y$ is an element of the interior of $V[y]$ with respect
to the topology on $X$ associated to $\mathcal{U}$ for each $y \in X$,
as in (\ref{x in Int(U[x])}).  Similarly, $y$ is an element of the
interior of $\widetilde{V}[y]$ for every $y \in X$, because
$\widetilde{V} \in \mathcal{U}$ too, as in (\ref{widetilde{U} in
  mathcal{U}}).  It follows that $(y, y)$ is an element of the
interior of $\widetilde{V}[y] \times V[y]$ with respect to the
corresponding product topology on $X \times X$ for every $y \in X$.
This shows that the diagonal $\Delta$ is contained in the interior of
$U$ with respect to the product topology on $X \times X$.  As a
refinement of this argument, one can first find a $W \in \mathcal{U}$
such that
\begin{equation}
\label{W * W * W subseteq U}
        W * W * W \subseteq U.
\end{equation}
Using (\ref{(U * V) * W = U * (V * W) = ...}), we get that
\begin{equation}
\label{bigcup_{(y, z) in W} widetilde{W}[y] times W[z] subseteq U}
        \bigcup_{(y, z) \in W} \widetilde{W}[y] \times W[z] \subseteq U,
\end{equation}
so that $W$ is contained in the interior of $U$ with respect to the
product topology on $X \times X$.

        Suppose that $E \subseteq X \times X$, and let $\overline{E}$
be the closure of $E$ with respect to the product topology on
$X \times X$, using the topology on $X$ associated to $\mathcal{U}$.
One can check that
\begin{equation}
\label{overline{E} = bigcap_{V_1, V_2 in mathcal{U}} V_1 * E * V_2}
        \overline{E} = \bigcap_{V_1, V_2 \in \mathcal{U}} V_1 * E * V_2,
\end{equation}
in analogy with (\ref{overline{A} = bigcap_{V in mathcal{U}} V[A]}),
and using (\ref{(U * V) * W = U * (V * W) = ...}).  In particular, if
$W \in \mathcal{U}$ is as in (\ref{W * W * W subseteq U}), then we get
that
\begin{equation}
\label{overline{W} subseteq U}
        \overline{W} \subseteq U.
\end{equation}

\section{Compactness}
\label{compactness}
\setcounter{equation}{0}

        Let $(X, \mathcal{U})$ be a uniform space, and suppose that
$K \subseteq X$ is compact with respect to the topology on $X$ associated
to $\mathcal{U}$ as in Section \ref{associated topology}.  Let $W
\subseteq X$ be an open set with respect to the associated topology
that contains $K$ as a subset.  We would like to check that there is a
$V \in \mathcal{U}$ such that
\begin{equation}
\label{V[K] subseteq W}
        V[K] \subseteq W,
\end{equation}
where $V[K]$ is as defined in (\ref{U[A] = {y in X : there is an x in
    A such that (x, y) in U}}).  Because $K \subseteq W$ and $W$ is an
open set in $X$, we have that for each $x \in K$ there is a $U_x \in
\mathcal{U}$ such that
\begin{equation}
\label{U_x[x] subseteq W}
        U_x[x] \subseteq W,
\end{equation}
as in (\ref{U[x] subseteq A}).  Using the definition of a uniform
space, we get that for each $x \in K$ there is a $V_x \in \mathcal{U}$
that satisfies
\begin{equation}
\label{V_x * V_x subseteq U_x}
        V_x * V_x \subseteq U_x,
\end{equation}
as in (\ref{V * V subseteq U}).  Thus
\begin{equation}
\label{V_x[V_x[x]] = (V_x * V_x)[x] subseteq U_x[x] subseteq W}
        V_x[V_x[x]] = (V_x * V_x)[x] \subseteq U_x[x] \subseteq W
\end{equation}
for every $x \in K$, using (\ref{(U * V)[A] = V[U[A]]}) in the first
step.  Remember that $x$ is an element of the interior of $V_x[x]$ in
$X$ for each $x \in K$, as in (\ref{x in Int(U[x])}).  It follows that
there are finitely many elements $x_1, \ldots, x_n$ of $K$ such that
\begin{equation}
\label{K subseteq bigcup_{j = 1}^n V_{x_j}[x_j]}
        K \subseteq \bigcup_{j = 1}^n V_{x_j}[x_j],
\end{equation}
because $K$ is compact in $X$.  Put
\begin{equation}
\label{V = bigcap_{j = 1}^n V_{x_j}}
        V = \bigcap_{j = 1}^n V_{x_j},
\end{equation}
which is an element of $\mathcal{U}$, since $V_{x_j} \in \mathcal{U}$
for each $j = 1, \ldots, n$.  Observe that
\begin{equation}
\label{V[V_{x_j}[x_j]] subseteq V_{x_j}[V_{x_j}[x_j]] subseteq W}
        V[V_{x_j}[x_j]] \subseteq V_{x_j}[V_{x_j}[x_j]] \subseteq W
\end{equation}
for each $j = 1, \ldots, n$, using the fact that $V \subseteq V_{x_j}$
for each $j$ in the first step, and (\ref{V_x[V_x[x]] = (V_x * V_x)[x]
  subseteq U_x[x] subseteq W}) in the second step.  Combining this
with (\ref{K subseteq bigcup_{j = 1}^n V_{x_j}[x_j]}), we get that
\begin{equation}
\label{V[K] subseteq bigcup_{j = 1}^n V[V_{x_j}[x_j]] subseteq W}
        V[K] \subseteq \bigcup_{j = 1}^n V[V_{x_j}[x_j]] \subseteq W,
\end{equation}
as desired.

        If $V$ is any element of $\mathcal{U}$, then there is a
$V_1 \in \mathcal{U}$ such that
\begin{equation}
\label{V_1 * V_1 subseteq V}
        V_1 * V_1 \subseteq V,
\end{equation}
as in (\ref{V * V subseteq U}).  This implies that
\begin{equation}
\label{V_1[V_1[K]] = (V_1 * V_1)[K] subseteq V[K]}
        V_1[V_1[K]] = (V_1 * V_1)[K] \subseteq V[K]
\end{equation}
for any $K \subseteq X$, using (\ref{(U * V)[A] = V[U[A]]}) in the
first step.  It follows that
\begin{equation}
\label{overline{V_1[K]} subseteq V_1[V_1[K]] subseteq V[K]}
        \overline{V_1[K]} \subseteq V_1[V_1[K]] \subseteq V[K],
\end{equation}
using (\ref{overline{A} = bigcap_{V in mathcal{U}} V[A]}) in the first
step.  If $V$, $K$, and $W$ are as in (\ref{V[K] subseteq W}), then we
get that
\begin{equation}
\label{overline{V_1[K]} subseteq W}
        \overline{V_1[K]} \subseteq W.
\end{equation}
Remember that $K$ is automatically contained in the interior of
$V_1[K]$ for every $V_1 \in \mathcal{U}$, as in (\ref{A subseteq B
  subseteq U[A]}).

        Suppose now that $X$ is a topological space which is regular
in the strict sense, as in Section \ref{closure, regularity}.
If $K \subseteq X$ is compact, $W \subseteq X$ is an open set,
and $K \subseteq W$, then there is an open subset of $X$ that contains
$X$ and whose closure is contained in $W$.  More precisely, one can
first use regularity to cover $K$ by open sets whose closures are
contained in $W$, and then use compactness to reduce to a finite
subcover.

\section{Semi-ultrametrics}
\label{semi-ultrametrics}
\setcounter{equation}{0}

        A semimetric $d(\cdot, \cdot)$ on a set $X$ is said to be a
\emph{semi-ultrametric}\index{semi-ultrametrics} on $X$ if it satisfies
\begin{equation}
\label{d(x, z) le max(d(x, y), d(y, z))}
        d(x, z) \le \max(d(x, y), d(y, z))
\end{equation}
for every $x, y, z \in X$.  Of course, (\ref{d(x, z) le max(d(x, y),
  d(y, z))}) automatically implies the ordinary version (\ref{d(x, z)
  le d(x, y) + d(y, z)}) of the triangle inequality.  If a
semi-ultrametric $d(\cdot, \cdot)$ on $X$ is also a metric on $X$, so
that $d(\cdot, \cdot)$ satisfies (\ref{d(x, y) > 0}) too, then
$d(\cdot, \cdot)$ is said to be an
\emph{ultrametric}\index{ultrametrics} on $X$.  It is easy to see that
the discrete metric on any set is an ultrametric.

        Suppose that $d(\cdot, \cdot)$ is a semi-ultrametric on a set $X$,
and let $r > 0$ be given.  Observe that
\begin{equation}
\label{d(x, y) < r}
        d(x, y) < r
\end{equation}
defines an equivalence relation on $X$.  The corresponding equivalence
class in $X$ containing a point $x \in X$ is the same as the open ball
$B(x, r)$ centered at $x$ with radius $r$ with respect to $d(\cdot,
\cdot)$, as in (\ref{B(x, r) = B_d(x, r) = {y in X : d(x, y) < r}}).
If $x, y \in X$ satisfy (\ref{d(x, y) < r}), then the equivalence
classes containing $x$ and $y$ are the same, so that
\begin{equation}
\label{B(x, r) = B(y, r)}
        B(x, r) = B(y, r).
\end{equation}
Similarly, any two open balls in $X$ of radius $r$ with respect to
$d(\cdot, \cdot)$ are either the same or disjoint as subsets of $X$,
because $X$ is partitioned by these equivalence classes.  In
particular, the complement of any open ball in $X$ of radius $r$ with
respect to $d(\cdot, \cdot)$ can be expressed as a union of other open
balls of radius $r$.  This implies that for every $x \in X$,
\begin{equation}
\label{B(x, r) is a closed set in X}
        B(x, r) \hbox{ is a closed set in } X
\end{equation}
with respect to the usual topology on $X$ associated to $d(\cdot,
\cdot)$, because its complement is an open set.

        In the same way,
\begin{equation}
\label{d(x, y) le r}
        d(x, y) \le r
\end{equation}
defines an equivalence relation on $X$ for each $r \ge 0$.  As before,
the corresponding equivalence class in $X$ containing a point $x \in
X$ is the same as the closed ball $\overline{B}(x, r)$ centered at $x$
with radius $r$ with respect to $d(\cdot, \cdot)$, which was defined
in (\ref{overline{B}(x, r) = overline{B}_d(x, r) = {y in X : d(x, y)
    le r}}).  If $x, y \in X$ satisfy (\ref{d(x, y) le r}), then the
equivalence classes containing $x$ and $y$ are the same, which means
that
\begin{equation}
\label{overline{B}(x, r) = overline{B}(y, r)}
        \overline{B}(x, r) = \overline{B}(y, r).
\end{equation}
It follows that
\begin{equation}
\label{overline{B}(x, r) is an open set in X}
        \overline{B}(x, r) \hbox{ is an open set in } X
\end{equation}
for every $x \in X$ when $r > 0$, with respect to the usual topology
on $X$ associated to $d(\cdot, \cdot)$.  Note that if $r = 0$, then
(\ref{d(x, y) le r}) defines an equivalence relation on $X$ for every
semi-metric $d(\cdot, \cdot)$ on $X$.

        Now let $x \sim y$ be any equivalence relation on a set $X$,
and define $d(x, y)$ for $x, y \in X$ by
\begin{eqnarray}
\label{d(x, y) = 0 when x sim y, = 1 when x not sim y}
        d(x, y) & = & 0 \quad\hbox{when } x \sim y \\
                & = & 1 \quad\hbox{when } x \not\sim y. \nonumber
\end{eqnarray}
One can check that this defines a semi-ultrametric on $X$, which we
shall call the \emph{discrete semi-ultrametric}\index{discrete
  semi-ultrametrics} associated to this equivalence relation.  This
reduces to the discrete metric on $X$ when this equivalence relation
is simply equality of elements of $X$.  If $d(x, y)$ is any semimetric
on $X$ that only takes the values $0$ or $1$, then it is easy to see
that $d(x, y)$ is a semi-ultrametric on $X$.  In this case, $d(x, y)$
is the same as the discrete semi-ultrametric associated to the
equivalence relation (\ref{d(x, y) < r}) for any $0 < r \le 1$, which
is the same as the equivalence relation (\ref{d(x, y) le r}) for any
$0 \le r < 1$.

\section{Uniform continuity}
\label{uniform continuity}
\setcounter{equation}{0}

        Let $(X, \mathcal{U})$ and $(Y, \mathcal{V})$ be uniform spaces,
and let $f$ be a mapping from $X$ into $Y$.  This leads to a mapping
$f_2 : X \times X \to Y \times Y$ defined by
\begin{equation}
\label{f_2(x, x') = (f(x), f(x'))}
        f_2(x, x') = (f(x), f(x'))
\end{equation}
for every $x, x' \in X$.  We say that $f$ is \emph{uniformly
  continuous}\index{uniform continuity} if for each $V \in
\mathcal{V}$ we have that
\begin{equation}
\label{f_2^{-1}(V) in mathcal{U}}
        f_2^{-1}(V) \in \mathcal{U}.
\end{equation}
Equivalently, this means that for each $V \in \mathcal{V}$ there is a
$U \in \mathcal{U}$ such that
\begin{equation}
\label{f_2(U) subseteq V}
        f_2(U) \subseteq V.
\end{equation}
If $\mathcal{U}$ and $\mathcal{V}$ are determined by semimetrics on
$X$ and $Y$, respectively, as in Section \ref{uniform structures},
then uniform continuity can also be characterized equivalently in
terms of $\epsilon$'s and $\delta$'s in the usual way.  In the case of
arbitrary uniformities, it suffices to verify (\ref{f_2^{-1}(V) in
  mathcal{U}}) for all $V$ in a base for $\mathcal{V}$, in order to
show that $f$ is uniformly continuous.  There is an analogous
statement for sub-bases for $\mathcal{V}$, which will be defined in
Section \ref{sub-bases}.

        Let $U \subseteq X \times X$ be given, and put
\begin{equation}
\label{U[x] = {x' in X : (x, x') in U}}
        U[x] = \{x' \in X : (x, x') \in U\}
\end{equation}
for each $x \in X$, as in (\ref{U[x] = U[{x}] = {y in X : (x, y) in
    U}}).  Thus
\begin{equation}
\label{U = bigcup_{x in X} ({x} times U[x])}
        U = \bigcup_{x \in X} (\{x\} \times U[x]),
\end{equation}
which implies that
\begin{equation}
\label{f_2(U) = bigcup_{x in X} ({f(x)} times f(U[x]))}
        f_2(U) = \bigcup_{x \in X} (\{f(x)\} \times f(U[x])).
\end{equation}
If $V \subseteq Y \times Y$, then it follows that (\ref{f_2(U)
  subseteq V}) holds if and only if
\begin{equation}
\label{f(U[x]) subseteq V[f(x)]}
        f(U[x]) \subseteq V[f(x)]
\end{equation}
for every $x \in X$.

        Let $X$ and $Y$ be equipped with the topologies associated
to the uniformities $\mathcal{U}$ and $\mathcal{V}$, respectively, as
in Section \ref{associated topology}.  One can check that a mapping $f
: X \to Y$ is continuous at a point $x \in X$ with respect to these
topologies on $X$ and $Y$ if and only if for every $V \in \mathcal{V}$
there is a $U \in \mathcal{U}$ such that (\ref{f(U[x]) subseteq
  V[f(x)]}) holds.  This uses (\ref{x in Int(U[x])}) and its analogue
in $Y$, in addition to the basic definitions.  If $f$ is uniformly
continuous, then it follows that $f$ is continuous in the ordinary
sense, because of the relationship between (\ref{f_2(U) subseteq V})
and (\ref{f(U[x]) subseteq V[f(x)]}).

        As a slight variant of this argument, remember that $f : X \to Y$
is continuous if and only if for every open subset $W$ of $Y$,
$f^{-1}(W)$ is an open subset of $X$.  In this situation, $f^{-1}(W)$
is an open subset of $X$ if and only if for each $x \in f^{-1}(W)$
there is a $U \in \mathcal{U}$ such that
\begin{equation}
\label{U[x] subseteq f^{-1}(W)}
        U[x] \subseteq f^{-1}(W).
\end{equation}
Of course, if $x \in f^{-1}(W)$, then $f(x) \in W$.  Because $W$ is
supposed to be an open set in $Y$, there should be a $V \in
\mathcal{V}$ such that
\begin{equation}
\label{V[f(x)] subseteq W}
        V[f(x)] \subseteq W.
\end{equation}
If there is a $U \in \mathcal{U}$ that satisfies (\ref{f(U[x])
  subseteq V[f(x)]}), then we get that
\begin{equation}
\label{f(U[x]) subseteq W}
        f(U[x]) \subseteq W,
\end{equation}
which is equivalent to (\ref{U[x] subseteq f^{-1}(W)}).

        As another variant, let $x \in X$ and $V \subseteq Y \times Y$
be given, and observe that
\begin{eqnarray}
\label{f^{-1}(V[f(x)]) = ... = (f_2^{-1}(V))[x]}
        f^{-1}(V[f(x)]) & = & \{x' \in X : f(x') \in V[f(x)]\} \\
 & = & \{x' \in X : (f(x), f(x')) \in V\} = (f_2^{-1}(V))[x]. \nonumber
\end{eqnarray}
This permits one to deal with uniform continuity in terms of
(\ref{f_2^{-1}(V) in mathcal{U}}), as in the proof of Theorem 9 on
p181 of \cite{jk}.

        Let $(Z, \mathcal{W})$ be another uniform space, and let $g$
be a mapping from $Y$ into $Z$.  Also let $g_2$ be the corresponding
mapping from $Y \times Y$ into $Z \times Z$, as in (\ref{f_2(x, x') =
  (f(x), f(x'))}).  The composition $g \circ f$ of $f$ and $g$ maps
$X$ into $Z$, and leads to a mapping $(g \circ f)_2$ from $X \times X$
into $Z \times Z$ as in (\ref{f_2(x, x') = (f(x), f(x'))}) as well.
It is easy to see that
\begin{equation}
\label{(g circ f)_2 = g_2 circ f_2}
        (g \circ f)_2 = g_2 \circ f_2
\end{equation}
as mappings from $X \times X$ into $Z \times Z$.  If $f$ and $g$ are
both uniformly continuous, then one can use (\ref{(g circ f)_2 = g_2
  circ f_2}) that $g \circ f$ is uniformly continuous from $X$ into
$Z$.

        Let us say that $f : X \to Y$ is uniformly continuous
\emph{along}\index{uniform continuity!along a set} a subset $E$ of $X$
if for each $V \in \mathcal{V}$ there is a $U \in \mathcal{U}$ such
that (\ref{f(U[x]) subseteq V[f(x)]}) holds for every $x \in E$.
Equivalently, this means that
\begin{equation}
\label{(f(x), f(x')) in V}
        (f(x), f(x')) \in V
\end{equation}
for every $(x, x') \in U$ such that $x \in E$.  This holds with $E =
X$ if and only if $f$ is uniformly continuous as a mapping from $X$
into $Y$, as before.  If $f$ is uniformly continuous along any set $E
\subseteq X$, then $f$ is continuous as a mapping from $X$ into $Y$ at
every point in $E$, with respect to the topologies associated to the
given uniformities.  This follows from the earlier discussion of the
continuity of $f$ at a point $x \in X$ in terms of (\ref{f(U[x])
  subseteq V[f(x)]}).  This condition also implies that the
restriction of $f$ to $E$ is uniformly continuous with respect to the
uniformity induced on $E$ by $\mathcal{U}$ on $X$, which is defined in
Section \ref{induced uniform structures}.  If $f$ is uniformly
continuous along $E \subseteq X$ and $g : Y \to Z$ is uniformly
continuous along $f(E) \subseteq Y$, then it is easy to see that $g
\circ f$ is uniformly continuous along $E$ too, as a mapping from $X$
into $Z$.

\section{Compactness, continued}
\label{compactness, continued}
\setcounter{equation}{0}

        Let $(X, \mathcal{U})$ and $(Y, \mathcal{V})$ be uniform spaces
again, and suppose that $f$ is a continuous mapping from $X$ into $Y$
with respect to their associated topologies.  If $K \subseteq X$ is
compact, then we would like to verify that $f$ is uniformly continuous
along $K$, as in the preceding section.  In particular, if $X$ is
compact, then it follows that $f$ is uniformly continuous as a mapping
from $X$ into $Y$.

        Let $V \in \mathcal{V}$ be given, and let $V_1$ be an element
of $\mathcal{V}$ such that
\begin{equation}
\label{widetilde{V_1} * V_1 subseteq V}
        \widetilde{V_1} * V_1 \subseteq V.
\end{equation}
More precisely, one might as well choose $V_1 \in \mathcal{V}$ to be
symmetric and to satisfy $V_1 * V_1 \subseteq V$, using
(\ref{widetilde{U} in mathcal{U}}), (\ref{V * V subseteq U}), and
(\ref{U cap V in mathcal{U}}).  Also let $x \in K$ be given, and
remember that $f(x)$ is an element of the interior of $V_1[f(x)]$ in
$Y$, as in (\ref{x in Int(U[x])}).  Because $f$ is continuous at $x$,
there is a $U_{1, x} \in \mathcal{U}$ such that
\begin{equation}
\label{f(U_{1, x}[x]) subseteq V_1[f(x)]}
        f(U_{1, x}[x]) \subseteq V_1[f(x)],
\end{equation}
as in (\ref{f(U[x]) subseteq V[f(x)]}) and discussed in the paragraph
immediately after that.  Let $U_{2, x}$ be an element of $\mathcal{U}$
such that
\begin{equation}
\label{U_{2, x} * U_{2, x} subseteq U_{1, x}}
        U_{2, x} * U_{2, x} \subseteq U_{1, x},
\end{equation}
as in (\ref{V * V subseteq U}) again.  If $K$ is compact in $X$, then
there are finitely many elements $x_1, \ldots, x_n$ of $K$ such that
\begin{equation}
\label{K subseteq bigcup_{j = 1}^n U_{2, x_j}[x_j]}
        K \subseteq \bigcup_{j = 1}^n U_{2, x_j}[x_j].
\end{equation}
This uses the fact that $x$ is an element of the interior of $U_{2,
  x}[x]$ in $X$ for every $x \in K$, as in (\ref{x in Int(U[x])}).
Put
\begin{equation}
\label{U_2 = bigcup_{j = 1}^n U_{2, x_j}}
        U_2 = \bigcup_{j = 1}^n U_{2, x_j},
\end{equation}
which is an element of $\mathcal{U}$, since $U_{2, x_j} \in
\mathcal{U}$ for every $j = 1, \ldots, n$.  We would like to check that
\begin{equation}
\label{(f(x), f(x')) in V, 2}
        (f(x), f(x')) \in V
\end{equation}
for every $x \in K$ and $x' \in X$ such that $(x, x') \in U_2$, as in
(\ref{(f(x), f(x')) in V}).  Let $x \in K$ be given again, and choose
$j \in \{1, \ldots, n\}$ such that
\begin{equation}
\label{x in U_{2, x_j}[x_j]}
        x \in U_{2, x_j}[x_j],
\end{equation}
which is possible by (\ref{K subseteq bigcup_{j = 1}^n U_{2,
    x_j}[x_j]}).  If $x' \in X$ satisfies $(x, x') \in U_2 \subseteq
U_{2, x_j}$, then
\begin{equation}
\label{x' in U_{2, x_j}[U_{2, x_j}[x_j]] = ... subseteq U_{1, x_j}[x_j]}
 x' \in U_{2, x_j}[U_{2, x_j}[x_j]] = (U_{2, x_j} * U_{2, x_j})[x_j]
                                  \subseteq U_{1, x_j}[x_j],
\end{equation}
using (\ref{(U * V)[A] = V[U[A]]}) in the second step, and (\ref{U_{2,
    x} * U_{2, x} subseteq U_{1, x}}) in the third step.  Note that
\begin{equation}
\label{U_{2, x_j} subseteq U_{1, x_j}}
        U_{2, x_j} \subseteq U_{1, x_j},
\end{equation}
because of (\ref{U_{2, x} * U_{2, x} subseteq U_{1, x}}) and the fact
that $\Delta \subseteq U_{2, x_j}$, as in (\ref{Delta subseteq U}).
Thus (\ref{x in U_{2, x_j}[x_j]}) implies that
\begin{equation}
\label{x in U_{1, x_j}[x_j]}
        x \in U_{1, x_j}[x_j].
\end{equation}
Using (\ref{x' in U_{2, x_j}[U_{2, x_j}[x_j]] = ... subseteq U_{1,
    x_j}[x_j]}) and (\ref{x in U_{1, x_j}[x_j]}), we can apply
(\ref{f(U_{1, x}[x]) subseteq V_1[f(x)]}) to $x_j$ to get that
\begin{equation}
\label{f(x), f(x') in V_1[f(x_j)]}
        f(x), f(x') \in V_1[f(x_j)].
\end{equation}
This implies that
\begin{equation}
\label{(f(x), f(x')) in widetilde{V_1} * V_1 subseteq V}
        (f(x), f(x')) \in \widetilde{V_1} * V_1 \subseteq V,
\end{equation}
using (\ref{widetilde{V_1} * V_1 subseteq V}) in the second step.
Thus (\ref{(f(x), f(x')) in V, 2}) holds, as desired.

\section{Totally bounded sets}
\label{totally bounded sets}
\setcounter{equation}{0}

        Let $(X, \mathcal{U})$ be a uniform space.  A subset $E$
of $X$ is said to be \emph{totally bounded}\index{totally bounded sets}
if for each $U \in \mathcal{U}$ there is a finite set $A \subseteq X$
such that
\begin{equation}
\label{E subseteq U[A]}
        E \subseteq U[A],
\end{equation}
where $U[A]$ is as in (\ref{U[A] = {y in X : there is an x in A such
    that (x, y) in U}}).  Equivalently, this means that there are
finitely many elements $x_1, \ldots, x_n$ of $X$ such that
\begin{equation}
\label{E subseteq bigcup_{j = 1}^n U[x_j]}
        E \subseteq \bigcup_{j = 1}^n U[x_j],
\end{equation}
where $U[x_j]$ is as in (\ref{U[x] = U[{x}] = {y in X : (x, y) in
    U}}).  If $E$ is compact with respect to the topology on $X$
associated to $\mathcal{U}$ as in Section \ref{associated topology},
then $E$ is totally bounded.  Indeed, if $U \in \mathcal{U}$, then we
can cover $E$ by open subsets of $X$ of the form $\Int(U[x])$ with $x
\in E$, because of (\ref{x in Int(U[x])}), and then reduce to a finite
subcovering to get (\ref{E subseteq bigcup_{j = 1}^n U[x_j]}).

        Let $V \subseteq X \times X$ be given, and let us say that
$B \subseteq X$ is \emph{$V$-small}\index{small sets} if
\begin{equation}
\label{(x, y) in V}
        (x, y) \in V
\end{equation}
for every $x, y \in B$, which implies that every subset of $B$ is
$V$-small too.  Equivalently, $B$ is $V$-small when
\begin{equation}
\label{B times B subseteq V}
        B \times B \subseteq V.
\end{equation}
which is the same as saying that
\begin{equation}
\label{B subseteq V[x]}
        B \subseteq V[x]
\end{equation}
for every $x \in B$.  Suppose that $U \subseteq X \times X$ satisfies
\begin{equation}
\label{widetilde{U} * U subseteq V}
        \widetilde{U} * U \subseteq V,
\end{equation}
using the notation in (\ref{widetilde{U} = {(x, y) : (y, x) in U}})
and (\ref{U * V = ...}).  If
\begin{equation}
\label{B subseteq U[w]}
        B \subseteq U[w]
\end{equation}
for some $w \in X$, then it is easy to see that $B$ is $V$-small in
$X$.  Using this, one can check that $E \subseteq X$ is totally
bounded if and only if for each $V \in \mathcal{U}$ there are finitely
many $V$-small subsets of $X$ whose union contains $E$.  More
precisely, the ``if'' part of this statement follows from (\ref{B
  subseteq V[x]}).  In the other direction, if $V \in \mathcal{U}$ is
given, then there is a $U \in \mathcal{U}$ that satisfies
(\ref{widetilde{U} * U subseteq V}), by the definition of a
uniformity.  With this choice of $U$, (\ref{E subseteq bigcup_{j =
    1}^n U[x_j]}) implies that $E$ can be covered by finitely many
$V$-small subsets of $X$.

        Note that $V$, $\widetilde{V}$, and $V \cap \widetilde{V}$
determine the same collection of small subsets of $X$ for any $V
\subseteq X \times X$.  If $V_1, V_2 \subseteq X \times X$, then $B
\subseteq X$ is small with respect to $V_1 \cap V_2$ if and only if
$B$ is small with respect to both $V_1$ and $V_2$.  In particular, if
$B_1 \subseteq X$ is $V_1$-small, and $B_2 \subseteq X$ is
$V_2$-small, then $B_1 \cap B_2$ is small with respect to $V_1 \cap
V_2$.  If $E \subseteq X$ can be covered by finitely many sets that
are $V_1$ small, and by finitely many sets that are $V_2$ small, then
it follows that $E$ can be covered by finitely many sets that are
small with respect to $V_1 \cap V_2$, by taking intersections of the
various sets that are small with respect to $V_1$ and $V_2$.

        Let $(Y, \mathcal{V})$ be another uniform space, and suppose that
$f : X \to Y$ is uniformly continuous.  If $E \subseteq X$ is totally
bounded, then one can check that $f(E)$ is totally bounded in $Y$.

\section{Induced uniform structures}
\label{induced uniform structures}
\setcounter{equation}{0}

        Let $(Y, \mathcal{V})$ be a uniform space, and let $X$ be a
subset of $Y$.  Consider
\begin{equation}
\label{mathcal{U} = ...}
 \quad \mathcal{U} = \{U \subseteq X \times X : \hbox{ there is a } 
   V \in \mathcal{V} \hbox{ such that } U = V \cap (X \times X)\}.
\end{equation}
One can check that this is a uniformity on $X$, which may be described
as the uniformity induced on $X$ by $\mathcal{V}$ on
$Y$.\index{induced uniformities}\index{uniformities!induced} By
construction, the natural inclusion mapping from $X$ into $Y$ is
uniformly continuous with respect to $\mathcal{U}$ and $\mathcal{V}$.
If $\mathcal{B}(\mathcal{V})$ is a base for $\mathcal{V}$, then it is
easy to see that
\begin{equation}
\label{mathcal{B}(mathcal{U}) = ...}
        \mathcal{B}(\mathcal{U}) = \{V \cap (X \times X) :
                                      V \in \mathcal{B}(\mathcal{V})\}
\end{equation}
is a base for $\mathcal{U}$.

        Note that the restriction of a semimetric $d(\cdot, \cdot)$ on $Y$
to elements of $X$ defines a semimetric on $X$.  If $\mathcal{V}$ is
the uniformity determined on $Y$ by $d(\cdot, \cdot)$ as in Section
\ref{uniform structures}, then the induced uniformity on $X \subseteq
Y$ is the same as the uniformity determined on $X$ by the restriction
of $d(\cdot, \cdot)$ to elements of $X$.  One way to look at this is
to use bases for these uniformities consisting of sets of the form
(\ref{U(r) = U_d(r) = {(x, y) in X times X : d(x, y) < r}}) associated
to $d(\cdot, \cdot)$.  If $d(\cdot, \cdot)$ is a semi-ultrametric on $Y$,
then the restriction of $d(\cdot, \cdot)$ to elements of $X$ defines
a semi-ultrametric on $X$ too.

        Let $\mathcal{V}$ be any uniformity on $Y$ again, and let
$\mathcal{U}$ be the induced uniformity on $X \subseteq Y$.  This
leads to associated topologies on $X$ and $Y$, as in Section \ref{associated
topology}.  The associated topology on $Y$ also leads to an induced
topology on $X$ in the usual way, where $A \subseteq X$ is an open set
if there is an open set $B \subseteq Y$ such that
\begin{equation}
\label{A = B cap X}
        A = B \cap X.
\end{equation}
In this case, it is easy to see that $A$ is an open set with respect
to the topology on $X$ associated to $\mathcal{U}$, just by unwinding
the definitions.  The converse is a bit more complicated, as in the
context of semimetric spaces.

        Suppose that $A \subseteq X$ is an open set with respect to the
topology associated to $\mathcal{U}$.  Thus for each $x \in A$ there
is a $U_x \in \mathcal{U}$ such that
\begin{equation}
\label{U_x[x] subseteq A}
        U_x[x] \subseteq A,
\end{equation}
as in (\ref{U[x] subseteq A}), and where $U_x[x]$ is as in (\ref{U[x]
  = U[{x}] = {y in X : (x, y) in U}}).  By definition of
$\mathcal{U}$, for each $x \in A$ there is a $V_x \in \mathcal{V}$
such that
\begin{equation}
\label{U_x = V_x cap (X times X)}
        U_x = V_x \cap (X \times X).
\end{equation}
This implies that
\begin{equation}
\label{U_x[x] = V_x[x] cap X}
        U_x[x] = V_x[x] \cap X
\end{equation}
for each $x \in A$, where $V_x[x]$ is also defined as in (\ref{U[x] =
  U[{x}] = {y in X : (x, y) in U}}), but in $Y$ instead of $X$.  Let
$\Int(V_x[x])$ be the interior of $V_x[x]$ in $Y$ with respect to the
topology associated to $\mathcal{V}$ for each $x \in A$.  Put
\begin{equation}
\label{B = bigcup_{x in A} Int(V_x[x])}
        B = \bigcup_{x \in A} \Int(V_x[x]),
\end{equation}
which is also an open set in $Y$ with respect to the topology
associated to $\mathcal{V}$.  Observe that $A \subseteq B$, because $x
\in \Int(V_x[x])$ for each $x \in A$, as in (\ref{x in Int(U[x])}).
It follows that $A \subseteq B \cap X$, since $A \subseteq X$ by
hypothesis.  To get the opposite inclusion, one can use (\ref{U_x[x]
  subseteq A}), (\ref{U_x[x] = V_x[x] cap X}), and the fact that
$V_x[x] \subseteq \Int(V_x[x])$ for every $x \in A$ automatically.

\section{Induced uniform structures, continued}
\label{induced uniform structures, continued}
\setcounter{equation}{0}

        Let $X$ and $Y$ be sets, and let $f$ be a mapping from $X$ into $Y$.
This leads to a mapping $f_2$ from $X \times X$ into $Y \times Y$, as
in (\ref{f_2(x, x') = (f(x), f(x'))}).  If $U, U' \subseteq X \times
X$, then
\begin{equation}
\label{f_2(widetilde{U}) = widetilde{f_2(U)}}
        f_2(\widetilde{U}) = \widetilde{f_2(U)}
\end{equation}
and
\begin{equation}
\label{f_2(U * U') subseteq f_2(U) * f_2(U')}
        f_2(U * U') \subseteq f_2(U) * f_2(U'),
\end{equation}
where $\widetilde{U}$ and $U * U'$ are as in (\ref{widetilde{U} = {(x,
    y) : (y, x) in U}}) and (\ref{U * V = ...}).  If $f$ is injective,
then equality holds in (\ref{f_2(U * U') subseteq f_2(U) * f_2(U')}).
Similarly, if $V, V' \subseteq Y \times Y$, then
\begin{equation}
\label{f_2^{-1}(widetilde{V}) = widetilde{f_2^{-1}(V)}}
        f_2^{-1}(\widetilde{V}) = \widetilde{f_2^{-1}(V)}
\end{equation}
and
\begin{equation}
\label{f_2^{-1}(V) * f_2^{-1}(V') subseteq f_2^{-1}(V * V')}
        f_2^{-1}(V) * f_2^{-1}(V') \subseteq f_2^{-1}(V * V').
\end{equation}
If $f$ is surjective, then equality holds in (\ref{f_2^{-1}(V) *
  f_2^{-1}(V') subseteq f_2^{-1}(V * V')}).  Of course,
\begin{equation}
\label{f_2(Delta_X) subseteq Delta_Y}
        f_2(\Delta_X) \subseteq \Delta_Y,
\end{equation}
where $\Delta_X$, $\Delta_Y$ are as in (\ref{Delta = Delta_X = {(x, x)
    : x in X}}).

        If $\mathcal{V}$ is a uniformity on $Y$, then one can check that
\begin{equation}
\label{mathcal{B} = {f_2^{-1}(V) : V in mathcal{V}}}
        \mathcal{B} = \{f_2^{-1}(V) : V \in \mathcal{V}\}
\end{equation}
is a base for a uniformity $\mathcal{U}$ on $X$.  Equivalently,
$\mathcal{U}$ consists of the $U \subseteq X \times X$ for which
there is a $V \in \mathcal{V}$ such that
\begin{equation}
\label{f_2^{-1}(V) subseteq U}
        f_2^{-1}(V) \subseteq U.
\end{equation}
By construction, $f$ is uniformly continuous as a mapping from $X$
into $Y$ with respect to $\mathcal{U}$ and $\mathcal{V}$,
respectively.  If $f$ is injective, then (\ref{mathcal{B} =
  {f_2^{-1}(V) : V in mathcal{V}}}) is already a uniformity on $X$.
In particular, if $X \subseteq Y$ and $f$ is the natural inclusion
mapping from $X$ into $Y$, then (\ref{mathcal{B} = {f_2^{-1}(V) : V in
    mathcal{V}}}) is the same as the uniformity induced on $X$ by
$\mathcal{V}$ on $Y$ as in (\ref{mathcal{U} = ...}).

        Let $f$ be any mapping from a set $X$ into a set $Y$ again,
and let $\mathcal{B}_Y$ be a base for a uniformity $\mathcal{V}$
on $Y$.  In this case,
\begin{equation}
\label{mathcal{B}_X = {f_2^{-1}(V) : V in mathcal{B}_Y}}
        \mathcal{B}_X = \{f_2^{-1}(V) : V \in \mathcal{B}_Y\}
\end{equation}
is a base for the same uniformity $\mathcal{U}$ on $X$ as in the
previous paragraph.  There is an analogous statement for sub-bases,
which are defined in the next section.  Note that a mapping from
another uniform space into $X$ is uniformly continuous with respect to
$\mathcal{U}$ on $X$ if and only if the composition of this mapping
with $f$ is uniformly continuous as a mapping into $Y$, with respect
to the given uniformity $\mathcal{V}$ on $Y$.  The ``only if'' part of
this statement follows from the fact that compositions of uniformly
continuous mappings are uniformly continuous, as in Section
\ref{uniform continuity}, while the ``if'' part uses the way that
$\mathcal{U}$ is defined in terms of $\mathcal{V}$ here.

        If $d_Y(\cdot, \cdot)$ is a semimetric on $Y$, then
it is easy to see that
\begin{equation}
\label{d_X(x, x') = d_Y(f(x), f(x'))}
        d_X(x, x') = d_Y(f(x), f(x'))
\end{equation}
defines a semimetric on $X$.  Put
\begin{equation}
\label{U_{d_Y}(r) = {(y, y') in Y times Y : d_Y(y, y') < r}}
        U_{d_Y}(r) = \{(y, y') \in Y \times Y : d_Y(y, y') < r\}
\end{equation}
and
\begin{equation}
\label{U_{d_X}(r) = {(x, x') in X times X : d_X(x, x') < r}}
        U_{d_X}(r) = \{(x, x') \in X \times X : d_X(x, x') < r\}
\end{equation}
for each $r > 0$, as in (\ref{U(r) = U_d(r) = {(x, y) in X times X :
    d(x, y) < r}}).  If $f_2 : X \times Y \to Y \times Y$ is as in
(\ref{f_2(x, x') = (f(x), f(x'))}) again, then
\begin{equation}
\label{f_2^{-1}(U_{d_Y}(r)) = U_{d_X}(r)}
        f_2^{-1}(U_{d_Y}(r)) = U_{d_X}(r)
\end{equation}
for every $r > 0$.  Let $\mathcal{B}_{d_Y}$ be the collection of
subsets of $Y \times Y$ of the form (\ref{U_{d_Y}(r) = {(y, y') in Y
    times Y : d_Y(y, y') < r}}) for some $r > 0$, and let
$\mathcal{B}_{d_X}$ be the collection of subsets of $X \times X$ of
the form (\ref{U_{d_X}(r) = {(x, x') in X times X : d_X(x, x') < r}})
for some $r > 0$, as in Section \ref{uniform structures}.  These are
bases for the uniformities $\mathcal{U}_{d_Y}$ and $\mathcal{U}_{d_X}$
on $X$ and $Y$ associated to $d_Y(\cdot, \cdot)$ and $d_X(\cdot,
\cdot)$, respectively.  In this situation, $\mathcal{U}_{d_X}$ is the
same as the uniformity induced on $X$ by $\mathcal{U}_{d_Y}$ on $Y$
and the mapping $f$ as before, because of (\ref{f_2^{-1}(U_{d_Y}(r)) =
  U_{d_X}(r)}).  Note that (\ref{d_X(x, x') = d_Y(f(x), f(x'))}) is a
semi-ultrametric on $X$ when $d_Y(\cdot, \cdot)$ is a semi-ultrametric
on $Y$.

\section{Sub-bases}
\label{sub-bases}
\setcounter{equation}{0}

        Let $X$ be a set, let $\mathcal{U}$ be a uniformity on $X$,
and let $\mathcal{B}_0$ be a subcollection of $\mathcal{U}$.  Consider
the collection $\mathcal{B}$ of subsets of $X \times X$ that can be
expressed as the intersection of finitely many elements of
$\mathcal{B}_0$.  If $\mathcal{B}$ is a base for $\mathcal{U}$, as in
Section \ref{uniform structures}, then $\mathcal{B}_0$ is said to be a
\emph{sub-base}\index{sub-base for a uniformity} for $\mathcal{U}$.
Note that $\mathcal{B}_0 \ne \emptyset$ in this case, because
$\mathcal{U} \ne \emptyset$, and hence $\mathcal{B} \ne \emptyset$.
If $U \in \mathcal{B}_0$, then (\ref{Delta subseteq U}) holds, since
$U \in \mathcal{U}$.  There should also be a $V_0 \in \mathcal{B}$
that satisifes (\ref{V_0 subseteq widetilde{U}}), because $U \in
\mathcal{B}$.  Similarly, there should be a $V \in \mathcal{B}$ that
satisfies (\ref{V * V subseteq U}).

        Now let $\mathcal{B}_0$ be any nonempty collection of subsets of
$X \times X$,  and let $\mathcal{B}$ be the collection of subsets of
$X \times X$ that can be expressed as the intersection of finitely many
elements of $\mathcal{B}_0$, as before.  If $\mathcal{B}_0$ satisfies
the three conditions with respect to $\mathcal{B}$ mentioned at the
end of the preceding paragraph, then it is easy to see that
$\mathcal{B}$ has the analogous properties with respect to itself.  Of
course, the intersection of any two elements of $\mathcal{B}$ is
automatically an element of $\mathcal{B}$, by construction.  As in
Section \ref{uniform structures}, it follows that
(\ref{mathcal{U}(mathcal{B}) = ...}) is a uniformity on $X$, and that
$\mathcal{B}$ is a base for this uniformity.  Thus $\mathcal{B}_0$
is a sub-base for this uniformity.

        Let $I$ be a nonempty set, and suppose that $\mathcal{B}_i$
is a sub-base for a uniformity on $X$ for each $i \in I$.  Under these
conditions, one can check that
\begin{equation}
\label{bigcup_{i in I} mathcal{B}_i}
        \bigcup_{i \in I} \mathcal{B}_i
\end{equation}
satisfies the requirements of a sub-base for a uniformity on $X$
described in the previous paragraph.  Note that (\ref{bigcup_{i in I}
  mathcal{B}_i}) is not necessarily a base for a uniformity on $X$,
even when $\mathcal{B}_i$ is a base for a uniformity on $X$ for each
$i \in I$.  As a basic scenario, let $(Y_i, \mathcal{V}_i)$ be a
uniform space for each $i \in I$, and let $f_i$ be a mapping from $X$
into $Y_i$.  This leads to a mapping $f_{i, 2}$ from $X \times X$ into
$Y_i \times Y_i$ for each $i \in I$, as in (\ref{f_2(x, x') = (f(x),
  f(x'))}).  We also get a base
\begin{equation}
\label{mathcal{B}_i = {f_{i, 2}^{-1}(V_i) : V_i in mathcal{V}_i}}
        \mathcal{B}_i = \{f_{i, 2}^{-1}(V_i) : V_i \in \mathcal{V}_i\}
\end{equation}
for a uniformity on $X$ for each $i \in I$, as in (\ref{mathcal{B} =
  {f_2^{-1}(V) : V in mathcal{V}}}).  Thus (\ref{bigcup_{i in I}
  mathcal{B}_i}) is a sub-base for a uniformity $\mathcal{U}$ on $X$,
as before.  By construction, for each $i \in I$, $f_i$ is uniformly
continuous as a mapping from $X$ into $Y_i$ with respect to the
uniformities $\mathcal{U}$ and $\mathcal{V}_i$, respectively.

        Let $(X, \mathcal{U})$ and $(Y, \mathcal{V})$ be uniform
spaces, and let $f$ be a mapping from $X$ into $Y$.  In order to check
that $f$ is uniformly continuous, it suffices to verify that
(\ref{f_2^{-1}(V) in mathcal{U}}) holds for every $V$ in a sub-base
for $\mathcal{V}$, as mentioned in Section \ref{uniform continuity}.
In particular, suppose that there is a sub-base for $\mathcal{V}$
which is the union of a family of sub-bases for other uniformities on
$Y$, as in the preceding paragraph.  If $f$ is uniformly continuous as
a mapping into $Y$ with respect to each of these other uniformities on
$Y$, then it follows that $f$ is uniformly continuous as a mapping
into $Y$ with respect to $\mathcal{V}$.  The converse is trivial,
since the uniform continuity of $f$ as a mapping into $Y$ with respect
to $\mathcal{V}$ automatically implies that $f$ is uniformly
continuous as a mapping into $Y$ with respect to any uniformity
contained in $\mathcal{V}$.

\section{Cartesian products}
\label{cartesian products}
\setcounter{equation}{0}

        Let $I$ be a nonempty set, and suppose that $(Y_i, \mathcal{V}_i)$
is a uniform space for each $i \in I$.  Also let
\begin{equation}
\label{Y = prod_{i in I} Y_i}
        Y = \prod_{i \in I} Y_i
\end{equation}
be the Cartesian product of the $Y_i$'s, and let $p_j$ be the standard
coordinate projection from $Y$ onto $Y_j$ for each $j \in I$.  This
leads to a mapping $p_{j, 2}$ from $Y \times Y$ onto $Y_j \times Y_j$
for each $j \in I$, as in (\ref{f_2(x, x') = (f(x), f(x'))}).  Put
\begin{equation}
\label{mathcal{B}_j = {p_{j, 2}^{-1}(V_j) : V_j in mathcal{V}_j}}
        \mathcal{B}_j = \{p_{j, 2}^{-1}(V_j) : V_j \in \mathcal{V}_j\}
\end{equation}
for each $j \in I$, which is a base for a uniformity on $Y$, as in
(\ref{mathcal{B} = {f_2^{-1}(V) : V in mathcal{V}}}).  Thus
\begin{equation}
\label{bigcup_{j in I} mathcal{B}_j}
        \bigcup_{j \in I} \mathcal{B}_j
\end{equation}
is a sub-base for a uniformity $\mathcal{V}$ on $Y$, as in the
previous section.  This uniformity $\mathcal{V}$ is the \emph{product
  uniformity}\index{product uniformities}\index{uniformities!product}
on $Y$ associated to the given uniformities on the $Y_j$'s.  By
construction, $p_j$ is uniformly continuous as a mapping from $Y$ onto
$Y_j$ for each $j \in I$, with respect to the product uniformity on
$Y$ and $\mathcal{V}_j$ on $Y_j$.

        Let $(X, \mathcal{U})$ be another uniform space, and let $f_i$
be a mapping from $X$ into $Y_i$ for each $i \in I$.  This leads to a
mapping $f$ from $X$ into $Y$, with
\begin{equation}
\label{f_i = p_i circ f}
        f_i = p_i \circ f
\end{equation}
for each $i \in I$.  If $f_i$ is uniformly continuous for each $i$,
then one can check that $f$ is uniformly continuous with respect to
the product uniformity on $Y$.  More precisely, it suffices to verify
that (\ref{f_2^{-1}(V) in mathcal{U}}) holds for every $V$ in the
sub-base (\ref{bigcup_{j in I} mathcal{B}_j}) for the product
uniformity, which reduces to the uniform continuity of the $f_j$'s.
Conversely, if $f$ is uniformly continuous, then $f_i$ is uniformly
continuous for each $i \in I$, because it is the composition of
uniformly continuus mappings, as in (\ref{f_i = p_i circ f}).

        Of course,
\begin{equation}
\label{Y times Y = (prod_{i in I} Y_i) times (prod_{i in I} Y_i)}
 Y \times Y = \Big(\prod_{i \in I} Y_i\Big) \times \Big(\prod_{i \in I} Y_i\Big)
\end{equation}
can be identified with
\begin{equation}
\label{prod_{i in I} (Y_i times Y_i)}
        \prod_{i \in I} (Y_i \times Y_i)
\end{equation}
in a natural way.  Let $j \in I$ be given, and observe that $p_{j, 2}$
corresponds to the standard coordinate projection from (\ref{prod_{i
    in I} (Y_i times Y_i)}) onto $Y_j \times Y_j$.  If $V_j \subseteq
Y_j \times Y_j$, then
\begin{equation}
\label{p_{j, 2}^{-1}(V_j) subseteq Y times Y}
        p_{j, 2}^{-1}(V_j) \subseteq Y \times Y
\end{equation}
corresponds to the subset of (\ref{prod_{i in I} (Y_i times Y_i)})
which is the Cartesian product of $V_j$ with $Y_i \times Y_i$ for each
$i \in I$ different from $j$.  If $y \in Y$, then $V_j[p_j(y)]$ can be
defined as a subset of $Y_j$ as in (\ref{U[x] = U[{x}] = {y in X : (x,
    y) in U}}), and similarly $(p_{j, 2}^{-1}(V_j))[y]$ can be defined
as a subset of $Y$.  It is easy to see that
\begin{equation}
\label{(p_{j, 2}^{-1}(V_j))[y] = p_j^{-1}(V_j[p_j(y)])}
        (p_{j, 2}^{-1}(V_j))[y] = p_j^{-1}(V_j[p_j(y)]),
\end{equation}
which is the same as the Cartesian product of $V_j[p_j(y)]$ with $Y_i$
for each $i \in I$ with $i \ne j$.  

        Each $Y_i$ has a topology associated to the given uniformity
$\mathcal{V}_i$, as in Section \ref{associated topology}, and there is
a topology on $Y$ associated to the product uniformity as well.  There
is also the product topology on $Y$ corresponding to the topologies on
the $Y_i$'s just mentioned.  It is easy to see that an open subset of
$Y$ with respect to the product topology is an open set with respect
to the topology associated to the product uniformity, directly from
the definitions.  To show that an open subset of $Y$ with respect to
the topology associated to the product uniformity is an open set with
respect to the product topology, one can use (\ref{x in Int(U[x])}).
One can look at this in terms of local bases and sub-bases for the
relevant topologies, as on p182-3 of \cite{jk}.

\section{Cartesian products, continued}
\label{cartesian products, continued}
\setcounter{equation}{0}

        Let $(X, \mathcal{U})$ be a uniform space, and let us consider
the corresponding product uniformity on $X \times X$, as in the previous
section.  If $U, V \in \mathcal{U}$, then
\begin{equation}
\label{U_1 = ...}
 U_1 = \{((x_1, x_2), (x_1', x_2')) \in (X \times X) \times (X \times X) :
                                                (x_1, x_1') \in U\}
\end{equation}
and
\begin{equation}
\label{V_2 = ...}
 V_2 = \{((x_1, x_2), (x_1', x_2')) \in (X \times X) \times (X \times X) :
                                                 (x_2, x_2') \in V\}
\end{equation}
are sub-basic elements of the product uniformity on $X \times X$.
More precisely, $U_1$ and $V_2$ correspond to the first and second
coordinate projections from $X \times X$ onto $X$, as in
(\ref{mathcal{B}_j = {p_{j, 2}^{-1}(V_j) : V_j in mathcal{V}_j}}) and
(\ref{p_{j, 2}^{-1}(V_j) subseteq Y times Y}).  Thus
\begin{eqnarray}
\label{U_1 cap V_2 = ...}
\lefteqn{\qquad U_1 \cap V_2 =}  \\
 & & \{((x_1, x_2), (x_1', x_2')) \in (X \times X) \times (X \times X) :
                      (x_1, x_1') \in U, \ (x_2, x_2') \in V\}  \nonumber
\end{eqnarray}
is an element of the product uniformity on $X \times X$.  The
collection of these sets forms a base for the product uniformity on $X
\times X$.

        If $(x_1, x_2) \in X \times X$, then $U_1[(x_1, x_2)]$ and
$V_2[(x_1, x_2)]$ can be defined as subsets of $X \times X$ as in
(\ref{U[x] = U[{x}] = {y in X : (x, y) in U}}), and we have that
\begin{equation}
\label{U_1[(x_1, x_2)] = (U[x_1]) times X}
        U_1[(x_1, x_2)] = (U[x_1]) \times X
\end{equation}
and
\begin{equation}
\label{V_2[(x_1, x_2)] = X times (V[x_2])}
        V_2[(x_1, x_2)] = X \times (V[x_2]).
\end{equation}
Similarly,
\begin{equation}
\label{(U_1 cap V_2)[(x_1, x_2)] = (U[x_1]) times (V[x_2])}
        (U_1 \cap V_2)[(x_1, x_2)] = (U[x_1]) \times (V[x_2]).
\end{equation}
If $A \subseteq X \times X$, then $(U_1 \cap V_2)[A]$ can be defined
as a subset of $X \times X$ as in (\ref{U[A] = {y in X : there is an x
    in A such that (x, y) in U}}), which is equivalent to (\ref{U[A] =
  bigcup_{x in A} U[x]}).  In this situation, (\ref{U[A] = bigcup_{x
    in A} U[x]}) reduces to
\begin{equation}
\label{(U_1 cap V_2)[A] = bigcup_{(x_1, x_2) in A} (U[x_1]) times (V_2[x_2])}
 (U_1 \cap V_2)[A] = \bigcup_{(x_1, x_2) \in A} (U[x_1]) \times (V_2[x_2]),
\end{equation}
because of (\ref{(U_1 cap V_2)[(x_1, x_2)] = (U[x_1]) times (V[x_2])}).

        If $A = \Delta$, the diagonal in $X \times X$, then
(\ref{(U_1 cap V_2)[A] = bigcup_{(x_1, x_2) in A} (U[x_1]) times (V_2[x_2])})
reduces to
\begin{equation}
\label{(U_1 cap V_2)[Delta] = bigcup_{x in X} (U_1[x]) times (V_2[x])}
        (U_1 \cap V_2)[\Delta] = \bigcup_{x \in X} (U_1[x]) \times (V_2[x]).
\end{equation}
It follows that
\begin{equation}
\label{(U_1 cap V_2)[Delta] = widetilde{U} * V}
        (U_1 \cap V_2)[\Delta] = \widetilde{U} * V,
\end{equation}
where $\widetilde{U}$ is as defined in (\ref{widetilde{U} = {(x, y) :
    (y, x) in U}}), and then $\widetilde{U} * V$ is defined as in
(\ref{U * V = ...}).  If $W$ is any element of $\mathcal{U}$, then
there are $U, V \in \mathcal{U}$ such that
\begin{equation}
\label{widetilde{U} * V subseteq W}
        \widetilde{U} * V \subseteq W.
\end{equation}
More precisely, there is a $V \in \mathcal{U}$ such that $V * V
\subseteq W$, as in (\ref{V * V subseteq U}), and one can take
$U = \widetilde{V}$.

\section{Compatible semimetrics}
\label{compatible semimetrics}
\setcounter{equation}{0}

        Let $(X, \mathcal{U})$ be a uniform space again, and let
$d(\cdot, \cdot)$ be a semimetric on $X$.  Also let $\mathcal{U}_d$
be the uniformity on $X$ associated to $d(\cdot, \cdot)$ as in
Section \ref{uniform structures}.  Let us say that $d(\cdot, \cdot)$
is \emph{compatible}\index{compatible semimetrics} with $\mathcal{U}$
if
\begin{equation}
\label{mathcal{U}_d subseteq mathcal{U}}
        \mathcal{U}_d \subseteq \mathcal{U}.
\end{equation}
Equivalently, this means that the identity mapping on $X$ is uniformly
continuous as a mapping from $X$ equipped with $\mathcal{U}$ into $X$
equipped with $\mathcal{U}_d$.  If $U_d(r)$ is defined as in
(\ref{U(r) = U_d(r) = {(x, y) in X times X : d(x, y) < r}}) for each
$r > 0$, then (\ref{mathcal{U}_d subseteq mathcal{U}}) is the same as
saying that
\begin{equation}
\label{U_d(r) in mathcal{U}}
        U_d(r) \in \mathcal{U}
\end{equation}
for every $r > 0$.

        As in Theorem 11 on p183 of \cite{jk}, $d(\cdot, \cdot)$ is
compatible with $\mathcal{U}$ in the sense described in the preceding
paragraph if and only if $d(\cdot, \cdot)$ is uniformly continuous on
$X \times X$ with respect to the product uniformity corresponding to
$\mathcal{U}$ on $X$ as in the previous sections.  Here we use the
standard uniform structure on the real line, corresponding to the
standard Euclidean metric, for the range of $d(\cdot, \cdot)$ on $X
\times X$.  To prove the ``if'' part, suppose that $d(\cdot, \cdot)$
is uniformly continuous on $X \times X$, and let $r > 0$.  Under these
conditions, there are $U, V \in \mathcal{U}$ so that if $U_1$ and
$V_2$ are as in (\ref{U_1[(x_1, x_2)] = (U[x_1]) times X}) and
(\ref{V_2[(x_1, x_2)] = X times (V[x_2])}), respectively, then
\begin{equation}
\label{|d(x_1, x_2) - d(x_1', x_2')| < r}
        |d(x_1, x_2) - d(x_1', x_2')| < r
\end{equation}
for every $((x_1, x_2), (x_1', x_2')) \in U_1 \cap V_2$.  This uses
the fact that the sets $U_1 \cap V_2$ with $U, V \in \mathcal{U}$ form
a base for the product uniformity on $X \times X$, as in the previous
section.  If we restrict our attention to $x_1 = x_2$ in (\ref{|d(x_1,
  x_2) - d(x_1', x_2')| < r}), then we get that
\begin{equation}
\label{d(x_1', x_2') < r}
        d(x_1', x_2') < r
\end{equation}
for every $(x_1', x_2') \in X \times X$ for which there is an $x \in X$
such that
\begin{equation}
\label{((x, x), (x_1', x_2')) in U_1 cap V_2}
        ((x, x), (x_1', x_2')) \in U_1 \cap V_2.
\end{equation}
Equivalently, this means that (\ref{d(x_1', x_2') < r}) holds when
$(x_1', x_2')$ is an element of (\ref{(U_1 cap V_2)[Delta] = bigcup_{x
    in X} (U_1[x]) times (V_2[x])}), so that
\begin{equation}
\label{(U_1 cap V_2)[Delta] subseteq U_d(r)}
        (U_1 \cap V_2)[\Delta] \subseteq U_d(r),
\end{equation}
where $U_d(r)$ is as in (\ref{U(r) = U_d(r) = {(x, y) in X times X :
    d(x, y) < r}}).  It follows that
\begin{equation}
\label{widetilde{U} * V subseteq U_d(r)}
        \widetilde{U} * V \subseteq U_d(r),
\end{equation}
using also (\ref{(U_1 cap V_2)[Delta] = widetilde{U} * V}).  Remember
that
\begin{equation}
\label{Delta subseteq U, V}
        \Delta \subseteq U, V,
\end{equation}
because $U, V \in \mathcal{U}$, and hence $\Delta \subseteq
\widetilde{U}$ as well.  Thus (\ref{widetilde{U} * V subseteq U_d(r)})
implies that
\begin{equation}
\label{widetilde{U}, V subseteq U_d(r)}
        \widetilde{U}, V \subseteq U_d(r),
\end{equation}
which corresponds to taking $x$ equal to $x_1'$ or $x_2'$ in
(\ref{((x, x), (x_1', x_2')) in U_1 cap V_2}).  Note that
(\ref{widetilde{U}, V subseteq U_d(r)}) also implies that $U \subseteq
U_d(r)$, because $U_d(r)$ is symmetric.  It follows that (\ref{U_d(r)
  in mathcal{U}}) holds, as desired, by the definition of a
uniformity.  More precisely, this shows that $d(\cdot, \cdot)$ is
compatible with $\mathcal{U}$ on $X$ when $d(\cdot, \cdot)$ is
uniformly continuous along the diagonal $\Delta$ in $X \times X$,
which corresponds to the restriction to $x_1 = x_2$ in (\ref{|d(x_1,
  x_2) - d(x_1', x_2')| < r}).  Of course, if $d(\cdot, \cdot)$ is
uniformly continuous on $X \times X$, then $d(\cdot, \cdot)$ is
uniformly continuous along any subset of $X \times X$.

        Before proving the converse, let us record some simple estimates.
Using the triangle inequality twice, we get that
\begin{equation}
\label{d(x_1, x_2) le d(x_1', x_2') + d(x_1, x_1') + d(x_2, x_2')}
        d(x_1, x_2) \le d(x_1', x_2') + d(x_1, x_1') + d(x_2, x_2')
\end{equation}
for every $x_1, x_2, x_1', x_2' \in X$, and hence
\begin{equation}
\label{d(x_1, x_2) - d(x_1', x_2') le d(x_1, x_1') + d(x_2, x_2')}
        d(x_1, x_2) - d(x_1', x_2') \le d(x_1, x_1') + d(x_2, x_2').
\end{equation}
Similarly,
\begin{equation}
\label{d(x_1', x_2') - d(x_1, x_2) le d(x_1, x_1') + d(x_2, x_2')}
        d(x_1', x_2') - d(x_1, x_2) \le d(x_1, x_1') + d(x_2, x_2').
\end{equation}
Combining (\ref{d(x_1, x_2) - d(x_1', x_2') le d(x_1, x_1') + d(x_2,
  x_2')}) and (\ref{d(x_1', x_2') - d(x_1, x_2) le d(x_1, x_1') +
  d(x_2, x_2')}), we get that
\begin{equation}
\label{|d(x_1, x_2) - d(x_1', x_2')| le d(x_1, x_1') + d(x_2, x_2')}
        |d(x_1, x_2) - d(x_1', x_2')| \le d(x_1, x_1') + d(x_2, x_2').
\end{equation}
If $(x_1, x_1'), (x_2, x_2') \in U_d(r)$ for some $r > 0$, then it
follows that
\begin{equation}
\label{|d(x_1, x_2) - d(x_1', x_2')| < r + r = 2 r}
        |d(x_1, x_2) - d(x_1', x_2')| < r + r = 2 \, r,
\end{equation}
by the definition (\ref{U(r) = U_d(r) = {(x, y) in X times X : d(x, y)
    < r}}) of $U_d(r)$.

        Suppose now that $d(\cdot, \cdot)$ is compatible with
$\mathcal{U}$ on $X$, and let $r > 0$ be given.  Thus
(\ref{U_d(r) in mathcal{U}}) holds, and we can take $U = V = U_d(r)$
in (\ref{U_1 cap V_2 = ...}), to get an element of the product
uniformity on $X \times X$ associated to $\mathcal{U}$ on $X$.  In
this case, we have just seen that (\ref{|d(x_1, x_2) - d(x_1', x_2')|
  < r + r = 2 r}) holds when $((x_1, x_2), (x_1', x_2'))$ is an
element of (\ref{U_1 cap V_2 = ...}).  This implies that $d(\cdot,
\cdot)$ is uniformly continuous on $X \times X$ with respect to the
product uniformity associated to $\mathcal{U}$ on $X$, as desired.
In particular, $d(\cdot, \cdot)$ is uniformly continuous on $X \times X$
with respect to the product uniformity associated to the uniformity
$\mathcal{U}_d$ determined on $X$ by itself, as in Section \ref{uniform
structures}.

\section{Collections of semimetrics}
\label{collections of semimetrics}
\setcounter{equation}{0}

        Let $X$ be a set, and let $\mathcal{M}$ be a nonempty
collection of semimetrics on $X$.  Also let $U_d(r)$ be the subset of
$X \times X$ associated to $d \in \mathcal{M}$ and a positive real
number $r$ as in (\ref{U(r) = U_d(r) = {(x, y) in X times X : d(x, y)
    < r}}).  Thus
\begin{equation}
\label{mathcal{B}_d = {U_d(r) : r > 0}}
        \mathcal{B}_d = \{U_d(r) : r > 0\}
\end{equation}
is a base for a uniformity $\mathcal{U}_d$ on $X$ for each $d \in
\mathcal{M}$, as discussed in Section \ref{uniform structures}.  It
follows that
\begin{equation}
\label{bigcup_{d in mathcal{M}} mathcal{B}_d}
        \bigcup_{d \in \mathcal{M}} \mathcal{B}_d
\end{equation}
is a sub-base for a uniformity on $X$, as in Section \ref{sub-bases}.
Of course, each element of $\mathcal{M}$ is compatible with this
uniformity on $X$, by construction.

        Now let $I$ be nonempty set, let $Y_i$ be a set for each $i \in I$,
and let
\begin{equation}
\label{Y - prod_{i in I} Y_i}
        Y - \prod_{i \in I} Y_i
\end{equation}
be the corresponding Cartesian product, as in Section \ref{cartesian
  products}.  Suppose that $\mathcal{M}_j$ is a nonempty collection of
semimetrics on $Y_j$ for each $j \in I$, and let $\mathcal{V}_j$ be
the uniformity on $Y_j$ associated to $\mathcal{M}_j$ as in the
preceding paragraph.  If $j \in I$ and $d_j \in \mathcal{M}_j$, then
put
\begin{equation}
\label{widehat{d}_j(y, y') = d_j(p_j(y), p_j(y'))}
        \widehat{d}_j(y, y') = d_j(p_j(y), p_j(y'))
\end{equation}
for every $y, y' \in Y$, where $p_j$ is the standard coordinate
projection from $Y$ onto $Y_j$, as before.  Note that
(\ref{widehat{d}_j(y, y') = d_j(p_j(y), p_j(y'))}) defines a
semimetric on $Y$, as in (\ref{d_X(x, x') = d_Y(f(x), f(x'))}).  Let
$\widehat{\mathcal{M}}_j$ be the collection of semimetrics on $Y$ of
the form (\ref{widehat{d}_j(y, y') = d_j(p_j(y), p_j(y'))}) with
$d_j(\cdot, \cdot) \in \mathcal{M}_j$, and put
\begin{equation}
\label{widehat{mathcal{M}} = bigcup_{j in I} widehat{mathcal{M}}_j}
        \widehat{\mathcal{M}} = \bigcup_{j \in I} \widehat{\mathcal{M}}_j.
\end{equation}
This is a nonempty collection of semimetrics on $Y$, which determines
a uniform structure on $Y$ as in the previous paragraph.  In this
situation, the uniform structure on $Y$ associated to $\mathcal{M}$ is
the same as the product uniformity on $Y$ corresponding to the
uniformity $\mathcal{V}_j$ on $Y_j$ associated to $\mathcal{M}_j$ for
each $j \in I$.

        Suppose that $d_1, \ldots, d_n$ are finitely many semimetrics
on a set $X$.  Under these conditions, one can check that
\begin{equation}
\label{d(x, y) = max_{1 le j le n} d_j(x, y)}
        d(x, y) = \max_{1 \le j \le n} d_j(x, y)
\end{equation}
defines a semimetric on $X$ as well.  Put
\begin{equation}
\label{U_{d_j}(r) = {(x, y) in X times X : d_j(x, y) < r}}
        U_{d_j}(r) = \{(x, y) \in X \times X : d_j(x, y) < r\}
\end{equation}
for each $j = 1, \ldots, n$ and $r > 0$, and let $U_d(r)$ be the
analogous set corresponding to $d$, as in (\ref{U(r) = U_d(r) = {(x,
    y) in X times X : d(x, y) < r}}).  Observe that
\begin{equation}
\label{U_d(r) = bigcap_{j = 1}^n U_{d_j}(r)}
        U_d(r) = \bigcap_{j = 1}^n U_{d_j}(r)
\end{equation}
for every $r > 0$.  Let $\mathcal{B}_d$ be the collection of subsets
of $X \times X$ of the form $U_d(r)$ for some $r > 0$, as in Section
\ref{uniform structures}, and let $\mathcal{B}_{d_j}$ be the analogous
collection associated to $d_j$ for each $j = 1, \ldots, n$.  Thus
$\mathcal{B}_d$ is a base for the uniformity $\mathcal{U}_d$ on $X$
associated to $d$, as before, and similarly $\mathcal{B}_{d_j}$ is a
base for the uniformity $\mathcal{U}_{d_j}$ on $X$ associated to $d_j$
for each $j = 1, \ldots, n$.  We also have that
\begin{equation}
\label{bigcup_{j = 1}^n mathcal{B}_{d_j}}
        \bigcup_{j = 1}^n \mathcal{B}_{d_j}
\end{equation}
is a sub-base for a uniformity on $X$, as in (\ref{bigcup_{d in
    mathcal{M}} mathcal{B}_d}).  One can check that this is actually a
sub-base for $\mathcal{U}_d$, using (\ref{U_d(r) = bigcap_{j = 1}^n
  U_{d_j}(r)}).  If $d_1, \ldots, d_n$ are semi-ultrametrics on $X$,
then (\ref{d(x, y) = max_{1 le j le n} d_j(x, y)}) is a
semi-ultrametric on $X$ too.

        If $(X, \mathcal{U})$ is any uniform space, then it is well known
that there is a collection $\mathcal{M}$ of semimetrics on $X$ for which
$\mathcal{U}$ is the corresponding uniformity.

\section{Sequences of semimetrics}
\label{sequences of semimetrics}
\setcounter{equation}{0}

        Let $d(\cdot, \cdot)$ be a semimetric on a set $X$, and let
$t$ be a positive real number.  Under these conditions, it is easy
to see that
\begin{equation}
\label{d_t(x, y) = min(d(x, y), t)}
        d_t(x, y) = \min(d(x, y), t)
\end{equation}
also defines a semimetric on $X$.  Put
\begin{equation}
\label{U_{d_t}(r) = {(x, y) in X times X : d_t(x, y) < r}}
        U_{d_t}(r) = \{(x, y) \in X \times X : d_t(x, y) < r\}
\end{equation}
for each $r > 0$, and let $U_d(r)$ be the analogous subset of $X
\times X$ corresponding to $d(\cdot, \cdot)$, as in (\ref{U(r) =
  U_d(r) = {(x, y) in X times X : d(x, y) < r}}).  Observe that
\begin{equation}
\label{U_{d_t}(r) = U_d(r)}
        U_{d_t}(r) = U_d(r)
\end{equation}
when $r \le t$, and
\begin{equation}
\label{U_{d_t}(r) = X times X}
        U_{d_t}(r) = X \times X
\end{equation}
when $r > t$.  It follows that the uniformity on $X$ associated to
$d_t$ as in Section \ref{uniform structures} is the same as the
uniformity associated to $d$ for every $t > 0$.  In particular, the
corresponding topologies are the same.  Note that (\ref{d_t(x, y) =
  min(d(x, y), t)}) is a semi-ultrametric on $X$ when $d(\cdot,
\cdot)$ is a semi-ultrametric on $X$.

        Now let $d_1, d_2, d_3, \ldots$ be an infinite sequence of
semimetrics on $X$, and put
\begin{equation}
\label{d_j'(x, y) = min(d_j(x, y), 1/j)}
        d_j'(x, y) = \min(d_j(x, y), 1/j)
\end{equation}
for every positive integer $j$ and $x, y \in X$.  Thus $d_j'$ is also
a semimetric on $X$ for each $j$, as in the previous paragraph, and
$d_j'$ determines the same uniformity on $X$ as $d_j$.  Put
\begin{equation}
\label{d'(x, y) = max_{j ge 1} d_j'(x, y)}
        d'(x, y) = \max_{j \ge 1} d_j'(x, y)
\end{equation}
for each $x, y \in X$.  Of course, (\ref{d'(x, y) = max_{j ge 1}
  d_j'(x, y)}) is equal to $0$ when $d_j'(x, y) = 0$ for every $j$.
Otherwise, if $d_l'(x, y) > 0$ for some $l$, then (\ref{d'(x, y) =
  max_{j ge 1} d_j'(x, y)}) reduces to the maximum over finitely many
$j$, because (\ref{d_j'(x, y) = min(d_j(x, y), 1/j)}) is automatically
less than or equal to $1/j$.  This implies that the maximum in
(\ref{d'(x, y) = max_{j ge 1} d_j'(x, y)}) is always attained, and one
can check that (\ref{d'(x, y) = max_{j ge 1} d_j'(x, y)}) defines a
semimetric on $X$.  If $d_j$ is a semi-ultrametric on $X$ for each
$j$, then $d_j'$ is a semi-ultrametric on $X$ for each $j$, and
(\ref{d'(x, y) = max_{j ge 1} d_j'(x, y)}) is a semi-ultrametric on
$X$ as well.

        As usual, for each $r > 0$, we put
\begin{eqnarray}
\label{U_{d_j}(r) = {(x, y) in X times X : d_j(x, y) < r}, 2}
        U_{d_j}(r) & = & \{(x, y) \in X \times X : d_j(x, y) < r\}, \\
\label{U_{d_j'}(r) = {(x, y) in X times X : d_j'(x, y) < r}}
        U_{d_j'}(r) & = & \{(x, y) \in X \times X : d_j'(x, y) < r\}, 
\end{eqnarray}
and
\begin{equation}
\label{U_{d'}(r) = {(x, y) in X times X : d'(x, y) < r}}
        U_{d'}(r) = \{(x, y) \in X \times X : d'(x, y) < r\}.
\end{equation}
We also have that
\begin{eqnarray}
\label{U_{d_j'}(r) = U_{d_j}(r) when r le 1/j, = X times X when r > 1/j}
        U_{d_j'}(r) & = & U_{d_j}(r) \quad\hbox{ when } r \le 1/j \\
                   & = & X \times X \quad\hbox{when } r > 1/j, \nonumber
\end{eqnarray}
as in (\ref{U_{d_t}(r) = U_d(r)}) and (\ref{U_{d_t}(r) = X times X}).
By construction,
\begin{equation}
\label{U_{d'}(r) = bigcap_{j = 1}^infty U_{d_j'}(r)}
        U_{d'}(r) = \bigcap_{j = 1}^\infty U_{d_j'}(r)
\end{equation}
for every $r > 0$.  Combining this with (\ref{U_{d_j'}(r) = U_{d_j}(r)
  when r le 1/j, = X times X when r > 1/j}), we get that
\begin{equation}
\label{U_{d'}(r) = bigcap_{j = 1}^{[1/r]} U_{d_j}(r)}
        U_{d'}(r) = \bigcap_{j = 1}^{[1/r]} U_{d_j}(r)
\end{equation}
when $0 < r \le 1$, where $[1/r]$ is the largest positive integer
less than or equal to $1/r$.  If $r > 1$, then
\begin{equation}
\label{U_{d'}(r) = X times X}
        U_{d'}(r) = X \times X,
\end{equation}
because (\ref{d'(x, y) = max_{j ge 1} d_j'(x, y)}) is always less than
or equal to $1$.  Using (\ref{U_{d'}(r) = bigcap_{j = 1}^{[1/r]}
  U_{d_j}(r)}), one can check that the uniformity on $X$ associated to
$d'$ as in Section \ref{uniform structures} is the same as the
uniformity on $X$ associated to the sequence of semimetrics $d_1',
d_2', d_3', \ldots$, as in the previous section.  This is also the
same as the uniformity on $X$ associated to the initial sequence of
semimetrics $d_1, d_2, d_3, \ldots$, since the uniformities on $X$
associated to $d_j$ and to $d_j'$ are the same for each $j$, as
before.

        Let $d(\cdot, \cdot)$ be any semimetric on $X$ again, and
let $\mathcal{U}_d$ be the corresponding uniformity on $X$, as in
Section \ref{uniform structures}.  We have seen that the collection of
subsets of $X \times X$ of the form $U_d(r)$ as in (\ref{U(r) = U_d(r)
  = {(x, y) in X times X : d(x, y) < r}}) for some $r > 0$ forms a base
for $\mathcal{U}_d$.  More precisely, one can get a base for $\mathcal{U}_r$
using the sets $U_d(r)$ corresponding to a sequence of positive real
numbers $r$ converging to $0$.  In particular, there is a base for
$\mathcal{U}$ consisting of finitely or countably many subsets of
$X \times X$.  Conversely, if $\mathcal{U}$ is a uniformity on $X$,
and if there is a base for $\mathcal{U}$ with only finitely or countably
many elements, then it is well known that there is a semimetric on $X$
for which $\mathcal{U}$ is the associated uniformity, as on p186
of \cite{jk}.

\section{Collections of semi-ultrametrics}
\label{collections of semi-ultrametrics}
\setcounter{equation}{0}

        Let $X$ be a set, and let $U$ be a subset of $X \times X$.
Thus $U$ corresponds to a binary relation $x \sim y$ on $X$, as in
Section \ref{relations}.  Of course, $x \sim y$ is reflexive on $X$ if
and only if $U$ contains the diagonal $\Delta$ in (\ref{Delta =
  Delta_X = {(x, x) : x in X}}) as a subset.  Similarly, $x \sim y$ is
symmetric on $X$ if and only if $U$ is symmetric in the sense that
$\widetilde{U}$ in (\ref{widetilde{U} = {(x, y) : (y, x) in U}}) is
equal to $U$.  Transitivity of $x \sim y$ on $X$ may be expressed as
\begin{equation}
\label{U * U subseteq U}
        U * U \subseteq U,
\end{equation}
using the notation in (\ref{U * V = ...}).  If $x \sim y$ is reflexive
on $X$, so that $\Delta \subseteq X$, then we have that
\begin{equation}
\label{U subseteq U * U}
        U \subseteq U * U,
\end{equation}
by (\ref{U * Delta = Delta * U = U}).  It follows that
\begin{equation}
\label{U * U = U}
        U * U = U
\end{equation}
when $x \sim y$ is both transitive and reflexive on $X$.

        Let $U_1$, $U_2$ be subsets of $X \times X$, which correspond to
binary relations $x \sim_1 y$ and $x \sim_2 y$ on $X$, respectively,
as before.  Note that $U_1 \cap U_2$ corresponds to the binary
relation $x \sim y$ defined on $X$ by requiring that both $x \sim_1 y$
and $x \sim_2 y$ hold.  If $x \sim_1 y$ and $x \sim_2 y$ are
equivalence relations on $X$, then $x \sim y$ is an equivalence
relation on $X$ as well.

        Let $\mathcal{B}_0$ be a nonempty collection of subsets of
$X \times X$, each of which corresponds to an equivalence relation
on $X$.  Also let $\mathcal{B}$ be the collection of subsets of
$X \times X$ that can be expressed as the intersection of finitely
many elements of $\mathcal{B}_0$.  Every element of $\mathcal{B}$
corresponds to an equivalence relation on $X$ too, by the remarks
in the previous paragraph.  This means that for each $U \in \mathcal{B}$,
we have that $\Delta \subseteq U$, $\widetilde{U} = U$, and $U$ satisfies
(\ref{U * U subseteq U}).  Using this, it is easy to see that $\mathcal{B}$
is a base for a uniformity on $X$, since $\mathcal{B}$ is automatically
closed under finite intersections.

        Suppose that $\mathcal{M}$ is a nonempty collection of
semi-ultrametrics on $X$, and let $\mathcal{U}$ be the corresponding
uniformity on $X$, as in Section \ref{collections of semimetrics}.
Also let $U_d(r)$ be the subset of $X \times X$ associated to $d \in
\mathcal{M}$ and $r > 0$ as in (\ref{U(r) = U_d(r) = {(x, y) in X
    times X : d(x, y) < r}}).  Note that $U_d(r)$ corresponds to an
equivalence relation on $X$ for every $d \in \mathcal{M}$ and $r > 0$,
as in Section \ref{semi-ultrametrics}.  By construction, the collection
of subsets of $X \times X$ of the form $U_d(r)$ for some $d \in \mathcal{M}$
and $r > 0$ is a sub-base for the uniformity $\mathcal{U}$ on $X$
associated to $\mathcal{M}$.

        Let $x \sim y$ be any equivalence relation on $X$, and let
$d(x, y)$ be the corresponding discrete semi-ultrametric on $X$,
as in (\ref{d(x, y) = 0 when x sim y, = 1 when x not sim y}).
Also let $U_d(r)$ be the subset of $X \times X$ corresponding to $d(x,
y)$ and $r > 0$ as in (\ref{U(r) = U_d(r) = {(x, y) in X times X :
    d(x, y) < r}}).  Observe that
\begin{equation}
\label{U_d(r) = {(x, y) in X times X : x sim y}}
        U_d(r) = \{(x, y) \in X \times X : x \sim y\}
\end{equation}
when $0 < r \le 1$, and $U_d(r) = X \times X$ when $r > 1$.

        Let us now return to the situation where $\mathcal{B}_0$
is a nonempty collection of subsets of $X \times X$, each of which
corresponds to an equivalence relation on $X$.  The earlier remarks
imply that $\mathcal{B}_0$ is a sub-base for a uniformity $\mathcal{U}$
on $X$.
Let $\mathcal{M}_0$ be the
collection of discrete semi-ultrametrics on $X$ that correspond to
equivalence relations on $X$ associated to elements of $\mathcal{B}_0$,
as in (\ref{d(x, y) = 0 when x sim y, = 1 when x not sim y}).
Under these conditions, $\mathcal{U}$ is the same as the uniformity
on $X$ associated to $\mathcal{M}_0$ as in Section \ref{collections
of semimetrics}.

        If $\mathcal{B}_0$ has only finitely many elements, then
$\mathcal{M}_0$ has only finitely many elements too.  In this case,
the same uniformity on $X$ is determined by a single semi-ultrametric,
as in Section \ref{collections of semimetrics}.  Equivalently, one can
take the intersection of the elements of $\mathcal{B}_0$, to get
another subset of $X \times X$ that corresponds to an equivalence
relation on $X$.  The discrete semi-ultrametric on $X$ associated
to that subset of $X \times X$ determines the same uniformity on $X$.

        If $\mathcal{B}_0$ is countably infinite, then $\mathcal{M}_0$
is countably infinite as well.  In this case, one can get a
semi-ultrametric on $X$ that determines the same uniformity on $X$ as
in the previous section.  Conversely, let $d(\cdot, \cdot)$ be any
semi-ultrametric on $X$, and let $\mathcal{U}_d$ be the corresponding
uniformity on $X$, as in Section \ref{uniform structures}.  If
$U_d(r)$ is as in (\ref{U(r) = U_d(r) = {(x, y) in X times X : d(x, y)
    < r}}), then $U_d(r)$ corresponds to an equivalence relation on
$X$ for every $r > 0$, as in Section \ref{semi-ultrametrics}.  As in
the previous section, one can get a base for $\mathcal{U}_d$
consisting of the sets $U_d(r)$ for a sequence of positive real
numbers $r$ converging to $0$.

\section{$q$-Semimetrics}
\label{q-semimetrics}
\setcounter{equation}{0}

        Let $X$ be a set, and let $q$ be a positive real number.
A nonnegative real-valued function $d(x, y)$ defined on $X \times X$
is said to be a \emph{$q$-semimetric}\index{q-semimetrics@$q$-semimetrics}
if it satisfies (\ref{d(x, x) = 0}), (\ref{d(x, y) = d(y, x)}), and
\begin{equation}
\label{d(x, z)^q le d(x, y)^q + d(y, z)^q}
        d(x, z)^q \le d(x, y)^q + d(y, z)^q
\end{equation}
for every $x, y, z \in X$, instead of the ordinary triangle inequality
(\ref{d(x, z) le d(x, y) + d(y, z)}).  This is the same as saying that
$d(x, y)^q$ is an ordinary semimetric on $X$.  If $d(x, y)$ also
satisfies (\ref{d(x, y) > 0}), then it is said to be a
\emph{$q$-metric}\index{q-metrics@$q$-metrics} on $X$.  Of course,
(\ref{d(x, z)^q le d(x, y)^q + d(y, z)^q}) is equivalent to asking
that
\begin{equation}
\label{d(x, z) le (d(x, y)^q + d(y, z)^q)^{1/q}}
        d(x, z) \le (d(x, y)^q + d(y, z)^q)^{1/q}
\end{equation}
for every $x, y, z \in X$.

        Let $d(x, y)$ be a $q$-semimetric on $X$ for some $q > 0$,
and let $U_d(r) \subseteq X \times X$ be as in (\ref{U(r) = U_d(r) =
  {(x, y) in X times X : d(x, y) < r}}) for each $r > 0$.  This
satisfies (\ref{Delta subseteq U(r)}), (\ref{widetilde{U(r)} = U(r)}),
(\ref{U(r) subseteq U(t)}), and (\ref{bigcap_{r > 0} U(r) = {(x, y) in
    X times X : d(x, y) = 0}}), as in the case of ordinary semimetrics.
Instead of (\ref{U(r_1) * U(r_2) subseteq U(r_1 + r_2)}), we have that
\begin{equation}
\label{U_d(r_1) * U_d(r_2) subseteq U_d((r_1^q + r_2^q)^{1/q})}
        U_d(r_1) * U_d(r_2) \subseteq U_d((r_1^q + r_2^q)^{1/q})
\end{equation}
for every $r_1, r_2 > 0$, by (\ref{d(x, z) le (d(x, y)^q + d(y,
  z)^q)^{1/q}}).  Using this, one can define a uniformity
$\mathcal{U}_d$ on $X$ associated to $d(x, y)$ in the same way as for
ordinary semimetrics, as in Section \ref{uniform structures}.  If
\begin{equation}
\label{U_{d^q}(r) = {(x, y) in X times X : d(x, y)^q < r}}
        U_{d^q}(r) = \{(x, y) \in X \times X : d(x, y)^q < r\}
\end{equation}
is the analogous set associated to the semimetric $d(x, y)^q$ on $X$,
then
\begin{equation}
\label{U_{d^q}(r^q) = U_d(r)}
        U_{d^q}(r^q) = U_d(r)
\end{equation}
for every $r > 0$.  This implies that $\mathcal{U}_d$ is the same as
the uniformity $\mathcal{U}_{d^q}$ on $X$ associated to $d(x, y)^q$ as
in Section \ref{uniform structures}.  In particular, the corresponding
topologies on $X$ are the same.

        Let $\mathcal{M}$ be a nonempty collection of $q$-semimetrics
on $X$, where one can let $q > 0$ depend on the element of $\mathcal{M}$.
One can define a uniformity on $X$ using $\mathcal{M}$ in the same way
as in Section \ref{collections of semimetrics}.  Alternatively, one
can get a collection of ordinary semimetrics on $X$ by replacing
each $d \in \mathcal{M}$ with $d(x, y)^q$ for a suitable $q > 0$.
Such a collection of ordinary semimetrics leads to the same uniformity
on $X$, because of (\ref{U_{d^q}(r^q) = U_d(r)}).

        If $a$, $b$ are nonnegative real numbers, then
\begin{equation}
\label{max(a, b) le (a^q + b^q)^{1/q} le 2^{1/q} max(a, b)}
        \max(a, b) \le (a^q + b^q)^{1/q} \le 2^{1/q} \, \max(a, b)
\end{equation}
for every $q > 0$.  If $0 < q \le q' < \infty$, then we get that
\begin{equation}
\label{a^{q'} + b^{q'} le ... = (a^q + b^q)^{q'/q}}
 \qquad a^{q'} + b^{q'} \le \max(a, b)^{q' - q} \, (a^q + b^q)
                        \le (a^q + b^q)^{(q' - q)/q + 1} = (a^q + b^q)^{q'/q}.
\end{equation}
Hence
\begin{equation}
\label{(a^{q'} + b^{q'})^{1/q'} le (a^q + b^q)^{1/q}}
        (a^{q'} + b^{q'})^{1/q'} \le (a^q + b^q)^{1/q}.
\end{equation}
It follows that every $q'$-semimetric on $X$ is a $q$-semimetric on
$X$ when $q \le q'$, using also (\ref{d(x, z) le (d(x, y)^q + d(y,
  z)^q)^{1/q}}).  Note that a semi-ultrametric on $X$ is a
$q$-semimetric for each $q > 0$, because of the first inequality in
(\ref{max(a, b) le (a^q + b^q)^{1/q} le 2^{1/q} max(a, b)}).  Using
both inequalities in (\ref{max(a, b) le (a^q + b^q)^{1/q} le 2^{1/q}
  max(a, b)}), we have that
\begin{equation}
\label{lim_{q to infty} (a^q + b^q)^{1/q} = max(a, b)}
        \lim_{q \to \infty} (a^q + b^q)^{1/q} = \max(a, b)
\end{equation}
for every $a, b \ge 0$.  Thus one can think of a semi-ultrametric
as being a $q$-semimetric with $q = \infty$.

\section{$q$-Absolute value functions}
\label{q-absolute value functions}
\setcounter{equation}{0}

        Let $k$ be a field, and let $q$ be a positive real number.
A nonnegative real-valued function $|\cdot|$ on $k$ is said to be an
\emph{$q$-absolute value function}\index{q-absolute value
  functions@$q$-absolute value functions} on $k$ if it satisfies the
following three conditions.  First,
\begin{equation}
\label{|x| = 0 if and only if x = 0}
        |x| = 0 \hbox{ if and only if } x = 0.
\end{equation}
Second,
\begin{equation}
\label{|x y| = |x| |y|}
        |x \, y| = |x| \, |y|
\end{equation}
for every $x, y \in k$.  Third,
\begin{equation}
\label{|x + y|^q le |x|^q + |y|^q}
        |x + y|^q \le |x|^q + |y|^q
\end{equation}
for every $x, y \in k$.  As before, (\ref{|x + y|^q le |x|^q + |y|^q})
is equivalent to asking that
\begin{equation}
\label{|x + y| le (|x|^q + |y|^q)^{1/q}}
        |x + y| \le (|x|^q + |y|^q)^{1/q}
\end{equation}
for every $x, y \in k$.  A $q$-absolute value function on $k$ with $q
= 1$ is also known simply as an absolute value function\index{absolute
  value functions} on $k$.  Thus $|x|$ is a $q$-absolute value
function on $k$ if and only if $|x|^q$ is an absolute value function
on $k$.  Of course, the standard absolute value functions on the
fields ${\bf R}$\index{R@${\bf R}$} and ${\bf C}$\index{C@${\bf C}$}
of real and complex numbers are absolute value functions in this
sense.  If $0 < q \le q' < \infty$ and $|\cdot|$ is a $q'$-absolute
value function on a field $k$, then $|\cdot|$ is a $q$-absolute value
function on $k$ too, by (\ref{(a^{q'} + b^{q'})^{1/q'} le (a^q +
  b^q)^{1/q}}).

        Similarly, a nonnegative real-valued function $|\cdot|$ on
a field $k$ is said to be an \emph{ultrametric absolute value
  function}\index{ultrametric absolute value functions} on $k$ if it
satisfies (\ref{|x| = 0 if and only if x = 0}), (\ref{|x y| = |x|
  |y|}), and
\begin{equation}
\label{|x + y| le max(|x|, |y|)}
        |x + y| \le \max(|x|, |y|)
\end{equation}
for every $x, y \in k$.  An ultrametric absolute value function on $k$
is a $q$-absolute value function for every $q > 0$, by the first
inequality in (\ref{max(a, b) le (a^q + b^q)^{1/q} le 2^{1/q} max(a,
  b)}).  One can also think of an ultrametric absolute value function
on $k$ as being a $q$-absolute value function with $q = \infty$,
because of (\ref{lim_{q to infty} (a^q + b^q)^{1/q} = max(a, b)}).
The \emph{trivial absolute value function}\index{trivial absolute
  value function} is defined on any field $k$ by putting $|0| = 0$ and
\begin{equation}
\label{|x| = 1}
        |x| = 1
\end{equation}
 for every $x \in k$ with $x \ne 0$.  It is easy to see that this is
 an ultrametric absolute value function on $k$.

        Suppose that $|\cdot|$ is a nonnegative real-valued function
on a field $k$ that satisfies (\ref{|x| = 0 if and only if x = 0}) and
(\ref{|x y| = |x| |y|}).  If $1$ is the multiplicative identity
element in $k$, then $1 \ne 0$ in $k$, by the definition of a field,
and hence $|1| > 0$.  We also have that $1^2 = 1$ in $k$, so that $|1|
= |1^2| = |1|^2$, and hence
\begin{equation}
\label{|1| = 1}
        |1| = 1.
\end{equation}
If $x \in k$ satisfies $x^n = 1$ for some positive integer $n$, then
\begin{equation}
\label{|x|^n = |x^n| = |1| = 1}
        |x|^n = |x^n| = |1| = 1,
\end{equation}
so that $|x| = 1$.  In particular,
\begin{equation}
\label{|-1| = 1}
        |-1| = 1,
\end{equation}
which implies that
\begin{equation}
\label{|-x| = |(-1) x| = |-1| |x| = |x|}
        |-x| = |(-1) \, x| = |-1| \, |x| = |x|
\end{equation}
for every $x \in k$.

        If $|\cdot|$ is a $q$-absolute value function on $k$
for some $q > 0$, then it follows that
\begin{equation}
\label{d(x, y) = |x - y|}
        d(x, y) = |x - y|
\end{equation}
defines a $q$-metric on $k$.  Similarly, if $|\cdot|$ is an
ultrametric absolute value function on $k$, then (\ref{d(x, y) = |x -
  y|}) is an ultrametric on $k$.  The ultrametric corresponding to the
trivial absolute value function on $k$ as in (\ref{d(x, y) = |x - y|})
is the discrete metric.

\section{$q$-Seminorms}
\label{q-seminorms}
\setcounter{equation}{0}

        Let $k$ be a field, let $|\cdot|$ be a $q$-absolute value
function on $k$ for some $q > 0$, and let $V$ be a vector space over
$k$.  A nonnegative real-valued function $N$ on $V$ is said to be a
\emph{$q$-seminorm}\index{q-seminorms@$q$-seminorms} on $V$ if
\begin{equation}
\label{N(t v) = |t| N(v)}
        N(t \, v) = |t| \, N(v)
\end{equation}
for every $v \in V$ and $t \in k$, and
\begin{equation}
\label{N(v + w)^q le N(v)^q + N(w)^q}
        N(v + w)^q \le N(v)^q + N(w)^q
\end{equation}
for every $v, w \in V$.  Of course, (\ref{N(t v) = |t| N(v)}) implies
that $N(0) = 0$, by taking $t = 0$.  If we also have that
\begin{equation}
\label{N(v) > 0}
        N(v) > 0
\end{equation}
for every $v \in V$ with $v \ne 0$, then $N$ is said to be a
\emph{$q$-norm}\index{q-norms@$q$-norms} on $V$.

        If $q = 1$, then a $q$-seminorm may simply be called a
\emph{seminorm},\index{seminorms} and a $q$-norm may be called a
\emph{norm}.\index{norms}  Remember that $|\cdot|$ is a $q$-absolute
value function on $k$ if and only if $|\cdot|^q$ an absolute value
function on $k$.  Similarly, $N$ is a $q$-seminorm on $V$ with respect
to $|\cdot|$ on $k$ if and only if $N^q$ is a seminorm on $V$
with respect to $|\cdot|^q$ on $k$.  Of course, there is an analogous
statement for $q$-norms.

        As usual, (\ref{N(v + w)^q le N(v)^q + N(w)^q}) is the same as
asking that
\begin{equation}
\label{N(v + w) le (N(v)^q + N(w))^{1/q}}
        N(v + w) \le (N(v)^q + N(w))^{1/q}
\end{equation}
for every $v, w \in V$.  Suppose for the moment that $0 < q \le q' \le
\infty$, and that $|\cdot|$ is a $q'$-absolute value function on $k$.
This implies that $|\cdot|$ is a $q$-absolute value function on $k$ as
well, as in the previous section.  Similarly, if $N$ is a
$q'$-seminorm on $V$, then $N$ is a $q$-seminorm on $V$, because of
(\ref{(a^{q'} + b^{q'})^{1/q'} le (a^q + b^q)^{1/q}}).  The analogous
statement for $q$-norms follows by including the additional condition
(\ref{N(v) > 0}).

        Suppose now that $|\cdot|$ is an ultrametric absolute value
function on $k$.  A nonnegative real-valued function $N$ is said to
be a \emph{semi-ultranorm}\index{semi-ultranorms} on $V$ if $N$ satisfies
(\ref{N(t v) = |t| N(v)}) and
\begin{equation}
\label{N(v + w) le max(N(v), N(w))}
        N(v + w) \le \max(N(v), N(w))
\end{equation}
for every $v, w \in V$.  If $N$ also satisfies (\ref{N(v) > 0}), then
$N$ is said to be an \emph{ultranorm}\index{ultranorms} on $V$.  Note
that (\ref{N(v + w) le max(N(v), N(w))}) implies (\ref{N(v + w) le
  (N(v)^q + N(w))^{1/q}}) for every $q > 0$, by the first inequality
in (\ref{max(a, b) le (a^q + b^q)^{1/q} le 2^{1/q} max(a, b)}).  Thus
a semi-ultranorm on $V$ is a $q$-seminorm on $V$ for every $q > 0$,
and similarly for ultranorms.  This implicitly uses the fact that an
ultrametric absolute value function on $k$ is a $q$-absolute value
function on $k$ for every $q > 0$, as in the previous section.  One
can think of semi-ultranorms and ultranorms as being the $q = \infty$
versions of $q$-seminorms and $q$-norms, respectively, because of
(\ref{lim_{q to infty} (a^q + b^q)^{1/q} = max(a, b)}).

        Let $|\cdot|$ be a $q$-absolute value function on $k$ for
some $q > 0$ again.  If $N$ is a $q$-seminorm on $V$, then
\begin{equation}
\label{d(v, w) = N(v - w)}
        d(v, w) = N(v - w)
\end{equation}
is a $q$-semimetric on $V$.  If $N$ is a $q$-norm on $V$, then
(\ref{d(v, w) = N(v - w)}) is a $q$-metric on $V$.  Similarly, if
$|\cdot|$ is an ultrametric absolute value function on $k$, and if $N$
is a semi-ultranorm on $V$, then (\ref{d(v, w) = N(v - w)}) is a
semi-ultrametric on $V$.  Under the same conditions, if $N$ is an
ultranorm on $V$, then (\ref{d(v, w) = N(v - w)}) is an ultrametric on
$V$.

        Let $|\cdot|$ be the trivial absolute value function on $k$,
which is an ultrametric absolute value function on $k$, as in the
previous section.  The \emph{trivial ultranorm}\index{trivial ultranorm}
is defined by putting $N(0) = 0$ and
\begin{equation}
\label{N(v) = 1}
        N(v) = 1
\end{equation}
for every $v \in V$ with $v \ne 0$.  It is easy to see that this is an
ultranorm on $V$, for which the corresponding metric as in (\ref{d(v,
  w) = N(v - w)}) is the same as the discrete metric.

        We have included the requirement that $|\cdot|$ be a $q$-absolute
value function on $k$ in the definition of a $q$-seminorm for
convenience, and to avoid trivialities.  If $N(v) > 0$ for some $v \in
V$, then (\ref{N(t v) = |t| N(v)}) and (\ref{N(v + w)^q le N(v)^q +
  N(w)^q}) imply (\ref{|x + y|^q le |x|^q + |y|^q}), as well as
(\ref{|x y| = |x| |y|}).  Similarly, (\ref{N(t v) = |t| N(v)}) and
(\ref{N(v + w) le max(N(v), N(w))}) imply (\ref{|x + y| le max(|x|,
  |y|)}) in this case.

\part{Connectedness and dimension $0$}
\label{connectedness, dimension 0}

\section{Connected sets}
\label{connected sets}
\setcounter{equation}{0}

        As usual, a pair $A$, $B$ of subsets of a topological space $X$
are said to be \emph{separated}\index{separated sets} in $X$ if
\begin{equation}
\label{overline{A} cap B = A cap overline{B} = emptyset}
        \overline{A} \cap B = A \cap \overline{B} = \emptyset,
\end{equation}
where $\overline{A}$, $\overline{B}$ are the closures of $A$, $B$ in
$X$, respectively.  In particular, disjoint closed subsets of $X$ are
automatically separated, and one can check that disjoint open subsets
of $X$ are separated as well.  A set $E \subseteq X$ is said to be
\emph{connected}\index{connected sets} if it cannot be expressed as $A
\cup B$, where $A$, $B$ are nonempty separated subsets of $X$.
Suppose that $Y \subseteq X$ is equipped with the topology induced by
the given topology on $X$.  It is well known that $A, B \subseteq Y$
are separated as subsets of $Y$ with respect to the induced topology
if and only if $A$ and $B$ are separated in $X$.  This is because the
closures of $A$, $B$ in $Y$ with respect to the induced topology are
the same as the intersections of $Y$ with the closures of $A$, $B$ in
$X$.  It follows that $E \subseteq Y$ is connected with respect to the
induced topology on $Y$ if and only if $E$ is connected as a subset of
$X$.

        If $A$, $B$ are separated subsets of $X$ such that
\begin{equation}
\label{A cup B = X}
        A \cup B = X,
\end{equation}
then it is easy to see that $A$ and $B$ are each both open and closed
as subsets of $X$.  It follows that $X$ is connected as a subset of
itself if and only if it cannot be expressed as the union of two
nonempty disjoint open sets, which is the same as saying that $X$
cannot be expressed as the union of two nonempty disjoint closed sets.
Of course, $A \subseteq X$ is both open and closed in $X$ exactly when
$A$ and $X \setminus A$ are both open in $X$, which is the same as
saying that they are both closed in $X$.  Thus $X$ is connected if and
only if there are no nonempty proper subsets of $X$ that are both open
and closed.  If $E \subseteq X$, then the remarks in the previous
paragraph can be applied to $Y = E$, to get that $E$ is connected as a
subset of $X$ if and only if $E$ is connected as a subset of itself,
with respect to the induced topology.

        Let $Z$ be another topological space, and suppose that $f$
is a continuous mapping from $X$ into $Z$.  If $E$ is a connected
subset of $X$, then it is well known that $f(E)$ is connected in $Z$.
As in the previous paragraph, one may as well suppose that $E = X$,
since otherwise one can simply restrict $f$ to $E$.  Similarly,
one can restrict one's attention to the case where $f$ is surjective,
since otherwise one can replace $Z$ with $f(X)$, where $f(X)$ is equipped
with the topology induced by the one on $Z$.  With these reductions,
one can use the characterization of connectedness in terms of open
or closed sets, as in the preceding paragraph.

        Suppose that $E = A \cup B$, where $A$, $B$ are separated subsets
of $X$.  Observe that
\begin{equation}
\label{A = overline{A} cap E, B = overline{B} cap E}
        A = \overline{A} \cap E, \quad B = \overline{B} \cap E.
\end{equation}
If $E$ is compact in $X$, then it follows that $A$ and $B$ are compact
in $X$ too, because the intersection of a closed set and a compact set
is compact.  Of course, if $E$ is not connected, then we can take $A$
and $B$ to be nonempty.  If $X$ is Hausdorff, then it is well known
that compact subsets of $X$ are closed in $X$.

\section{$U$-Separated sets}
\label{U-separated sets}
\setcounter{equation}{0}

        Let $X$ be a set, let $A$, $B$ be subsets of $X$, and let
$U$ be a subset of $X \times X$.  Suppose that
\begin{equation}
\label{Delta subseteq U, 2}
        \Delta \subseteq U,
\end{equation}
where $\Delta$ is the diagonal in $X \times X$, as in (\ref{Delta =
  Delta_X = {(x, x) : x in X}}).  Let us say that $A$, $B$ are
\emph{$U$-separated}\index{U-separated sets@$U$-separated sets} in $X$
if for every $x \in A$ and $y \in B$ we have that
\begin{equation}
\label{(x, y) not in U}
        (x, y) \not\in U.
\end{equation}
Of course, this implies that $A$ and $B$ are disjoint, because of
(\ref{Delta subseteq U, 2}).

        Equivalently, $A$, $B$ are $U$-separated in $X$ when
\begin{equation}
\label{U cap (A times B) = emptyset}
        U \cap (A \times B) = \emptyset.
\end{equation}
This is the same as saying that
\begin{equation}
\label{U[A] cap B = emptyset}
        U[A] \cap B = \emptyset,
\end{equation}
where $U[A]$ is as in (\ref{U[A] = {y in X : there is an x in A such
    that (x, y) in U}}).  We can also reformulate (\ref{U cap (A times
  B) = emptyset}) as saying that
\begin{equation}
\label{U subseteq (X times X) setminus (A times B)}
        U \subseteq (X \times X) \setminus (A \times B).
\end{equation}

        Note that $A$, $B$ are $U$-separated if and only if $B$, $A$
are $\widetilde{U}$-separated, where $\widetilde{U}$ is as in
(\ref{widetilde{U} = {(x, y) : (y, x) in U}}).  This is the same as
saying that (\ref{U cap (A times B) = emptyset}) holds if and only if
\begin{equation}
\label{widetilde{U} cap (B times A) = emptyset}
        \widetilde{U} \cap (B \times A) = \emptyset,
\end{equation}
which is equivalent to
\begin{equation}
\label{widetilde{U} subseteq (X times X) setminus (B times A)}
        \widetilde{U} \subseteq (X \times X) \setminus (B \times A).
\end{equation}
As before, (\ref{widetilde{U} cap (B times A) = emptyset}) is also
equivalent to
\begin{equation}
\label{A cap (widetilde{U}[B]) = emptyset}
        A \cap (\widetilde{U}[B]) = \emptyset,
\end{equation}
so that (\ref{U[A] cap B = emptyset}) and (\ref{A cap
  (widetilde{U}[B]) = emptyset}) are equivalent as well.

        If $U$ is symmetric, in the sense that $\widetilde{U} = U$,
then the condition that $A$, $B$ be $U$-separated is symmetric in $A$
and $B$.  In this case, the combination of (\ref{U cap (A times B) =
  emptyset}) and (\ref{widetilde{U} cap (B times A) = emptyset}) is
equivalent to
\begin{equation}
\label{U cap ((A times B) cup (B times A)) = emptyset}
        U \cap ((A \times B) \cup (B \times A)) = \emptyset,
\end{equation}
which is the same as saying that
\begin{equation}
\label{U subseteq (X times X) setminus ((A times B) cup (B times A))}
 U \subseteq (X \times X) \setminus ((A \times B) \cup (B \times A)).
\end{equation}

        Let $d(\cdot, \cdot)$ be a $q$-semimetric on $X$ for some
$q > 0$.  We say that $A, B \subseteq X$ are
\emph{$r$-separated}\index{r-separated sets@$r$-separated sets}
in $X$ with respect to $d(\cdot, \cdot)$ if
\begin{equation}
\label{d(x, y) ge r}
        d(x, y) \ge r
\end{equation}
for every $x \in A$ and $y \in B$.  If $U = U_d(r)$ is as in
(\ref{U(r) = U_d(r) = {(x, y) in X times X : d(x, y) < r}}), then this
is the same as asking that $A$, $B$ be $U$-separated, as before.

        Let $U$ be any subset of $X \times X$ that satisfies
(\ref{Delta subseteq U, 2}) again, and let $U_1$, $U_2$ be subsets of 
$X \times X$ such that
\begin{equation}
\label{U_1 * widetilde{U_2} subseteq U}
        U_1 * \widetilde{U_2} \subseteq U.
\end{equation}
Observe that
\begin{equation}
\label{widetilde{U_2}[U_1[A]] = (U_1 * widetilde{U_2})[A] subseteq U[A]}
        \widetilde{U_2}[U_1[A]] = (U_1 * \widetilde{U_2})[A] \subseteq U[A]
\end{equation}
for every $A \subseteq X$, using (\ref{(U * V)[A] = V[U[A]]}) in the
first step, and (\ref{U_1 * widetilde{U_2} subseteq U}) in the second
step.  If $A, B \subseteq X$, then the equivalence between (\ref{U[A]
  cap B = emptyset}) and (\ref{A cap (widetilde{U}[B]) = emptyset})
implies that
\begin{equation}
\label{(widetilde{U_2}[U_1[A]]) cap B = emptyset}
 (\widetilde{U_2}[U_1[A]]) \cap B = \emptyset
\end{equation}
holds if and only if
\begin{equation}
\label{(U_1[A]) cap (U_2[B]) = emptyset}
        (U_1[A]) \cap (U_2[B]) = \emptyset.
\end{equation}
More precisely, $U_1[A]$ plays the role here that $A$ had before, and
$\widetilde{U_2}$ plays the role that $U$ had before.  If $A$, $B$ are
$U$-separated in $X$, then (\ref{U[A] cap B = emptyset}) and
(\ref{widetilde{U_2}[U_1[A]] = (U_1 * widetilde{U_2})[A] subseteq
  U[A]}) imply that (\ref{(widetilde{U_2}[U_1[A]]) cap B = emptyset})
holds, so that (\ref{(U_1[A]) cap (U_2[B]) = emptyset}) holds as well.
If we have that
\begin{equation}
\label{U_1 * widetilde{U_2} = U}
        U_1 * \widetilde{U_2} = U
\end{equation}
instead of (\ref{U_1 * widetilde{U_2} subseteq U}), then the inclusion
in the second step in (\ref{widetilde{U_2}[U_1[A]] = (U_1 *
  widetilde{U_2})[A] subseteq U[A]}) can be replaced with an equality.
In this case, (\ref{(U_1[A]) cap (U_2[B]) = emptyset}) implies that
$A$, $B$ are $U$-separated in $X$, by essentially the same argument.

        Suppose that $U_{1, 1}, U_{1, 2} \subseteq X \times X$ satisfy
\begin{equation}
\label{U_{1, 1} * U_{1, 2} subseteq U_1}
        U_{1, 1} * U_{1, 2} \subseteq U_1,
\end{equation}
so that
\begin{equation}
\label{U_{1, 2}[U_{1, 1}[A]] subseteq U_1[A]}
        U_{1, 2}[U_{1, 1}[A]] \subseteq U_1[A],
\end{equation}
as in (\ref{(U * V)[A] = V[U[A]]}).  Under these conditions,
(\ref{(U_1[A]) cap (U_2[B]) = emptyset}) implies that
\begin{equation}
\label{(U_{1, 2}[U_{1, 1}[A]]) cap (U_2[B]) = emptyset}
        (U_{1, 2}[U_{1, 1}[A]]) \cap (U_2[B]) = \emptyset.
\end{equation}
This means that $U_{1, 1}[A]$ and $U_2[B]$ are $U_{1, 2}$-separated in
$X$, at least if $\Delta \subseteq U_{1, 2}$, so that this condition
is defined.

\section{Uniformly separated sets}
\label{uniformly separated sets}
\setcounter{equation}{0}

        Let $(X, \mathcal{U})$ be a uniform space, and let $A$, $B$
be subsets of $X$.  Let us say that $A$, $B$ are \emph{uniformly
  separated}\index{uniformly separated sets}\index{separated
  sets!uniformly} in $X$ if there is a $U \in \mathcal{U}$ such that
$A$, $B$ are $U$-separated, as in the previous section.  Remember that
every $U \in \mathcal{U}$ satisfies (\ref{Delta subseteq U, 2}), by
definition of a uniformity.

        As in the previous section, $A$ and $B$ are $U$-separated in $X$
if and only if $B$ and $A$ are $\widetilde{U}$-separated in $X$.  Thus
the property of being uniformly separated in $X$ is symmetric in $A$
and $B$, because of (\ref{widetilde{U} in mathcal{U}}).  Of course, if
$A$, $B$ are $U$-separated for some $U \subseteq X \times X$, then
$A$, $B$ satisfy the analogous condition with respect to any subset of
$U$.  If $U \in \mathcal{U}$, then $U \cap \widetilde{U}$ is a
symmetric element of $\mathcal{U}$ that is contained in $U$.  Thus we
may as well ask that $U$ be symmetric in the definition of uniformly
separated subsets of $X$.

        We have seen that $A$ and $B$ are $U$-separated in $X$ if and
only if (\ref{U subseteq (X times X) setminus (A times B)}) holds.
If $U \in \mathcal{U}$, then it follows that
\begin{equation}
\label{(X times X) setminus (A times B)}
        (X \times X) \setminus (A \times B)
\end{equation}
is an element of $\mathcal{U}$.  Conversely, if (\ref{(X times X)
  setminus (A times B)}) is an element of $\mathcal{U}$, then we
can simply take $U$ to be (\ref{(X times X) setminus (A times B)}),
so that $A$ and $B$ are automatically $U$-separated in $X$.
Thus $A$ and $B$ are uniformly separated in $X$ with respect to
$\mathcal{U}$ if and only if (\ref{(X times X) setminus (A times B)})
is an element of $\mathcal{U}$.

        Using (\ref{widetilde{U} in mathcal{U}}), we get that
(\ref{(X times X) setminus (A times B)}) is an element of $\mathcal{U}$
if and only if
\begin{equation}
\label{(X times X) setminus (B times A)}
        (X \times X) \setminus (B \times A)
\end{equation}
is an element of $\mathcal{U}$.  Of course, this corresponds to the
fact that the property of being uniformly separated is symmetric in
$A$ and $B$.  In this case, we get that
\begin{eqnarray}
\label{... = (X times X) setminus ((A times B) cup (B times A))}
\lefteqn{((X \times X) \setminus (A \times B)) \cap 
                         ((X \times X) \setminus (B \times A))} \\
 & = & (X \times X) \setminus ((A \times B) \cup (B \times A)) \nonumber
\end{eqnarray}
is an element of $\mathcal{U}$ as well.

        Let $\overline{A}$, $\overline{B}$ be the closures of $A, B \subseteq
X$ with respect to the topology on $X$ associated to $\mathcal{U}$
as in Section \ref{associated topology}.  Thus
\begin{equation}
\label{overline{A} subseteq U[A], overline{B} subseteq widetilde{U}[B]}
 \overline{A} \subseteq U[A], \quad \overline{B} \subseteq \widetilde{U}[B]
\end{equation}
for every $U \in \mathcal{U}$, as in Section \ref{closure,
  regularity}.  Suppose that $A$, $B$ are uniformly separated in $X$,
so that there is a $U \in \mathcal{U}$ such that (\ref{U[A] cap B =
  emptyset}) and (\ref{A cap (widetilde{U}[B]) = emptyset}) hold.
This implies that $A$, $B$ are separated in $X$ with respect to
associated topology as in (\ref{overline{A} cap B = A cap overline{B}
  = emptyset}), by (\ref{overline{A} subseteq U[A], overline{B}
  subseteq widetilde{U}[B]}).

        In fact, we have that
\begin{equation}
\label{overline{A} cap overline{B} = emptyset}
        \overline{A} \cap \overline{B} = \emptyset
\end{equation}
under these conditions.  To see this, we can first use the definition
of a uniformity, to get $U_1, U_2 \in \mathcal{U}$ satisfying
(\ref{U_1 * widetilde{U_2} subseteq U}).  If $A$, $B$ are
$U$-separated, then it follows that (\ref{(U_1[A]) cap (U_2[B]) =
  emptyset}) holds as well.  This implies (\ref{overline{A} cap
  overline{B} = emptyset}), for the same reasons as before.

        Alternatively, (\ref{overline{A} cap overline{B} = emptyset})
is the same as saying that
\begin{equation}
\label{(overline{A} times overline{B}) cap Delta = emptyset}
        (\overline{A} \times \overline{B}) \cap \Delta = \emptyset,
\end{equation}
where $\Delta$ is as in (\ref{Delta = Delta_X = {(x, x) : x in X}}),
as usual.  Of course, $\overline{A} \times \overline{B}$ is the same
as the closure of $A \times B$ in $X \times X$, with respect to the
product topology corresponding to the topology on $X$ associated to
$\mathcal{U}$.  If $A$, $B$ are $U$-separated in $X$ and $\Delta$ is
contained in the interior of $U$ with respect to the product topology
on $X \times X$, then (\ref{(overline{A} times overline{B}) cap Delta
  = emptyset}) follows from (\ref{U cap (A times B) = emptyset}).  If
$U \in \mathcal{U}$, then one can check that $\Delta$ is contained in
the interior of $U$ with respect to the product topology on $X \times
X$.  More precisely, the interior of $U$ with respect to the product
topology on $X \times X$ is an element of $\mathcal{U}$ too, as in
Theorem 6 on p179 of \cite{jk}.

        As another variant, suppose that $A$ and $B$ are uniformly
separated in $X$ again, and let $U$, $U_1$, and $U_2$ be elements of
$\mathcal{U}$ that satisfy (\ref{U[A] cap B = emptyset}) and (\ref{U_1
  * widetilde{U_2} subseteq U}), as before.  Similarly, there are
$U_{1, 1}, U_{1, 2} \in \mathcal{U}$ that satisfy (\ref{U_{1, 1} *
  U_{1, 2} subseteq U_1}), by the definition of a uniformity.  One
might as well take $U_{1, 1} = U_{1, 2}$ here, as in (\ref{V * V
  subseteq U}).  If $A$, $B$ are $U$-separated, then (\ref{(U_1[A])
  cap (U_2[B]) = emptyset}) and hence (\ref{(U_{1, 2}[U_{1, 1}[A]])
  cap (U_2[B]) = emptyset}) hold, so that $U_{1, 1}[A]$ and $U_2[B]$
are $U_{1, 2}$-separated in $X$.  In particular, this means that
$U_{1, 1}[A]$ and $U_2[B]$ are uniformly separated in $X$.
As in Section \ref{closure, regularity} again, we have that
\begin{equation}
\label{overline{A} subseteq U_{1, 1}[A], overline{B} subseteq U_2[B]}
 \overline{A} \subseteq U_{1, 1}[A], \quad \overline{B} \subseteq U_2[B].
\end{equation}
This implies that $\overline{A}$, $\overline{B}$ are uniformly
separated in $X$.

        Suppose that $K$, $E$ are disjoint subsets of $X$ such that
$K$ is compact and $E$ is closed with respect to the topology on $X$
associated to $\mathcal{U}$.  Thus $X \setminus E$ is an open set
that contains $K$, and hence there is a $U \in \mathcal{U}$ such that
\begin{equation}
\label{U[K] subseteq X setminus E}
        U[K] \subseteq X \setminus E,
\end{equation}
as in (\ref{V[K] subseteq W}).  This is the same as saying that
$U[K]$ is disjoint from $E$, so that $K$, $E$ are $U$-separated,
as in (\ref{U[A] cap B = emptyset}).  This shows that $K$ and $E$
are uniformly separated in $X$ under these conditions.

        Let $Y$ be a subset of $X$, equipped with the uniform
structure induced by $\mathcal{U}$ on $X$ as in Section \ref{induced
  uniform structures}.  If $A$, $B$ are subsets of $Y$, then one can
check that $A$ and $B$ are uniformly separated in $Y$ with respect to
the induced uniform structure if and only if $A$ and $B$ are uniformly
separated in $X$ with respect to $\mathcal{U}$.

\section{Equivalence relations}
\label{equivalence relations}
\setcounter{equation}{0}

        Let $X$ be a set, and suppose for the moment that
$U \subseteq X \times X$ corresponds to an equivalence relation
on $X$.  If $A \subseteq X$ and $U[A]$ is as in (\ref{U[A] = {y in X :
    there is an x in A such that (x, y) in U}}), then $U[A]$ is the
same as the union of the equivalence classes in $X$ associated to $U$
that contain an element of $A$.  Thus
\begin{equation}
\label{U[U[A]] = U[A]}
        U[U[A]] = U[A],
\end{equation}
which can also be viewed in terms of (\ref{(U * V)[A] = V[U[A]]}) and
(\ref{U * U = U}).  In particular, $A$ is itself a union of
equivalence classes in $X$ associated to $U$ if and only if
\begin{equation}
\label{U[A] = A}
        U[A] = A.
\end{equation}

        As in Section \ref{U-separated sets}, $A, B \subseteq X$ are
$U$-separated in $X$ when (\ref{U[A] cap B = emptyset}) holds.  In
this case, this is the same as saying that $B$ is disjoint from the
equivalence classes in $X$ associated to $U$ that contain an element
of $A$.  Because $U$ is symmetric, this is equivalent to
\begin{equation}
\label{A cap (U[B]) = emptyset}
        A \cap (U[B]) = \emptyset,
\end{equation}
which means that $A$ is disjoint from the equivalence classes in $X$
associated to $U$ that contain an element of $B$.  This implies that
\begin{equation}
\label{(U[A]) cap (U[B]) = emptyset}
        (U[A]) \cap (U[B]) = \emptyset,
\end{equation}
which is the same as saying that the equivalence classes in $X$
associated to $U$ that contain elements of $A$ and $B$, respectively,
are disjoint.  Note that (\ref{(U[A]) cap (U[B]) = emptyset})
corresponds to (\ref{(U_1[A]) cap (U_2[B]) = emptyset}) with
$U_1 = U_2 = U$.

        Let $d(\cdot, \cdot)$ be a semi-ultrametric on $X$, and let $U_d(r)$
be as in (\ref{U(r) = U_d(r) = {(x, y) in X times X : d(x, y) < r}})
for each $r > 0$.  Thus $U_d(r)$ corresponds to an equivalence relation
on $X$ for every $r > 0$, as in Section \ref{semi-ultrametrics}.  The
equivalence classes in $X$ associated to $U_d(r)$ are the same as open
balls in $X$ of radius $r$ with respect to $d(\cdot, \cdot)$.  By the
remarks in the preceding paragraph, any two distinct equivalence classes
in $X$ with respect to $U_d(r)$ are $U_d(r)$-separated in $X$.  In this
case, this means that any two open balls in $X$ of radius $r$ with
respect to $d(\cdot, \cdot)$ that correspond to distinct subsets of $X$
are $r$-separated in $X$ with respect to $d(\cdot, \cdot)$.

        Similarly,
\begin{equation}
\label{d(x, y) le r, 2}
        d(x, y) \le r
\end{equation}
defines an equivalence relation on $X$ for every $r \ge 0$ when
$d(\cdot, \cdot)$ is a semi-ultrametric on $X$, as in Section
\ref{semi-ultrametrics}.  The corresponding equivalence classes in $X$
are the closed balls in $X$ with respect to $d(\cdot, \cdot)$ with
radius $r$.  As before, if $B_1$ and $B_2$ are distinct subsets of $X$
that can be expressed as closed balls of radius $r$ with respect to
$d(\cdot, \cdot)$, then
\begin{equation}
\label{d(x, y) > r}
        d(x, y) > r
\end{equation}
for every $x \in B_1$ and $y \in B_2$.  This is the same as saying
that $B_1$ and $B_2$ are separated with respect to the equivalence
relation (\ref{d(x, y) le r, 2}).  This also works for any semimetric
$d(\cdot, \cdot)$ on $X$ when $r = 0$.

        Suppose that $U \subseteq X \times X$ corresponds to an
equivalence relation on $X$ again, and let $U_0$ be another subset
of $X \times X$ that contains the diagonal $\Delta$.  If
\begin{equation}
\label{U_0 subseteq U}
        U_0 \subseteq U,
\end{equation}
then any pair of distinct equivalence classes in $X$ with respect to
$U$ are $U_0$-separated in $X$.  Conversely, if every pair of distinct
equivalence classes in $X$ with respect to $U$ are $U_0$-separated in
$X$, then (\ref{U_0 subseteq U}) holds.  If $\mathcal{U}$ is a
uniformity on $X$, $U_0 \in \mathcal{U}$, and (\ref{U_0 subseteq U})
holds, then $U \in \mathcal{U}$ too.  Thus $U \in \mathcal{U}$ when
distinct equivalence classes in $X$ with respect to $U$ are uniformly
separated with respect to a single $U_0 \in \mathcal{U}$.

        Let $A$ be any subset of $X$, and consider
\begin{equation}
\label{U_A = (A times A) cup ((X setminus A) times (X setminus A))}
        U_A = (A \times A) \cup ((X \setminus A) \times (X \setminus A))
\end{equation}
as a subset of $X \times X$.  This corresponds to the equivalence
relation $x \sim_A y$ on $X$ which is satisfied when either both $x$
and $y$ are elements of $A$ or both $x$ and $y$ are not in $A$.
Observe that
\begin{equation}
\label{(X times X) setminus U_A = ...}
 (X \times X) \setminus U_A = (A \times (X \setminus A)) \cup
                                          ((X \setminus A) \times A).
\end{equation}
This means that $U_A$ is the same as (\ref{... = (X times X) setminus
  ((A times B) cup (B times A))}), with $B = X \setminus A$.  It
follows that $A$ is uniformly separated from $X \setminus A$ with
respect to a uniformity $\mathcal{U}$ on $X$ if and only if
\begin{equation}
\label{U_A in mathcal{U}}
        U_A \in \mathcal{U},
\end{equation}
as in the previous section.

        Let $Y$ be another set, let $f$ be a mapping from $X$ into $Y$,
and let $f_2$ be the corresponding mapping from $X \times X$ into $Y
\times Y$, as in (\ref{f_2(x, x') = (f(x), f(x'))}).  If $V \subseteq
Y \times Y$ corresponds to an equivalence relation on $Y$, then
\begin{equation}
\label{f_2^{-1}(V)}
        f_2^{-1}(V)
\end{equation}
corresponds to an equivalence relation on $X$.  In particular, if $X
\subseteq Y$, then we can take $f(x) = x$ for each $x \in X$.
In this case, (\ref{f_2^{-1}(V)}) reduces to
\begin{equation}
\label{V cap (X times X)}
        V \cap (X \times X).
\end{equation}

\section{$U$-Chains}
\label{U-chains}
\setcounter{equation}{0}

        Let $X$ be a set, and let $U$ be a subset of $X \times X$.
Suppose that
\begin{equation}
\label{Delta subseteq U, 3}
        \Delta \subseteq U
\end{equation}
and
\begin{equation}
\label{widetilde{U} = U, 2}
        \widetilde{U} = U,
\end{equation}
where $\Delta$ and $\widetilde{U}$ are as in (\ref{Delta = Delta_X =
  {(x, x) : x in X}}) and (\ref{widetilde{U} = {(x, y) : (y, x) in
    U}}), respectively.  This is the same as saying that the binary
relation on $X$ corresponding to $U$ is reflexive and symmetric.

        Let $n$ be a positive integer, and let $x_1, \ldots, x_n$ be
a finite sequence of elements of $X$ of length $n$.  Let us say that
$x_1, \ldots, x_n$ is a \emph{$U$-chain}\index{U-chains@$U$-chains}
of length $n$ in $X$ if
\begin{equation}
\label{(x_j, x_{j + 1}) in U}
        (x_j, x_{j + 1}) \in U
\end{equation}
for each $j = 1, \ldots, n - 1$.  This condition is considered to be
vacuous when $n = 1$.

        If $n$ is a nonnegative integer, then let $U^n$ be the set
of $(x, x') \in X \times X$ for which there is a $U$-chain
$x_1, \ldots, x_{n + 1}$ of length $n + 1$ such that $x_1 = x$ and
$x_{n + 1} = x'$.  Thus $U^0 = \Delta$ and $U^1 = U$, by construction.
If $k$, $l$ are nonnegative integers, then it is easy to see that
\begin{equation}
\label{U^k * U^l = U^{k + l}}
        U^k * U^l = U^{k + l},
\end{equation}
where the left side of (\ref{U^k * U^l = U^{k + l}}) is as defined in
(\ref{U * V = ...}).  Note that
\begin{equation}
\label{U^n subseteq U^{n + 1}}
        U^n \subseteq U^{n + 1}
\end{equation}
for each $n \ge 0$, because of (\ref{Delta subseteq U, 3}), and using
either (\ref{U^k * U^l = U^{k + l}}) or simply the definition of
$U^n$.  One can also check that $U^n$ is symmetric for every $n \ge
0$, which is to say that
\begin{equation}
\label{widetilde{(U^n)} = U^n}
        \widetilde{(U^n)} = U^n,
\end{equation}
using (\ref{widetilde{U} = U, 2}).

        Put
\begin{equation}
\label{widehat{U} = bigcup_{n = 0}^infty U^n}
        \widehat{U} = \bigcup_{n = 0}^\infty U^n,
\end{equation}
which is a subset of $X \times X$ that contains $U$, and hence
$\Delta$.  Equivalently, $\widehat{U}$ consists of the $(x, x') \in X
\times X$ such that $x$ can be connected to $x'$ by a $U$-chain
of elements of $X$ of some finite length.  Observe that
\begin{equation}
\label{widehat{U} * widehat{U} subseteq widehat{U}}
        \widehat{U} * \widehat{U} \subseteq \widehat{U},
\end{equation}
and in fact
\begin{equation}
\label{widehat{U} * widehat{U} = widehat{U}}
        \widehat{U} * \widehat{U} = \widehat{U}.
\end{equation}
More precisely,
\begin{equation}
\label{widehat{U} subseteq widehat{U} * widehat{U}}
        \widehat{U} \subseteq \widehat{U} * \widehat{U}
\end{equation}
because $\Delta \subseteq \widehat{U}$.  The reverse inclusion
basically says that one can combine two $U$-chains in $X$ to get
another $U$-chain in $X$ when the first $U$-chain ends at the point
where the second $U$-chain begins.  One can also derive
(\ref{widehat{U} * widehat{U} = widehat{U}}) from (\ref{U^k * U^l =
  U^{k + l}}).  It is easy to see that $\widehat{U}$ is symmetric too,
so that
\begin{equation}
\label{widetilde{(widehat{U})} = widehat{U}}
        \widetilde{(\widehat{U})} = \widehat{U},
\end{equation}
because of (\ref{widetilde{(U^n)} = U^n}).

        Equivalently, (\ref{widehat{U} * widehat{U} subseteq widehat{U}})
says that the binary relation on $X$ corresponding to $\widehat{U}$ is
transitive.  It follows that this is an equivalence relation on $X$,
since it is reflexive and symmetric as well.  If the binary relation
on $X$ corresponding to $U$ is already transitive, and hence an
equivalence relation, then
\begin{equation}
\label{U^n = U}
        U^n = U
\end{equation}
for every $n \ge 1$, so that
\begin{equation}
\label{widehat{U} = U}
        \widehat{U} = U.
\end{equation}

        Let $U$ be as before, not necessarily corresponding to a
transitive relation on $X$.  Let $x \in X$ be given, and consider the
set $\widehat{U}[x]$, as in (\ref{U[x] = U[{x}] = {y in X : (x, y) in
    U}}).  This is the same as the set of points in $X$ that can be
connected to $x$ by a $U$-chain of finite length.  This can also be
described as the equivalence class in $X$ containing $x$ determined by
the equivalence relation corresponding to $\widehat{U}$.  Let us say
that $X$ is
\emph{$U$-connected}\index{U-connectedness@$U$-connectedness} if
\begin{equation}
\label{widehat{U} = X times X}
        \widehat{U} = X \times X,
\end{equation}
which means that every pair of elements of $X$ can be connected by a
$U$-chain of finite length.  This is the same as saying that
$\widehat{U}[x] = X$ for every $x \in X$.  If $\widehat{U}[x] = X$ for
any $x \in X$, then this holds for every $x \in X$, and hence $X$ is
$U$-connected, because $\widehat{U}$ corresponds to an equivalence
relation on $X$.

        Similarly, let us say that $E \subseteq X$ is
\emph{$U$-connected}\index{U-connectedness@$U$-connectedness} if
every pair of elements of $E$ can be connected by a $U$-chain of
elements of $E$ of finite length.  This corresponds to replacing
$X$ with $E$ in the previous discussion, and replacing $U$ with
\begin{equation}
\label{U cap (E times E)}
        U \cap (E \times E).
\end{equation}
Note that (\ref{U cap (E times E)}) contains the diagonal in $E \times
E$ and is symmetric, because of (\ref{Delta subseteq U, 3}) and
(\ref{widetilde{U} = U, 2}).

        Let $d(\cdot, \cdot)$ be a semimetric on $X$, or a $q$-semimetric
on $X$ for some $q > 0$, and let $r$ be a positive real number.  If $U
= U_d(r)$ is as in (\ref{U(r) = U_d(r) = {(x, y) in X times X : d(x,
    y) < r}}), then $U$ satisfies (\ref{Delta subseteq U, 3}) and
(\ref{widetilde{U} = U, 2}), by hypothesis.  In this case, a $U$-chain
in $X$ may be described as an
\emph{$r$-chain}\index{r-chains@$r$-chains} in $X$ with respect to
$d(\cdot, \cdot)$.  If $d(\cdot, \cdot)$ is a semi-ultrametric on $X$,
then $U_d(r)$ corresponds to an equivalence relation on $X$, as in
Section \ref{semi-ultrametrics}.

\section{Chain connectedness}
\label{chain connectedness}
\setcounter{equation}{0}

        Let $(X, \mathcal{U})$ be a uniform space, and let $U$ be an
element of $\mathcal{U}$ that is also symmetric.  Thus $U$ satisfies
(\ref{Delta subseteq U, 3}) too, by definition of a uniformity.  As in
the previous section, $X$ is said to be $U$-connected if
(\ref{widehat{U} = X times X}) holds.  Let us say that $X$ is
\emph{chain connected}\index{chain connectedness} as a uniform space
with respect to $\mathcal{U}$ if $X$ is $U$-connected for every $U \in
\mathcal{U}$ such that $U$ is symmetric.  As before, this means that
every pair of elements of $X$ can be connected by a $U$-chain of elements
of $X$ of finite length.

        Let us check that $X$ is not chain connected if and only if
there is a proper nonempty subset $A$ of $X$ such that $A$ is uniformly
separated from $X \setminus A$ with respect to $\mathcal{U}$.
Suppose that $A$ is a proper nonempty subset of $X$ that is uniformly
separated from $X \setminus A$ with respect to $\mathcal{U}$, so that
there is a $U \in \mathcal{U}$ such that $A$ and $X \setminus A$ are
$U$-separated.  We may as well ask that $U$ be symmetric too, as in
Section \ref{uniformly separated sets}.  It is easy to see that there
is no $U$-chain of elements of $X$ that connects an element of $A$ to
an element of $X \setminus A$, because $A$ and $X \setminus A$ are
$U$-separated in $X$.  This implies that $X$ is not $U$-connected,
and hence that $X$ is not chain connected with respect to $\mathcal{U}$,
because $A \ne X, \emptyset$ by hypothesis.

        Conversely, suppose that $X$ is not chain connected with respect
to $\mathcal{U}$.  This means that there is a $U \in \mathcal{U}$ such
that $U$ is symmetric, and that $X$ is not $U$-connected.  More
precisely, it follows that there are elements $x$, $y$ of $X$ that
cannot be connected by a $U$-chain of elements of $X$ of finite
length.  Let $A$ be the set of $z \in X$ for which there is a
$U$-chain of elements of $X$ of finite length from $x$ to $z$.  Thus
$x \in A$ automatically, since $U$ satisfies (\ref{Delta subseteq U,
  3}) by the definition of a uniformity.  We also have that $y \not\in
A$, by hypothesis, so that $A \ne X$.  It is easy to see that $A$ is
$U$-separated from $X \setminus A$ under these conditions, by
construction.  It follows that $A$ is uniformly separated from $X
\setminus A$ with respect to $\mathcal{U}$, as desired.

        Equivalently, $X$ is not chain connected with respect to
$\mathcal{U}$ if there is a $U \in \mathcal{U}$ such that $U$
is symmetric and
\begin{equation}
\label{widehat{U} ne X times X}
        \widehat{U} \ne X \times X,
\end{equation}
where $\widehat{U}$ is as in (\ref{widehat{U} = bigcup_{n = 0}^infty
  U^n}).  Remember that $U \subseteq \widehat{U}$ by construction, so
that $U \in \mathcal{U}$ implies that $\widehat{U} \in \mathcal{U}$.
If $X$ is not chain connected with respect to $\mathcal{U}$, then it
follows that $\widehat{U}$ is an element of $\mathcal{U}$ that
corresponds to an equivalence relation on $X$ and satisfies
(\ref{widehat{U} ne X times X}).  Of course, if $U$ already
corresponds to an equivalence relation on $X$, then $\widehat{U} = U$,
as in (\ref{widehat{U} = U}).

        Thus $X$ is not chain connected with respect to $\mathcal{U}$
if and only if there is a $U \in \mathcal{U}$ such that $U$ corresponds
to an equivalence relation on $X$ and
\begin{equation}
\label{U ne X times X}
        U \ne X \times X.
\end{equation}
This means that $X$ is chain connected with respect to $\mathcal{U}$
if and only if the only $U \in \mathcal{U}$ that corresponds to an
equivalence relation on $X$ is $U = X \times X$.

        Let us see how this reformulation of chain connectedness is
related to the previous one.  Let $A$ be a subset of $X$, and let
$U_A$ be the corresponding subset of $X \times X$ defined in (\ref{U_A
  = (A times A) cup ((X setminus A) times (X setminus A))}).  If $A$
is uniformly separated from its complement $X \setminus A$ in $X$ with
respect to $\mathcal{U}$, then $U_A$ is an element of $\mathcal{U}$,
as in (\ref{U_A in mathcal{U}}).  We have also seen that $U_A$
automatically corresponds to an equivalence relation on $X$.  If $A$
is a nonempty proper subset of $X$, then
\begin{equation}
\label{U_A ne X times X}
        U_A \ne X \times X,
\end{equation}
as in (\ref{U ne X times X}).

        Conversely, suppose that there is a $U \in \mathcal{U}$ that
corresponds to an equivalence relation on $X$ and satisfies (\ref{U ne
  X times X}), as before.  In particular, (\ref{U ne X times X})
implies that $X \ne \emptyset$, and we let $A$ be the equivalence
class associated to any element of $X$.  Thus $A \ne \emptyset$, by
construction, and $A \ne X$, because of (\ref{U ne X times X}) again.
We also have that $A$ and $X \setminus A$ are $U$-separated in $X$, as
in Section \ref{equivalence relations}.  It follows that $A$ and $X
\setminus A$ are uniformly separated in $X$ with respect to
$\mathcal{U}$, because $U \in \mathcal{U}$ by hypothesis.

\section{Chain connectedness, continued}
\label{chain connectedness, continued}
\setcounter{equation}{0}

        Let $(X, \mathcal{U})$ be a uniform space again, and let $E$
be a subset of $X$.  Let us say that $E$ is \emph{chain
  connected}\index{chain connectedness} in $X$ if for each $U \in
\mathcal{U}$ such that $U$ is symmetric, we have that $E$ is
$U$-connected in $X$.  This means that every pair of elements of $E$
can be connected by a $U$-chain of elements of $E$ of finite length,
as in Section \ref{U-chains}.  In particular, this reduces to the
definition of chain connectedness in the previous section when
$E = X$.

        Let $Y$ be a subset of $X$, and let $\mathcal{V}$ be the
uniformity induced on $Y$ by $\mathcal{U}$ on $X$, as in Section
\ref{induced uniform structures}.  Thus $\mathcal{V}$ consists
of the $V \subseteq Y \times Y$ that can be expressed as
\begin{equation}
\label{V = U cap (Y times Y)}
        V = U \cap (Y \times Y),
\end{equation}
with $U \in \mathcal{U}$.  Clearly (\ref{V = U cap (Y times Y)})
implies that
\begin{equation}
\label{widetilde{V} = widetilde{U} cap (Y times Y)}
        \widetilde{V} = \widetilde{U} \cap (Y \times Y),
\end{equation}
using the notation in (\ref{widetilde{U} = {(x, y) : (y, x) in U}}).
If $U$ is symmetric and $V$ is as in (\ref{V = U cap (Y times Y)}),
then it follows that $V$ is symmetric too.  In the other direction, if
$V \in \mathcal{V}$ is symmetric, then we can choose $U \in
\mathcal{U}$ so that (\ref{V = U cap (Y times Y)}) holds and $U$ is
symmetric, by replacing $U$ with $U \cap \widetilde{U}$ if necessary.

        If $E \subseteq Y$, then it is easy to see that $E$ is chain
connected with respect to $\mathcal{V}$ on $Y$ if and only if $E$ is
chain connected with respect to $\mathcal{U}$ on $X$.  This uses
the remarks in the previous paragraph, to get that symmetric elements
of $\mathcal{V}$ correspond to symmetric elements of $\mathcal{U}$.
In particular, we can take $Y = E$, to reduce to the situation
considered in the preceding section.

        Equivalently, $E \subseteq X$ is not chain connected if and
only if there are nonempty subsets $A$, $B$ of $X$ that are uniformly
separated with respect to $\mathcal{U}$ and satisfy
\begin{equation}
\label{E = A cup B}
        E = A \cup B.
\end{equation}
This can be verified in essentially the same way as in the previous
section, when $E = X$.  Alternatively, if $E \subseteq Y \subseteq X$,
then one can check that $E$ has this property with respect to
$\mathcal{U}$ on $X$ if and only if $E$ has the analogous property
with respect to the induced uniformity $\mathcal{V}$ on $Y$.  This
uses the corresponding statement for uniformly separated sets
mentioned in Section \ref{uniformly separated sets}.  Since the other
formulation of chain connectedness has the same feature, as in the
preceding paragraph, the equivalence between the two formulations can
be reduced to the case where $E = X$.

        If $A, B \subseteq X$ are uniformly separated with respect to
$\mathcal{U}$, then $A$, $B$ are separated subsets of $X$ with respect to
the topology associated to $\mathcal{U}$, as in Section \ref{uniformly
  separated sets}.  If $E \subseteq X$ is not chain connected with
respect to $\mathcal{U}$, then it follows that $E$ is not connected
with respect to this topology on $X$.  Equivalently, if $E$ is
connected with respect to this topology, then $E$ is chain connected
with respect to $\mathcal{U}$.  It is easy to give examples of subsets
of the real line that are chain connected but not connected, using the
uniformity determined by the standard Euclidean metric on ${\bf R}$.

        Let $(Z, \mathcal{W})$ be another uniform space, and let $f$
be a uniformly continuous mapping from $X$ into $Z$, as in Section
\ref{uniform continuity}.  If $A$, $B$ are uniformly separated subsets
of $Z$, then one can check that $f^{-1}(A)$ and $f^{-1}(B)$ are
uniformly separated in $X$.  If $E \subseteq X$ is chain connected,
then it follows that $f(E)$ is chain connected in $Z$.  One can also
look at this more directly in terms of the initial definition of
chain connectedness.

        Suppose that $E \subseteq X$ is not connected with respect to
the topology on $X$ associated to $\mathcal{U}$ as in Section
\ref{associated topology}, so that $E$ can be expressed as $A \cup B$,
where $A$, $B$ are nonempty separated subsets of $X$.  If $E$ is
compact in $X$, then $A$, $B$ are compact too, as in Section
\ref{connected sets}.  Note that the closure $\overline{B}$ of $B$ in
$X$ is a closed set in $X$ that is disjoint from $A$, because $A$ and
$B$ are supposed to be separated in $X$.  Using the compactness of
$A$, we get that $A$ and $\overline{B}$ are uniformly separated in $X$
with respect to $\mathcal{U}$, as in Section \ref{uniformly separated
  sets}.  In particular, this implies that $A$ and $B$ are uniformly
separated in $X$.  It follows that $E$ is not chain connected in $X$
with respect to $\mathcal{U}$, because $A$ and $B$ are nonempty.
Equivalently, if $E \subseteq X$ is compact with respect to the
topology on $X$ associated to $\mathcal{U}$, and if $E$ is chain
connected with respect to $\mathcal{U}$, then $E$ is connected with
respect to the topology on $X$ associated to $\mathcal{U}$.

\section{Connectedness and closure}
\label{connectedness, closure}
\setcounter{equation}{0}

        If $E$ is a connected subset of a topological space $X$, then
it is well known that the closure $\overline{E}$ of $E$ in $X$ is
connected too.  Equivalently, if $\overline{E}$ is not connected, then
$E$ is not connected.  Of course, if $\overline{E}$ is not connected,
then there are nonempty separated sets $A_1, B_1 \subseteq X$ such
that
\begin{equation}
\label{A_1 cup B_1 = overline{E}}
        A_1 \cup B_1 = \overline{E}.
\end{equation}
If we put
\begin{equation}
\label{A = A_1 cap E, B = B_1 cap E}
        A = A_1 \cap E, \quad B = B_1 \cap E,
\end{equation}
then it is easy to see that $A$, $B$ are separated subsets of $X$
whose union is $E$.  The remaining point is to verify that $A, B \ne
\emptyset$ under these conditions.

        Let $(X, \mathcal{U})$ be a uniform space, and let $E$ be a
subset of $X$ again.  Also let $\overline{E}$ be the closure of $E$
with respect to the topology on $X$ associated to $\mathcal{U}$ as in
Section \ref{associated topology}.  If $\overline{E}$ is not chain
connected in $X$ with respect to $\mathcal{U}$, then there are
nonempty uniformly separated subsets $A_1$, $B_1$ of $X$ that satisfy
(\ref{A_1 cup B_1 = overline{E}}).  If $A, B \subseteq X$ are as in
(\ref{A = A_1 cap E, B = B_1 cap E}), then $A$, $B$ are uniformly
separated sets whose union is $E$.  As before, one can check that
$A, B \ne \emptyset$, to get that $E$ is not chain connected
with respect to $\mathcal{U}$.

        Equivalently, if $E$ is chain connected with respect to
$\mathcal{U}$, then $\overline{E}$ is chain connected with respect
to $\mathcal{U}$.  Alternatively, let $U$ be an arbitrary symmetric
element of $\mathcal{U}$.  If $x$ is any element of $\overline{E}$,
then there is a $w \in E$ such that
\begin{equation}
\label{(x, w) in U}
        (x, w) \in U,
\end{equation}
because of (\ref{x in Int(U[x])}).  If $E$ is $U$-connected, as in
Section \ref{U-chains}, then one can use this to verify that
$\overline{E}$ is $U$-connected too.  If $E$ is chain connected in
$X$, then it follows that $\overline{E}$ is chain connected in $X$ as
well.

        Suppose now that $E$ is not chain connected in $X$, so that
there are nonempty uniformly separated sets $A, B \subseteq X$ whose
union is equal to $E$.  This implies that the closures $\overline{A}$,
$\overline{B}$ of $A$, $B$ are uniformly separated in $X$, as in
Section \ref{uniformly separated sets}.  Of course, $\overline{A},
\overline{B} \ne \emptyset$, because $A, B \ne \emptyset$, and it is
well known that
\begin{equation}
\label{overline{E} = overline{A} cup overline{B}}
        \overline{E} = \overline{A} \cup \overline{B}
\end{equation}
when $E = A \cup B$.  It follows that $\overline{E}$ is not chain
connected in $X$ when $E$ is not chain connected in $X$, which is the
same as saying that $E$ is chain connected in $X$ when $\overline{E}$
is chain connected in $X$.  Note that the analogous statement for
ordinary connectedness in topological spaces does not hold.

\section{Two related uniformities}
\label{two related uniformities}
\setcounter{equation}{0}

        Let $(X, \mathcal{U})$ be a uniform space, and put
\begin{equation}
\label{mathcal{B}_{eq} = ...}
 \mathcal{B}_{eq} = \{U \in \mathcal{U}:
                    U \hbox{ corresponds to an equivalence relation on } X\}.
\end{equation}
Note that
\begin{equation}
\label{X times X in mathcal{B}_{eq}}
        X \times X \in \mathcal{B}_{eq},
\end{equation}
because $X \times X \in \mathcal{U}$ by the definition of a
uniformity, and $X \times X$ corresponds trivially to an equivalence
relation on $X$.  In particular, $\mathcal{B}_{eq} \ne \emptyset$.  If
$U_1, U_2 \in \mathcal{B}_{eq}$, then it is easy to see that
\begin{equation}
\label{U_1 cap U_2 in mathcal{B}_{eq}}
        U_1 \cap U_2 \in \mathcal{B}_{eq},
\end{equation}
using the corresponding property (\ref{U cap V in mathcal{U}}) for
$\mathcal{U}$, and the analogous statement for equivalence relations,
mentioned in Section \ref{collections of semi-ultrametrics}.  It
follows that $\mathcal{B}_{eq}$ is a base for a uniformity
\begin{equation}
\label{mathcal{U}_{eq}}
        \mathcal{U}_{eq}
\end{equation}
on $X$, as in Section \ref{collections of semi-ultrametrics} again.

        Similarly, put
\begin{eqnarray}
\label{mathcal{B}_{eq, 2} = ...}
 \mathcal{B}_{eq, 2} & = & \{U \in \mathcal{B}_{eq} :
 \hbox{ there are only finitely many equivalence} \\
 & & \qquad\qquad\quad\hbox{ classes in } X \hbox{ corresponding to } U\},
                                                   \nonumber
\end{eqnarray}
which is to say that there are only finitely many equivalence classes
in $X$ determined by the equivalence relation on $X$ corresponding to
each $U \in \mathcal{B}_{eq, 2}$.  As before,
\begin{equation}
\label{X times X in mathcal{B}_{eq, 2}}
        X \times X \in \mathcal{B}_{eq, 2},
\end{equation}
because $X$ itself is the only equivalence class determined by the
trivial equivalence relation corresponding to $X \times X$.  If $U_1,
U_2 \in \mathcal{B}_{eq, 2}$, then one can check that
\begin{equation}
\label{U_1 cap U_2 in mathcal{B}_{eq, 2}}
        U_1 \cap U_2 \in \mathcal{B}_{eq, 2},
\end{equation}
as in (\ref{U_1 cap U_2 in mathcal{B}_{eq}}).  More precisely, the
number of equivalence classes in $X$ corresponding to $U_1 \cap U_2$
is less than or equal to the product of the numbers of equivalence
classes in $X$ corresponding to $U_1$ and $U_2$, respectively.  It
follows again that $\mathcal{B}_{eq, 2}$ is a base for a uniformity
\begin{equation}
\label{mathcal{U}_{eq, 2}}
        \mathcal{U}_{eq, 2}
\end{equation}
on $X$, as in Section \ref{collections of semi-ultrametrics}.

        By construction,
\begin{equation}
\label{mathcal{B}_{eq, 2} subseteq mathcal{B}_{eq} subseteq mathcal{U}}
        \mathcal{B}_{eq, 2} \subseteq \mathcal{B}_{eq} \subseteq \mathcal{U},
\end{equation}
and hence
\begin{equation}
\label{mathcal{U}_{eq, 2} subseteq mathcal{U}_{eq} subseteq mathcal{U}}
        \mathcal{U}_{eq, 2} \subseteq \mathcal{U}_{eq} \subseteq \mathcal{U}.
\end{equation}
In particular, this implies that the topology on $X$ associated to
$\mathcal{U}_{eq, 2}$ as in Section \ref{associated topology} is
contained in the topology on $X$ associated to $\mathcal{U}_{eq}$,
which is itself contained in the topology on $X$ associated to
$\mathcal{U}$.  The topology on $X$ associated to $\mathcal{U}_{eq,
  2}$ is actually the same as the topology on $X$ associated to
$\mathcal{U}_{eq}$, as we shall see in a moment.

        Let $U_A$ be as in (\ref{U_A = (A times A) cup ((X setminus A) times 
(X setminus A))}) for each $A \subseteq X$, which corresponds to an
equivalence relation on $X$.  If $A$ is a proper nonempty subset of $X$,
then there are exactly two equivalence classes in $X$ corresponding to
$U_A$, which are $A$ and $X \setminus A$.  Otherwise, if $A = \emptyset$
or $A = X$, then $U_A = X \times X$.  Remember that $A$ is uniformly
separated from $X \setminus A$ with respect to $\mathcal{U}$ if and
only if $U_A \in \mathcal{U}$, as in (\ref{U_A in mathcal{U}}).
In this case, we get that
\begin{equation}
\label{U_A in mathcal{B}_{eq, 2}}
        U_A \in \mathcal{B}_{eq, 2},
\end{equation}
since there are at most two equivalence classes in $X$ corresponding
to $U_A$.

        In particular, if $A \subseteq X$ is uniformly separated from
$X \setminus A$ with respect to $\mathcal{U}$, then $A$ is uniformly
separated from $X \setminus A$ with respect to $\mathcal{U}_{eq, 2}$,
by (\ref{U_A in mathcal{B}_{eq, 2}}).  Of course, if $A$ is uniformly
separated from $X \setminus A$ with respect to $\mathcal{U}_{eq, 2}$,
then $A$ is uniformly separated from $X \setminus A$ with respect to
$\mathcal{U}_{eq}$, because of the first inclusion in
(\ref{mathcal{U}_{eq, 2} subseteq mathcal{U}_{eq} subseteq
  mathcal{U}}).  Similarly, if $A$ is uniformly separated with respect
to $X \setminus A$ with respect to $\mathcal{U}_{eq}$, then $A$ is
uniformly separated from $X \setminus A$ with respect to
$\mathcal{U}$, by the second inclusion in (\ref{mathcal{U}_{eq, 2}
  subseteq mathcal{U}_{eq} subseteq mathcal{U}}).  Thus $\mathcal{U}$,
$\mathcal{U}_{eq}$, and $\mathcal{U}_{eq, 2}$ determine the same
collection of subsets $A$ of $X$ that are uniformly separated from
their complements in $X$.

        Suppose that $A \subseteq X$ is an equivalence class corresponding
to an element $V$ of $\mathcal{B}_{eq}$.  This implies that $A$ is
uniformly separated from $X \setminus A$ with respect to
$\mathcal{U}$, so that (\ref{U_A in mathcal{B}_{eq, 2}}) holds.  Using
this, one can check that the topology on $X$ associated to
$\mathcal{U}_{eq}$ is the same as the topology on $X$ associated to
$\mathcal{U}_{eq, 2}$.

        Of course, (\ref{U_A in mathcal{B}_{eq, 2}}) says that
\begin{equation}
\label{{U_A : A subseteq X is uniformly separated from X setminus A ...}}
 \qquad \{U_A : A \subseteq X \hbox{ is uniformly separated from }
                X \setminus A \hbox{ with respect to } \mathcal{U}\}
\end{equation}
is contained in $\mathcal{B}_{eq, 2}$.  If $V$ is any element of
$\mathcal{B}_{eq, 2}$ and $A \subseteq X$ is an equivalence class
associated to $V$, then $U_A$ is an element of (\ref{{U_A : A subseteq
    X is uniformly separated from X setminus A ...}}), as in the
previous paragraph.  In this case, there are only finitely many such
equivalence classes, by the definition (\ref{mathcal{B}_{eq, 2} =
  ...}) of $\mathcal{B}_{eq, 2}$.  This implies that $V$ can be
expressed as the intersection of finitely many elements of (\ref{{U_A
    : A subseteq X is uniformly separated from X setminus A ...}}).
It follows that (\ref{{U_A : A subseteq X is uniformly separated from
    X setminus A ...}})  is a sub-base for $\mathcal{U}_{eq, 2}$.

        Observe that
\begin{equation}
\label{X is totally bounded with respect to mathcal{U}_{eq, 2}}
        X \hbox{ is totally bounded with respect to } \mathcal{U}_{eq, 2}
\end{equation}
automatically.  More precisely, in order to verify that $X$ is totally
bounded with respect to $\mathcal{U}_{eq, 2}$, it suffices to check
that $X$ satisfies the condition in Section \ref{totally bounded sets}
for each $U \in \mathcal{B}_{eq, 2}$, because $\mathcal{B}_{eq, 2}$ is
a base for $\mathcal{U}_{eq, 2}$.  In this case, this condition
follows from the requirement that there be only finitely many
equivalence classes in $X$ corresponding to $U$.  If $X$ is totally
bounded with respect to $\mathcal{U}_{eq}$, then it is easy to see
that there can only be finitely many equivalence classes in $X$
corresponding to any $U \in \mathcal{B}_{eq}$.  This implies that
\begin{equation}
\label{mathcal{U}_{eq} = mathcal{U}_{eq, 2}}
        \mathcal{U}_{eq} = \mathcal{U}_{eq, 2},
\end{equation}
which holds in particular when $X$ is totally bounded with respect to
$\mathcal{U}$.

\section{Topological dimension $0$}
\label{topological dimension 0}
\setcounter{equation}{0}

        A topological space $X$ is said to have \emph{topological 
dimension $0$}\index{topological dimension 0@topological dimension $0$}
at a point $x \in X$ if for every open set $W \subseteq X$ with $x \in W$
there is an open set $U \subseteq X$ such that $x \in U$, $U \subseteq W$,
and $U$ is also a closed set in $X$.  Equivalently, this means that
there is a local base for the topology of $X$ at $x$ consisting of
subsets of $X$ that are both open and closed.

        If $X$ has topological dimension $0$ at every point $x \in X$,
then one simply says that $X$ has topological dimension
$0$,\index{topological dimension 0@topological dimension $0$}
at least in a strict sense.  This is the same as saying that there
is a base for the topology of $X$ consisting of subsets of $X$
that are both open and closed.  Sometimes one also requires $X$
to be nonempty to have topological dimension $0$, in order to define
related conditions in other dimensions inductively, but we shall not
pursue this here.  If $X$ has topological dimension $0$ in this
strict sense, then $X$ is regular in the strict sense, as in Section
\ref{closure, regularity}.  In particular, if $X$ satisfies the
first or $0$th separation condition as well, then $X$ is Hausdorff,
as before.

        Suppose now that $(X, \mathcal{U})$ is a uniform space, so that
$X$ is equipped with the topology associated to $\mathcal{U}$ as in
Section \ref{associated topology} too.  If $A \subseteq X$ is
uniformly separated from $X \setminus A$ with respect to
$\mathcal{U}$, then $A$ and $X \setminus A$ are separated in $X$ with
respect to this topology on $X$, as in Section \ref{uniformly
  separated sets}, which means that $A$ is both open and closed in
$X$.  Let us say that $X$ is \emph{strongly $0$-dimensional}\index{strongly
0-dimensional uniform spaces@strongly $0$-dimensional uniform
spaces}\index{uniform spaces!strongly 0-dimensional@strongly $0$-dimensional}
at a point $x \in X$ if for each open set $W \subseteq X$ with $x \in
W$ there is a set $U \subseteq X$ such that $x \in U$, $U \subseteq W$,
and $U$ is uniformly separated from $X \setminus U$ in $X$ with respect
to $\mathcal{U}$.  This implies that $X$ has topological dimension $0$
at $x$, since $U$ is both open and closed in $X$ under these conditions.
Equivalently, $X$ is strongly $0$-dimensional at $x$ if there is a
local base for the topology of $X$ at $x$ consisting of subsets of $X$
that are uniformly separated from their complements in $X$ with respect
to $\mathcal{U}$.

        If $X$ is strongly $0$-dimensional at every point $x \in X$, then
let us say that $X$ is strongly $0$-dimensional\index{strongly 0-dimensional
uniform spaces@strongly $0$-dimensional uniform spaces}\index{uniform
spaces!strongly 0-dimensional@strongly $0$-dimensional} as a uniform space.
This implies that $X$ has topological dimension $0$, as before.
Equivalently, $X$ is strongly $0$-dimensional if there is a base for
the topology of $X$ consisting of subsets of $X$ that are uniformly
separated from their complements in $X$ with respect to $\mathcal{U}$.

        Let us say that $X$ is \emph{uniformly
$0$-dimensional}\index{uniformly 0-dimensional uniform
spaces@uniformly $0$-dimensional uniform spaces}\index{uniform
spaces!uniformly 0-dimensional@uniformly $0$-dimensional} if there
is a base $\mathcal{B}$ for $\mathcal{U}$ such that each $U \in
\mathcal{B}$ corresponds to an equivalence relation on $X$.  It
suffices to have a sub-base for $\mathcal{U}$ with this property,
since the intersection of two subsets of $X \times X$ corresponding to
equivalence relations on $X$ corresponds to an equivalence relation on
$X$ as well.  If $X$ is uniformly $0$-dimensional with respect to
$\mathcal{U}$, then one can check that $X$ is strongly $0$-dimensional
with respect to $\mathcal{U}$.  This uses the fact that if $U \in
\mathcal{U}$ corresponds to an equivalence relation on $X$, then each
equivalence class in $X$ determined by $U$ is $U$-separated from its
complement.  Note that $X$ is uniformly $0$-dimensional with respect
to $\mathcal{U}$ if and only if $\mathcal{U}$ corresponds to a
nonempty collection of semi-ultrametrics on $X$ as in Section
\ref{collections of semimetrics}, by the discussion in Section
\ref{collections of semi-ultrametrics}.

        If $\mathcal{U}_{eq}$ and $\mathcal{U}_{eq, 2}$ are the uniformities
obtained from $\mathcal{U}$ as in the previous section, then $X$ is
automatically uniformly $0$-dimensional with respect to both
$\mathcal{U}_{eq}$ and $\mathcal{U}_{eq, 2}$.  Similarly,
\begin{equation}
\label{mathcal{U} = mathcal{U}_{eq}}
        \mathcal{U} = \mathcal{U}_{eq}
\end{equation}
if and only if $X$ is uniformly $0$-dimensional with respect to
$\mathcal{U}$.  In particular, $X$ is automatically strongly
$0$-dimensional with respect to both $\mathcal{U}_{eq}$ and
$\mathcal{U}_{eq, 2}$, as in the preceding paragraph.  Remember that
the topologies on $X$ associated to $\mathcal{U}_{eq}$ and
$\mathcal{U}_{eq, 2}$ are the same, as in the previous section.  The
condition that $X$ be strongly $0$-dimensional with respect to
$\mathcal{U}$ is equivalent to saying that the topology on $X$
associated to $\mathcal{U}$ is the same as the topology on $X$
associated to $\mathcal{U}_{eq}$ or $\mathcal{U}_{eq, 2}$.

        A topological space $X$ is said to be \emph{totally
separated}\index{totally separated topological spaces} if for every
$x, y \in X$ with $x \ne y$ there is an open set $U \subseteq X$
such that $x \in U$, $y \not\in U$, and $U$ is a closed set in $X$
too.  In particular, this implies that $X$ is Hausdorff.  If $X$ has
topological dimension $0$ and satisfies the first or $0$th separation
condition, then it is easy to see that $X$ is totally separated.

        Let us say that a uniform space $(X, \mathcal{U})$ is
\emph{strongly totally separated}\index{strongly totally separated
uniform spaces}\index{uniform spaces!strongly totally separated}
if for each $x, y \in X$ with $x \ne y$ there is a set $U \subseteq X$
such that $x \in U$, $y \not\in U$, and $U$ is uniformly separated
from $X \setminus U$ in $X$ with respect to $\mathcal{U}$.  This
implies that $U$ is both open and closed in $X$ with respect to the
topology associated to $\mathcal{U}$, as before.  If $X$ is strongly
totally separated with respect to $\mathcal{U}$, then it follows that
$X$ is totally separated with respect to the topology associated to
$\mathcal{U}$.  If $X$ is strongly $0$-dimensional with respect to
$\mathcal{U}$, and if $X$ satisfies the first or $0$th separation
condition with respect to the topology associated to $\mathcal{U}$,
then $X$ is strongly totally separated with respect to $\mathcal{U}$.
Note that $X$ is strongly totally separated with respect to
$\mathcal{U}$ if and only if $X$ is Hausdorff with respect to the
topology associated to $\mathcal{U}_{eq}$ or $\mathcal{U}_{eq, 2}$.

        As a basic example, the set ${\bf Q}$\index{Q@${\bf Q}$} of
rational numbers has topological dimension $0$ with respect to the
standard Euclidean topology.  However, ${\bf Q}$ is also chain
connected with respect to the uniformity determined by the standard
Euclidean metric.  This implies that there are no nonempty proper
subsets of ${\bf Q}$ that are uniformly separated from their
complements in ${\bf Q}$ with respect to this uniformity.  In
particular, this means that ${\bf Q}$ is not strongly $0$-dimensional
at any point with respect to this uniformity.  Similarly, ${\bf Q}$ is
not strongly totally separated with respect to this uniformity.

\part{Topological groups}
\label{topological groups}

\section{Basic notions}
\label{basic notions}
\setcounter{equation}{0}

        Let $G$ be a group, with the group operations expressed
multiplicatively.  Suppose that $G$ is also equipped with a topology,
so that $G \times G$ may be equipped with the corresponding product
topology.  If the group operations on $G$ are continuous, then $G$ is
said to be a \emph{topological group}.\index{topological groups} More
precisely, this means that
\begin{equation}
\label{(x, y) mapsto x y}
        (x, y) \mapsto x \, y
\end{equation}
should be continuous as a mapping from $G \times G$ into $G$.  Similarly,
the mapping
\begin{equation}
\label{x mapsto x^{-1}}
        x \mapsto x^{-1}
\end{equation}
should be continuous on $G$ as well.  This implies that (\ref{x mapsto
  x^{-1}}) should be a homeomorphism on $G$, since this mapping is its
own inverse.  Sometimes one also asks that $\{e\}$ be a closed set in
$G$, where $e$ is the identity element in $G$.

        Put
\begin{equation}
\label{L_a(x) = a x}
        L_a(x) = a \, x
\end{equation}
and
\begin{equation}
\label{R_a(x) = x a}
        R_a(x) = x \, a
\end{equation}
for every $a, x \in G$.  These define the \emph{left}\index{left
  translations} and \emph{right translation}\index{right translations}
mappings on $G$ associated to $a \in G$.  If $G$ is a topological
group, then $L_a$ and $R_a$ define continuous mappings on $G$ for each
$a \in G$, which corresponds to continuity of multiplication on $G$ in
each variable separately.  Note that $L_a$ and $R_a$ are one-to-one
mappings from $G$ onto itself for each $a \in G$, with inverse
mappings $L_{a^{-1}}$ and $R_{a^{-1}}$, respectively.  If $G$ is a
topological group, then it follows that $L_a$ and $R_a$ are
homeomorphisms on $G$ for every $a \in G$.  Of course, if $G$ is
commutative, then $L_a$ is the same as $R_a$ for every $a \in G$.  If
$\{e\}$ is a closed set in $G$, then continuity of left or right
translations implies that $G$ satisfies the first separation
condition.

        If $a \in G$ and $E \subseteq G$, then we put
\begin{equation}
\label{a E = L_a(E) and E a = R_a(E)}
        a \, E = L_a(E) \quad\hbox{and}\quad E \, a = R_a(E).
\end{equation}
Also put
\begin{equation}
\label{E^{-1} = {x^{-1} : x in E}}
        E^{-1} = \{x^{-1} : x \in E\},
\end{equation}
which is the image of $E$ under (\ref{x mapsto x^{-1}}).  If $A, B
\subseteq G$, then we put
\begin{equation}
\label{A B = {a b : a in A, b in B}}
        A \, B = \{a \, b : a \in A, \, b \in B\}.
\end{equation}
Equivalently,
\begin{equation}
\label{A B = bigcup_{a in A} a B = bigcup_{b in B} A b}
        A \, B = \bigcup_{a \in A} a \, B = \bigcup_{b \in B} A \, b.
\end{equation}
If $G$ is a topological group and either $A$ or $B$ is an open set,
then it follows that $A \, B$ is an open set too, by continuity of
translations, and because a union of open sets is an open set.

        Continuity of (\ref{(x, y) mapsto x y}) as a mapping from
$G \times G$ into $G$ at $e$ means that for each open set $W \subseteq G$
with $e \in W$, there are open sets $U, V \subseteq G$ that both
contain $e$ and satisfy
\begin{equation}
\label{U V subseteq W}
        U \, V \subseteq W.
\end{equation}
Suppose that $G$ is equipped with a topology such that the left and
right translation mappings $L_a$ and $R_a$ are continuous for every $a
\in G$, and hence are homeomorphisms.  If (\ref{(x, y) mapsto x y}) is
continuous as a mapping from $G \times G$ into $G$ at $(e, e)$, then
it is easy to see that (\ref{(x, y) mapsto x y}) is continuous
everywhere on $G \times G$.  Similarly, if (\ref{x mapsto x^{-1}}) is
continuous at $e$, then (\ref{x mapsto x^{-1}}) is continuous
everywhere on $G$ under these conditions.

\section{Associated uniformities}
\label{associated uniformities}
\setcounter{equation}{0}

        Let $G$ be a group, and let $A$ be a subset of $G$.  Put
\begin{eqnarray}
\label{A_L = {(x, y) in G times G : x^{-1} y in A} = ...}
 A_L & = & \{(x, y) \in G \times G : x^{-1} \, y \in A\} \\
 & = & \{(x, y) \in G \times G : y = x \, a \hbox{ for some } a \in A\}
                                                 \nonumber
\end{eqnarray}
and
\begin{eqnarray}
\label{A_R = {(x, y) in G times G : y x^{-1} in A} = ...}
 A_R & = & \{(x, y) \in G \times G : y \, x^{-1} \in A\} \\
 & = & \{(x, y) \in G \times G : y = a \, x \hbox{ for some } a \in A\}.
                                                 \nonumber
\end{eqnarray}
Of course, $A_L = A_R$ when $G$ is commutative.  If $\Delta =
\Delta_G$ is the diagonal in $G \times G$, as in (\ref{Delta = Delta_X
  = {(x, x) : x in X}}), then
\begin{equation}
\label{Delta subseteq A_L, A_R}
        \Delta \subseteq A_L, A_R
\end{equation}
exactly when $e \in A$.  We also have that
\begin{equation}
\label{widetilde{(A_L)} = (A^{-1})_L and widetilde{(A_R)} = (A^{-1})_R}
        \widetilde{(A_L)} = (A^{-1})_L \quad\hbox{and}\quad
         \widetilde{(A_R)} = (A^{-1})_R,
\end{equation}
where the left sides of these equations are defined as in
(\ref{widetilde{U} = {(x, y) : (y, x) in U}}).  By construction,
\begin{equation}
\label{A_L[x] = x A and A_R[x] = A x}
        A_L[x] = x \, A \quad\hbox{and}\quad A_R[x] = A \, x
\end{equation}
for every $x \in G$, using the notation in (\ref{U[x] = U[{x}] = {y in
    X : (x, y) in U}}).  Similarly,
\begin{equation}
\label{A_L[E] = E A and A_R[E] = A E}
        A_L[E] = E \, A \quad\hbox{and}\quad A_R[E] = A \, E
\end{equation}
for every $E \subseteq G$, using the notation in (\ref{U[A] = {y in X
    : there is an x in A such that (x, y) in U}}).  If $B$ is another
subset of $G$, then
\begin{equation}
\label{A_L * B_L = (A B)_L and A_R * B_R = (B A)_R}
 A_L * B_L = (A \, B)_L \quad\hbox{and}\quad A_R * B_R = (B \, A)_R,
\end{equation}
where the left sides of these equations are as defined in (\ref{U * V
  = ...}).

        Suppose now that $G$ is a topological group.  One can check
that
\begin{equation}
\label{mathcal{B}_L = {W_L : W subseteq G is an open set with e in W}}
 \mathcal{B}_L = \{W_L : W \subseteq G \hbox{ is an open set with } e \in W\}
\end{equation}
is a base for a uniformity on $G$, which is known as the \emph{left
  uniformity}\index{left uniformity} on $G$.  Similarly,
\begin{equation}
\label{mathcal{B}_R = {W_R : W subseteq G is an open set with e in W}}
 \mathcal{B}_R = \{W_R : W \subseteq G \hbox{ is an open set with } e \in W\}
\end{equation}
is a base for a uniformity on $G$, which is known as the \emph{right
  uniformity}\index{right uniformity} on $G$.  Note that $W_R$ is
defined a bit differently on p210 of \cite{jk}, and corresponds to
$(W^{-1})_R$ here.  This does not affect (\ref{mathcal{B}_R = {W_R : W
    subseteq G is an open set with e in W}}), because $W$ is an open
subset of $G$ that contains $e$ if and only if $W^{-1}$ has the same
properties.

        It is easy to see that the given topology on $G$ is the same
as the topology associated to the corresponding left or right
uniformity on $G$.  In particular, this implies that $G$ is regular as
a topological space in the strict sense, as in Section \ref{closure,
  regularity}.  If $\{e\}$ is a closed set in $G$, then $G$ satisfies
the first separation condition, as in the previous section, and hence
$G$ is Hausdorff as a topological space.

        Suppose that $G_1$ and $G_2$ are topological groups, and that
$\phi$ is a continuous group homomorphism from $G_1$ into $G_2$.
Under these conditions, it is easy to see that $\phi$ is uniformly
continuous as a mapping from $G_1$ into $G_2$, with respect to their
corresponding left uniformities.  Similarly, $\phi$ is uniformly
continuous with respect to the corresponding right uniformities on
$G_1$ and $G_2$.  More precisely, if $\phi$ is continuous at the
identity element in $G_1$, then $\phi$ is continuous, and uniformly
continuous with respect to the left and right uniformities on $G_1$
and $G_2$, respectively.

\section{Some additional properties}
\label{some additional properties}
\setcounter{equation}{0}

        Let $G$ be a group again, and let $L_g$, $R_g$ be the
left and right translation mappings on $G$ corresponding to $g \in G$
as in (\ref{L_a(x) = a x}) and (\ref{R_a(x) = x a}).  This leads to
mappings $L_{g, 2}$ and $R_{g, 2}$ from $G \times G$ into itself as in
(\ref{f_2(x, x') = (f(x), f(x'))}), so that
\begin{equation}
\label{L_{g, 2}(x, y) = (L_g(x), L_g(y)) = (g x, g y)}
        L_{g, 2}(x, y) = (L_g(x), L_g(y)) = (g \, x, g \, y)
\end{equation}
and
\begin{equation}
\label{R_{g, 2}(x, y) = (R_g(x), R_g(y)) = (x g, y g)}
        R_{g, 2}(x, y) = (R_g(x), R_g(y)) = (x \, g, y \, g)
\end{equation}
for every $g, x, y \in G$.  If $A \subseteq G$ and $A_L, A_R \subseteq
G \times G$ are as in (\ref{A_L = {(x, y) in G times G : x^{-1} y in
    A} = ...}), (\ref{A_R = {(x, y) in G times G : y x^{-1} in A} =
  ...}), then it is easy to see that
\begin{equation}
\label{L_{g, 2}(A_L) = A_L and R_{g, 2}(A_R) = A_R}
        L_{g, 2}(A_L) = A_L \quad\hbox{and}\quad R_{g, 2}(A_R) = A_R
\end{equation}
for every $g \in G$.  In particular, if $G$ is a topological group,
then each element of $\mathcal{B}_L$ in (\ref{mathcal{B}_L = {W_L : W
    subseteq G is an open set with e in W}}) is invariant under $L_{g,
  2}$ for every $g \in G$, and each element of $\mathcal{B}_R$ in
(\ref{mathcal{B}_R = {W_R : W subseteq G is an open set with e in W}})
is invariant under $R_{g, 2}$ for every $g \in G$.

        Similarly, one can check that
\begin{equation}
\label{R_{g, 2}(A_L) = (g^{-1} A g)_L and L_{g, 2}(A_R) = (g A g^{-1})_R}
        R_{g, 2}(A_L) = (g^{-1} \, A \, g)_L \quad\hbox{and}\quad
         L_{g, 2}(A_R) = (g \, A \, g^{-1})_R
\end{equation}
for every $g \in G$ and $A \subseteq G$.  If $G$ is a topological
group, then it follows that $R_g$ is uniformly continuous with respect
to the left uniformity on $G$ for every $g \in G$, and that $L_g$ is
uniformly continuous with respect to the right uniformity on $G$ for
every $g \in G$.

        Let $j$ denote the mapping (\ref{x mapsto x^{-1}}) that sends
an element of $G$ to its inverse, so that
\begin{equation}
\label{j(A) = A^{-1}}
        j(A) = A^{-1}
\end{equation}
for each $A \subseteq G$, by construction.  Also let $j_2$ be the
corresponding mapping from $G \times G$ into itself, as in
(\ref{f_2(x, x') = (f(x), f(x'))}), so that
\begin{equation}
\label{j_2(x, y) = (x^{-1}, y^{-1})}
        j_2(x, y) = (x^{-1}, y^{-1})
\end{equation}
for every $x, y \in G$.  Observe that
\begin{equation}
\label{j_2(A_L) = (A^{-1})_R and j_2(A_R) = (A^{-1})_L}
        j_2(A_L) = (A^{-1})_R \quad\hbox{and}\quad j_2(A_R) = (A^{-1})_L
\end{equation}
for every $A \subseteq G$.  In particular, if $G$ is a topological
group, then $j_2$ sends elements of $\mathcal{B}_L$ to elements of
$\mathcal{B}_R$, and vice-versa.  This implies that $j$ is uniformly
continuous as a mapping from $G$ with the left uniformity to $G$ with
the right uniformity, and vice-versa.

        If $g \in G$, then
\begin{equation}
\label{C_g(x) = g x g^{-1}}
        C_g(x) = g \, x \, g^{-1}
\end{equation}
defines an (inner) automorphism on $G$, which is conjugation by $g$.
Of course,
\begin{equation}
\label{C_g(A) = g A g^{-1}}
        C_g(A) = g \, A \, g^{-1}
\end{equation}
for each $A \subseteq G$.  Let $C_{g, 2}$ be the mapping from $G
\times G$ into itself corresponding to $C_g$ as in (\ref{f_2(x, x') =
  (f(x), f(x'))}), so that
\begin{equation}
\label{C_{g, 2}(x, y) = (C_g(x), C_g(y)) = (g x g^{-1}, g y g^{-1})}
        C_{g, 2}(x, y) = (C_g(x), C_g(y)) = (g \, x \, g^{-1}, g \, y \, g^{-1})
\end{equation}
for every $x, y \in G$.    Observe that
\begin{equation}
\label{C_{g, 2}(A_L) = (C_g(A))_L and C_{g, 2}(A_R) = (C_g(A))_R}
 C_{g, 2}(A_L) = (C_g(A))_L \quad\hbox{and}\quad C_{g, 2}(A_R) = (C_g(A))_R
\end{equation}
for every $A \subseteq G$, which can also be derived from (\ref{L_{g,
    2}(A_L) = A_L and R_{g, 2}(A_R) = A_R}) and (\ref{R_{g, 2}(A_L) =
  (g^{-1} A g)_L and L_{g, 2}(A_R) = (g A g^{-1})_R}).  If $G$ is a
topological group, then it follows that $C_g$ is uniformly continuous
with respect to the left and right uniformities on $G$, which could be
obtained from the earlier statements for left and right translations
as well.

\section{Translation-invariant relations}
\label{translation-invariant relations}
\setcounter{equation}{0}

        Let $G$ be a group, and let $U \subseteq G \times G$ be given.
Let us say that $U$ is invariant under left translations on $G$ if
\begin{equation}
\label{L_{g, 2}(U) = U}
        L_{g, 2}(U) = U
\end{equation}
for every $g \in G$, where $L_{g, 2}$ is as in (\ref{L_{g, 2}(x, y) =
  (L_g(x), L_g(y)) = (g x, g y)}).  Similarly, $U$ is invariant under
right translations on $G$ if
\begin{equation}
\label{R_{g, 2}(U) = U}
        R_{g, 2}(U) = U
\end{equation}
for every $g \in G$, where $R_{g, 2}$ is as in (\ref{R_{g, 2}(x, y) =
  (R_g(x), R_g(y)) = (x g, y g)}).  Put
\begin{equation}
\label{A = U[e]}
        A = U[e],
\end{equation}
where the right side of (\ref{A = U[e]}) is defined as in (\ref{U[x] =
  U[{x}] = {y in X : (x, y) in U}}).  If $U$ is invariant under left
translations on $G$, then it is easy to see that
\begin{equation}
\label{U = A_L}
        U = A_L,
\end{equation}
where $A_L$ is as in (\ref{A_L = {(x, y) in G times G : x^{-1} y in A}
  = ...}).  In the same way, if $U$ is invariant under right
translations on $G$, then
\begin{equation}
\label{U = A_R}
        U = A_R,
\end{equation}
where $A_R$ is as in (\ref{A_R = {(x, y) in G times G : y x^{-1} in A}
  = ...}).  Of course, for any $A \subseteq G$, we have seen that
$A_L$ and $A_R$ are invariant under left and right translations on
$G$, respectively, as in (\ref{L_{g, 2}(x, y) = (L_g(x), L_g(y)) = (g
  x, g y)}) and (\ref{R_{g, 2}(x, y) = (R_g(x), R_g(y)) = (x g, y
  g)}).

        Let us say that $A \subseteq G$ is invariant under conjugation if
\begin{equation}
\label{C_g(A) = A}
        C_g(A) = A
\end{equation}
for every $g \in G$, where $C_g$ is as in (\ref{C_g(x) = g x g^{-1}}).
Similarly, we say that $U \subseteq G \times G$ is invariant under
conjugations on $G$ when
\begin{equation}
\label{C_{g, 2}(U) = U}
        C_{g, 2}(U) = U
\end{equation}
for every $g \in G$, where $C_{g, 2}$ is as in (\ref{C_{g, 2}(x, y) =
  (C_g(x), C_g(y)) = (g x g^{-1}, g y g^{-1})}).  If $A \subseteq G$
is invariant under conjugation in the first sense, then $A_L$ and
$A_R$ are invariant under conjugation in the second sense, by
(\ref{C_{g, 2}(A_L) = (C_g(A))_L and C_{g, 2}(A_R) = (C_g(A))_R}).
This is basically the same as saying that $A_L$ and $A_R$ are
invariant under right and left translations, respectively, as in
(\ref{R_{g, 2}(A_L) = (g^{-1} A g)_L and L_{g, 2}(A_R) = (g A
  g^{-1})_R}).  In this case, we also have that
\begin{equation}
\label{A_L = A_R}
        A_L = A_R.
\end{equation}
If $U \subseteq G \times G$ is invariant under conjugations and $A$ is
as in (\ref{A = U[e]}), then $A$ is invariant under conjugation as a
subset of $G$.  If $U$ is invariant under both left and right
translations, then $U$ is invariant under conjugations in particular.
If $U$ is invariant under either left or right translations, and if
$U$ is invariant under conjugations, then $U$ is invariant under both
left and right translations.

        If $A$ is a subgroup of $G$, then $A_L$ and $A_R$ correspond
to equivalence relations on $G$.  The equivalence classes associated
to these equivalence relations are the usual left and right cosets of
$A$ in $G$, respectively.  Remember that $A$ is said to be normal in
$G$ exactly when $A$ is invariant under conjugation, in which case the
left and right cosets of $A$ in $G$ are the same.  This implies that
(\ref{A_L = A_R}) holds, and the corresponding equivalence relation on
$G$ is invariant under both left and right translations, as in the
previous paragraph.

        Suppose now that $U \subseteq G \times G$ is invariant
under left or right translations, so that $U$ can be given as in
(\ref{U = A_L}) or (\ref{U = A_R}), respectively, with $A$ as in
(\ref{A = U[e]}).  In both cases, if $U$ also corresponds to an
equivalence relation on $G$, then one can check that $A$ is a subgroup
of $G$.  If $U$ is invariant under both left and right translations,
and hence under conjugations, then $A$ is a normal subgroup of $G$.

\section{Translation-invariant semimetrics}
\label{translation-invariant semimetrics}
\setcounter{equation}{0}

        Let $G$ be a group, and let $d(x, y)$ be a $q$-semimetric
on $G$ for some $q > 0$.  We say that $d(\cdot, \cdot)$ is invariant
under left translations on $G$ if
\begin{equation}
\label{d(a x, a y) = d(x, y)}
        d(a \, x, a \, y) = d(x, y)
\end{equation}
for every $a, x, y \in G$.  Similarly, $d(\cdot, \cdot)$ is invariant
under right translations on $G$ if
\begin{equation}
\label{d(x a, y a) = d(x, y)}
        d(x \, a, y \, a) = d(x, y)
\end{equation}
for every $a, x, y \in G$.  In both cases, we get that
\begin{equation}
\label{d(x, e) = d(x^{-1}, e)}
        d(x, e) = d(x^{-1}, e)
\end{equation}
for every $x \in G$, using the symmetry condition (\ref{d(x, y) = d(y,
  x)}).  Note that $d(x, y)$ is invariant under left translations on
$G$ if and only if
\begin{equation}
\label{d(x^{-1}, y^{-1})}
        d(x^{-1}, y^{-1})
\end{equation}
is invariant under right translations on $G$.

        If
\begin{equation}
\label{d(a x a^{-1}, a y a^{-1}) = d(x, y)}
        d(a \, x \, a^{-1}, a \, y \, a^{-1}) = d(x, y)
\end{equation}
for every $a, x, y \in G$, then $d(\cdot, \cdot)$ is said to be
invariant under conjugations on $G$.  In particular, this implies that
\begin{equation}
\label{d(a x a^{-1}, e) = d(x, e)}
        d(a \, x \, a^{-1}, e) = d(x, e)
\end{equation}
for every $a, x \in G$.  If $d(\cdot, \cdot)$ is invariant under left
or right translations on $G$ and satisfies (\ref{d(a x a^{-1}, e) =
  d(x, e)}), then one can check that $d(\cdot, \cdot)$ is invariant
under conjugations, and in fact under both left and right translations
on $G$.  Of course, if $d(\cdot, \cdot)$ is invariant under both left
and right translations on $G$, then $d(\cdot, \cdot)$ is invariant
under conjugations on $G$.  In this case, one can verify that
\begin{equation}
\label{d(x^{-1}, y^{-1}) = d(x, y)}
        d(x^{-1}, y^{-1}) = d(x, y)
\end{equation}
for every $x, y \in G$ too.

        Let $r > 0$ be given, and let $U_d(r) \subseteq G \times G$
be as in (\ref{U(r) = U_d(r) = {(x, y) in X times X : d(x, y) < r}}).
If $d(\cdot, \cdot)$ is invariant under left or right translations on
$G$, then $U_d(r)$ has the analogous property, as defined in the
preceding section.  Similarly, if $d(\cdot, \cdot)$ is invariant under
conjugations on $G$, then $U_d(r)$ is invariant under conjugations as
well.  This implies that the open ball $B_d(e, r)$ defined as in
(\ref{B(x, r) = B_d(x, r) = {y in X : d(x, y) < r}}) is invariant
under conjugations as a subset of $G$, since
\begin{equation}
\label{B_d(e, r) = (U_d(r))[e]}
        B_d(e, r) = (U_d(r))[e],
\end{equation}
as in (\ref{A = U[e]}).  More precisely, $B_d(e, r)$ is invariant
under conjugations when (\ref{d(a x a^{-1}, e) = d(x, e)}) holds.  If
$d(\cdot, \cdot)$ is a semi-ultrametric on $G$, then $U_d(r)$
corresponds to an equivalence relation on $G$, as in Section
\ref{semi-ultrametrics}.  If $d(\cdot, \cdot)$ is a semi-ultrametric
on $G$ that is invariant under left or right translations on $G$, then
it follows that $B_d(e, r)$ is a subgroup of $G$, as in the previous
section.  If $d(\cdot, \cdot)$ is a semi-ultrametric on $G$ that is
invariant under both left and right translations on $G$, then $B_d(e,
r)$ is a normal subgroup of $G$.

        If $A$ is a subgroup of $G$, then $A_L$ and $A_R$ in
(\ref{A_L = {(x, y) in G times G : x^{-1} y in A} = ...}) and
(\ref{A_R = {(x, y) in G times G : y x^{-1} in A} = ...})
correspond to equivalence relations on $G$, as in the previous
section.  These equivalence relations are also invariant under left
and right translations, respectively.  It is easy to see that the
corresponding discrete semi-ultrametrics on $G$ as in (\ref{d(x, y) =
  0 when x sim y, = 1 when x not sim y}) are invariant under left and
right translations as well.  If $A$ is a normal subgroup of $G$, then
$A_L = A_R$, and the corresponding discrete semi-ultrametric is
invariant under both left and right translations on $G$.

\section{Translation-invariance and topology}
\label{translation-invariance, topology}
\setcounter{equation}{0}

        Suppose now that $G$ is a topological group.  Let $A$
be an open subset of $G$ that contains $e$, and let $A_L$, $A_R$
be as in (\ref{A_L = {(x, y) in G times G : x^{-1} y in A} = ...}),
(\ref{A_R = {(x, y) in G times G : y x^{-1} in A} = ...}), respectively.
Thus $A_L$, $A_R$ are elements of $\mathcal{B}_L$, $\mathcal{B}_R$
in (\ref{mathcal{B}_L = {W_L : W subseteq G is an open set with e in W}}),
(\ref{mathcal{B}_R = {W_R : W subseteq G is an open set with e in W}}),
respectively.  In particular, $A_L$ and $A_R$ are elements of the
corresponding left and right uniformities on $G$, respectively.
If $A$ is a subgroup of $G$, then $A_L$, $A_R$ correspond to
equivalence relations on $G$, as before.  Note that a subgroup $A$
of $G$ is an open subset of $G$ when $e$ is an element of the
interior of $A$, by continuity of translations.

        Conversely, if $U \subseteq G \times G$ is invariant under
left or right translations, then we have seen that $U$ can be
expressed as $A_L$ or $A_R$, respectively, where $A$ is as in (\ref{A
  = U[e]}).  If $U$ corresponds to an equivalence relation on $G$,
then $A$ is a subgroup of $G$.  If $U$ is an element of the left or
right uniformity on $G$, then it is easy to see that $e$ is an element
of the interior of (\ref{A = U[e]}) in $G$.  This uses the way that
the left and right uniformities on $G$ are defined in terms of the
given topology on $G$.  If $U$ has each of these three properties,
then it follows that $A$ is an open subgroup of $G$.

        Let $\mathcal{A}$ be a collection of open subsets of $G$
that contain $e$ as an element.  Put
\begin{equation}
\label{mathcal{A}_L = {A_L : A in mathcal{A}}}
        \mathcal{A}_L = \{A_L : A \in \mathcal{A}\}
\end{equation}
and
\begin{equation}
\label{mathcal{A}_R = {A_R : A in mathcal{A}}}
        \mathcal{A}_R = \{A_R : A \in \mathcal{A}\},
\end{equation}
so that $\mathcal{A}_L \subseteq \mathcal{B}_L$ and $\mathcal{A}_R
\subseteq \mathcal{B}_R$, by definition of $\mathcal{B}_L$ and
$\mathcal{B}_R$.  If $\mathcal{A}$ is the collection of all open
subsets of $G$ that contain $e$, then $\mathcal{A}_L = \mathcal{B}_L$
and $\mathcal{A}_R = \mathcal{B}_R$.  Similarly, if $\mathcal{A}$ is a
local base for the topology of $G$ at $e$, then $\mathcal{A}_L$ and
$\mathcal{A}_R$ form bases for the left and right uniformities on $G$,
respectively.  Conversely, if $\mathcal{A}_L$ or $\mathcal{A}_R$ is a
base for the left or right uniformity on $G$, respectively, then it is
easy to see that $\mathcal{A}$ has to be a local base for the topology
of $G$ at $e$.

        If $\mathcal{A}$ is a local base for the topology of $e$
consisting of open subgroups of $G$, then $\mathcal{A}_L$ and
$\mathcal{A}_R$ are bases for the left and right uniformities on $G$,
respectively, whose elements correspond to equivalence relations on
$G$.  This implies that $G$ is uniformly $0$-dimensional with respect
to the left and right uniformities, as in Section \ref{topological
  dimension 0}.  The converse will be discussed in the next section.

        It is well known that there is a collection of semimetrics on $G$
that are invariant under left translations on $G$, and for which the
corresponding topology on $G$ is the given topology.  This implies
that the uniformity on $G$ determined by this collection of
semimetrics is the same as the left uniformity on $G$, because of
invariance under left translations.  If there is a local base for the
topology of $G$ at $e$ with only finitely or countably many elements,
then there is a single semimetric on $G$ that is invariant under left
translations and determines the same topology on $G$, and for which
the corresponding uniformity is hence the left uniformity.  Of course,
there are analogous statements for right translations and the right
uniformity.

        Suppose that the topology on $G$ is determined by a
nonempty collection $\mathcal{M}$ of semi-ultrametrics on $G$ that are
invariant under left or right translations.  Of course, the open
balls in $G$ centered at $e$ with respect to elements of $\mathcal{M}$
are open subsets of $G$ under these conditions, and they are also
subgroups of $G$, as in the previous section.  This implies that
the open subgroups of $G$ form a local base for the topology of
$G$ at $e$.  If the elements of $\mathcal{M}$ are invariant under
both left and right translations, then the corresponding open balls
in $G$ centered at $e$ are normal subgroups of $G$.  In this case,
the open normal subgroups of $G$ form a local base for the topology
of $G$ at $e$.

        Conversely, suppose that $\mathcal{A}$ is a collection of
open subgroups of $G$ that is a local base for the topology of $G$ at
$e$.  If $A \in \mathcal{A}$, then we get subsets $A_L$, $A_R$ of $G
\times G$ corresponding to equivalence relations on $G$ as before.
The discrete semi-ultrametrics on $G$ associated to these equivalence
relations are invariant under left and right translations,
respectively, as in the previous section.  By doing this for each $A
\in \mathcal{A}$, we get collections of semi-ultrametrics on $G$ that
are invariant under left or right translations and which determine the
same topology on $G$.  If the elements of $\mathcal{A}$ are normal
subgroups of $G$, then the corresponding discrete semi-ultrametrics on
$G$ are invariant under both left and right translations, as before.

        If $d_1, d_2, \ldots, d_n$ are finitely many semimetrics on $G$,
each of which is invariant under left translations, then their sum and
maximum are invariant under left translations too.  If $d$ is a
semimetric on $G$ that is invariant under left translations and $t$ is
a positive real number, then the semimetric $d_t$ defined on $G$ as in
(\ref{d_t(x, y) = min(d(x, y), t)}) is invariant under left
translations as well.  Similarly, if $d_1, d_2, d_3, \ldots$ is an
infinite sequence of semimetrics on $G$, each of which is invariant
under left translations, then the semimetric $d'$ defined on $G$ as in
(\ref{d'(x, y) = max_{j ge 1} d_j'(x, y)}) is also invariant under
left translations.  There are analogous statements for right translations,
as usual.

\section{Translation-invariant relations, continued}
\label{translation-invariant relations, continued}
\setcounter{equation}{0}

        Let $G$ be a group, and let $U$ be a subset of $G \times G$.
Observe that
\begin{equation}
\label{bigcap_{g in G} L_{g, 2}(U) = ...}
  \quad \bigcap_{g \in G} L_{g, 2}(U)
              = \{(x, y) \in G \times G : (g^{-1} \, x, g^{-1} \, y) \in U
                                               \hbox{ for every } g \in G\} 
\end{equation}
where $L_{g, 2}$ is as in (\ref{L_{g, 2}(x, y) = (L_g(x), L_g(y)) = (g
  x, g y)}) for each $g \in G$.  Equivalently,
\begin{equation}
\label{bigcap_{g in G} L_{g, 2}(U) = ..., 2}
  \qquad \bigcap_{g \in G} L_{g, 2}(U)
              = \{(x, y) \in G \times G : (h^{-1}, h^{-1} \, x^{-1} \, y) \in U
                                                 \hbox{ for every } h \in G\},
\end{equation}
where $h$ corresponds to $x^{-1} \, g$ in the right side of
(\ref{bigcap_{g in G} L_{g, 2}(U) = ...}).  This implies that
\begin{equation}
\label{bigcap_{g in G} L_{g, 2}(U) = (bigcap_{h in G} h U[h^{-1}])_L}
 \bigcap_{g \in G} L_{g, 2}(U) = \Big(\bigcap_{h \in G} h \, U[h^{-1}]\Big)_L,
\end{equation}
where the right side of (\ref{bigcap_{g in G} L_{g, 2}(U) = (bigcap_{h
    in G} h U[h^{-1}])_L}) is defined as in (\ref{A_L = {(x, y) in G
    times G : x^{-1} y in A} = ...}), as in Section
\ref{translation-invariant relations}.  Of course, the left side of
(\ref{bigcap_{g in G} L_{g, 2}(U) = (bigcap_{h in G} h U[h^{-1}])_L})
is invariant under left translations on $G$ as a subset of $G \times
G$ by construction.  One can also check directly that
\begin{equation}
\label{(bigcap_{g in G} L_{g, 2}(U))[e] = bigcap_{g in G} g U[g^{-1}]}
 \Big(\bigcap_{g \in G} L_{g, 2}(U)\Big)[e] = \bigcap_{g \in G} g \, U[g^{-1}],
\end{equation}
where the definition (\ref{U[x] = U[{x}] = {y in X : (x, y) in U}}) is
used on both sides of the equation.  Thus (\ref{bigcap_{g in G} L_{g,
    2}(U) = (bigcap_{h in G} h U[h^{-1}])_L}) is the same as (\ref{U =
  A_L}) in this situation.

        Suppose now that $G$ is a topological group, and that $U$ is an
element of the corresponding left uniformity on $G$, as in
Section \ref{associated uniformities}.  This means that there is an
open set $W \subseteq G$ such that $e \in W$ and
\begin{equation}
\label{W_L subseteq U}
        W_L \subseteq U,
\end{equation}
where $W_L$ is as in (\ref{A_L = {(x, y) in G times G : x^{-1} y in A}
  = ...}) again.  It follows that
\begin{equation}
\label{W_L subseteq bigcap_{g in G} L_{g, 2}(U)}
        W_L \subseteq \bigcap_{g \in G} L_{g, 2}(U),
\end{equation}
because $W_L$ is automatically invariant under left translations, as
in (\ref{L_{g, 2}(A_L) = A_L and R_{g, 2}(A_R) = A_R}).  This is
basically the same as saying that
\begin{equation}
\label{W subseteq bigcap_{g in G} g U[g^{-1}]}
        W \subseteq \bigcap_{g \in G} g \, U[g^{-1}],
\end{equation}
because of (\ref{bigcap_{g in G} L_{g, 2}(U) = (bigcap_{h in G} h
  U[h^{-1}])_L}).  Note that (\ref{W_L subseteq bigcap_{g in G} L_{g,
    2}(U)}) implies that (\ref{bigcap_{g in G} L_{g, 2}(U) = ...}) is
an element of the left uniformity on $G$ too.

        If $U$ corresponds to an equivalence relation on $G$, then
it is easy to see that (\ref{bigcap_{g in G} L_{g, 2}(U) = ...})
corresponds to an equivalence relation on $G$ as well.  This implies
that (\ref{(bigcap_{g in G} L_{g, 2}(U))[e] = bigcap_{g in G} g
  U[g^{-1}]}) is a subgroup of $G$, as in Section
\ref{translation-invariant relations}, because (\ref{bigcap_{g in G}
  L_{g, 2}(U) = ...}) is invariant under left translations.  If $G$ is
a topological group, and if $U$ is an element of the left uniformity,
then $e$ is an element of the interior of this subgroup, by (\ref{W
  subseteq bigcap_{g in G} g U[g^{-1}]}).  It follows that
(\ref{(bigcap_{g in G} L_{g, 2}(U))[e] = bigcap_{g in G} g U[g^{-1}]})
is an open subgroup of $G$ under these conditions.  Note that
(\ref{bigcap_{g in G} L_{g, 2}(U) = ...}) is automatically contained in $U$.

        To summarize, if $U \subseteq G \times G$ is an element of the
left uniformity on $G$, and if $U$ corresponds to an equivalence
relation on $G$, then there is an open subgroup $A$ of $L$ such that
\begin{equation}
\label{A_L subseteq U}
        A_L \subseteq U.
\end{equation}
More precisely, one can take $A$ to be (\ref{(bigcap_{g in G} L_{g,
    2}(U))[e] = bigcap_{g in G} g U[g^{-1}]}), so that $A_L$ is as in
(\ref{bigcap_{g in G} L_{g, 2}(U) = (bigcap_{h in G} h U[h^{-1}])_L}).
If $\mathcal{A}$ is the collection of open subgroups of $G$ and
$\mathcal{A}_L$ is as in (\ref{mathcal{A}_L = {A_L : A in
    mathcal{A}}}), then it follows that $\mathcal{A}_L$ is a base for
the uniformity on $G$ associated to the left uniformity on $G$ as in
(\ref{mathcal{U}_{eq}}).  In particular, if $G$ is uniformly
$0$-dimensional with respect to the left uniformity, as in Section
\ref{topological dimension 0}, then it follows that there is a local
base for the topology of $G$ at $e$ consisting of open subgroups of
$G$.  Of course, one can deal with right translations and the right
uniformity on $G$ in the same way.

\section{Uniform separation}
\label{uniform separation}
\setcounter{equation}{0}

        Let $G$ be a group, let $W$ be a subset of $G$ that contains $e$,
and let $W_L$, $W_R$ be defined as in (\ref{A_L = {(x, y) in G times G
    : x^{-1} y in A} = ...}), (\ref{A_R = {(x, y) in G times G : y
    x^{-1} in A} = ...}), as usual.  Note that $W_L$ and $W_R$ contain
the diagonal in $G \times G$ under these conditions, as in
(\ref{Delta subseteq A_L, A_R}).  Remember that
\begin{equation}
\label{W_L[E] = E W and W_R[E] = W E}
        W_L[E] = E \, W \quad\hbox{and}\quad W_R[E] = W \, E
\end{equation}
for every $E \subseteq G$, as in (\ref{A_L[E] = E A and A_R[E] = A
  E}), and using the notation in (\ref{U[A] = {y in X : there is an x
    in A such that (x, y) in U}}).  It follows that $E_1, E_2
\subseteq G$ are $W_L$-separated in the terminology of Section
\ref{U-separated sets} exactly when
\begin{equation}
\label{(E_1 W) cap E_2 = emptyset}
        (E_1 \, W) \cap E_2 = \emptyset,
\end{equation}
as in (\ref{U[A] cap B = emptyset}).  Similarly, $E_1$, $E_2$ are
$W_R$-separated exactly when
\begin{equation}
\label{(W E_1) cap E_2 = emptyset}
        (W \, E_1) \cap E_2 = \emptyset.
\end{equation}

        Suppose now that $G$ is a topological group, so that uniform
separation of subsets of $G$ with respect to the left and right
uniformities can be defined as in Section \ref{uniformly separated sets}.
More precisely, $E_1, E_2 \subseteq G$ are uniformly separated with
respect to the left uniformity exactly when (\ref{(E_1 W) cap E_2 = emptyset})
holds for some open set $W \subseteq G$ that contains $e$.  Similarly,
$E_1$, $E_2$ are uniformly separated with respect to the right uniformity
on $G$ exactly when (\ref{(W E_1) cap E_2 = emptyset}) holds for some
open set $W \subseteq G$ that contains $e$.

        Let $E$ be a subset of $G$, and let us apply the previous remarks
to $E$ and its complement in $G$.  It follows that $E$ is uniformly
separated from its complement with respect to the left uniformity on $G$
if and only if
\begin{equation}
\label{E W subseteq E}
        E \, W \subseteq E
\end{equation}
for some open set $W \subseteq G$ that contains $e$.  Similarly, $E$
is uniformly separated from its complement with respect to the right
uniformity on $G$ if and only if
\begin{equation}
\label{W E subseteq E}
        W \, E \subseteq E
\end{equation}
for some open set $W \subseteq G$ that contains $e$.  In both cases,
the inclusion is actually an equality, because $e \in W$.  We may also
ask that $W$ be symmetric about $e$, in the sense that
\begin{equation}
\label{W^{-1} = W}
        W^{-1} = W,
\end{equation}
since otherwise we can replace $W$ with $W \cap W^{-1}$.

        If $W$ is any subset of $G$, then let $W^n$ be the subset of $G$
consisting of products of $n$ elements of $W$ for each positive
integer $n$, so that
\begin{equation}
\label{W^{n + 1} = W^n W = W W^n}
        W^{n + 1} = W^n \, W = W \, W^n
\end{equation}
for every $n$.  If $e \in W$ and $W$ satisfies (\ref{W^{-1} = W}),
then it is easy to see that
\begin{equation}
\label{widehat{W} = bigcup_{n = 1}^infty W^n}
        \widehat{W} = \bigcup_{n = 1}^\infty W^n
\end{equation}
is a subgroup of $G$.  If $W$ is also an open set in $G$, then
$\widehat{W}$ is an open subgroup of $G$.  If $W$ satisfies (\ref{E W
  subseteq E}) for some $E \subseteq G$, then we have that
\begin{equation}
\label{E W^n subseteq E}
        E \, W^n \subseteq E
\end{equation}
for every $n \ge 1$, and hence
\begin{equation}
\label{E widehat{W} subseteq E}
        E \, \widehat{W} \subseteq E.
\end{equation}
Similarly, (\ref{W E subseteq E}) implies that
\begin{equation}
\label{widehat{W} E subseteq E}
        \widehat{W} \, E \subseteq E.
\end{equation}

        This shows that if $E$ is uniformly separated from its complement
with respect to the left uniformity on $G$, then $E$ is
$A_L$-separated from its complement for some open subgroup $A$ of $G$,
where $A_L$ is as in (\ref{A_L = {(x, y) in G times G : x^{-1} y in A}
  = ...}).  Similarly, if $E$ is uniformly separated from its
complement with respect to the right uniformity on $G$, then $E$ is
$A_R$-separated from its complement for some open subgroup $A$ of $G$,
where $A_R$ is as in (\ref{A_R = {(x, y) in G times G : y x^{-1} in A}
  = ...}).  In both cases, if we also have that $e \in E$, then we get
that $A \subseteq E$.  If $G$ is strongly $0$-dimensional at $e$ with
respect to the left or right uniformity, as in Section
\ref{topological dimension 0}, then it follows that there is a local
base for the topology of $G$ at $e$ consisting of open subgroups.

        If $A$ is any open subgroup of $G$, then it is easy to see
that $A$ is $A_L$ and $A_R$-separated from its complement in $G$.
If there is a local base for the topology of $G$ at $e$ consisting
of open subgroups, then it follows that $G$ is strongly $0$-dimensional
at $e$ with respect to the left and right uniformities.  This implies
that $G$ is strongly $0$-dimensional at every point with respect to the
left and right uniformities, using invariance under translations.
More precisely, $G$ is uniformly $0$-dimensional with respect to the
left and right uniformities under these conditions, as in Section
\ref{translation-invariance, topology}.

        Alternatively, if $E$ is a subset of $G$, then  let $U_E$
be the subset of $G \times G$ that corresponds to the equivalence
relation whose equivalence classes are $E$ and its complement.
As in Section \ref{equivalence relations}, $E$ is uniformly separated
from its complement with respect to the left uniformity on $G$ if and
only if $U_E$ is an element of the left uniformity on $G$.  In this
case, the discussion in the previous section implies that there is
an open subgroup $A$ of $G$ such that $A_L$ is contained in $U_E$.
It follows that $E$ is $A_L$-separated from its complement, as before.
Similarly, if $E$ is uniformly separated from its complement with
respect to the right uniformity, then one can use an analogous argument
to get an open subgroup $A$ of $G$ such that $E$ is $A_R$-separated
from its complement.

\section{Open subgroups}
\label{open subgroups}
\setcounter{equation}{0}

        Let $G$ be a group, and let $\mathcal{A}$ be a nonempty
collection of subgroups of $G$.  Suppose that $\mathcal{A}$ is
compatible with conjugations on $G$, in the sense that for each $A \in
\mathcal{A}$ and $g \in G$ there is a $B \in \mathcal{A}$ such that
\begin{equation}
\label{B subseteq g A g^{-1}}
        B \subseteq g \, A \, g^{-1}.
\end{equation}
In particular, this holds when $\mathcal{A}$ is invariant under
conjugations on $G$, in the sense that
\begin{equation}
\label{g A g^{-1} in mathcal{A}}
        g \, A \, g^{-1} \in \mathcal{A}
\end{equation}
for every $A \in \mathcal{A}$ and $g \in G$.  Of course, if the
elements of $\mathcal{A}$ are normal subgroups of $G$, then this
condition holds automatically.  If $G$ is a topological group, then
the collection of all open subgroups of $G$ has this property,
because of continuity of translations.

        Let $\tau_\mathcal{A}$ be the collection of subsets $W$ of $G$
such that for each $x \in W$ there are finitely many elements
$A_1, \ldots, A_n$ of $\mathcal{A}$ with the property that
\begin{equation}
\label{x (bigcap_{j = 1}^n A_j) subseteq W}
        x \, \Big(\bigcap_{j = 1}^n A_j\Big) \subseteq W.
\end{equation}
This is equivalent to saying that for each $x \in W$ there are
finitely many elements $A_1, \ldots, A_n$ of $\mathcal{A}$ such that
\begin{equation}
\label{(bigcap_{j = 1}^n A_j) x subseteq W}
        \Big(\bigcap_{j = 1}^n A_j\Big) \, x \subseteq W,
\end{equation}
because $\mathcal{A}$ is supposed to be compatible with conjugations,
as in the preceding paragraph.  Let us say that $\mathcal{A}$ behaves
well with respect to finite intersections if the intersection of
finitely many elements of $\mathcal{A}$ always contains another
element of $\mathcal{A}$ as a subset.  In this case, it suffices to
take $n = 1$ in (\ref{x (bigcap_{j = 1}^n A_j) subseteq W}) and
(\ref{(bigcap_{j = 1}^n A_j) x subseteq W}).  If $G$ is a topological
group, then the collection of all open subgroups of $G$ is closed
under finite intersections.

        It is easy to see that $\tau_\mathcal{A}$ defines a topology
on $G$.  Note that
\begin{equation}
\label{mathcal{A} subseteq tau_mathcal{A}}
        \mathcal{A} \subseteq \tau_\mathcal{A},
\end{equation}
because the elements of $\mathcal{A}$ are subgroups of $G$.  By
construction, $\mathcal{A}$ is a local sub-base for $\tau_\mathcal{A}$
at $e$.  If $\mathcal{A}$ behaves well with respect to finite
intersections, then $\mathcal{A}$ is a local base for
$\tau_\mathcal{A}$ at $e$.  One can also check that $G$ is a
topological group with respect to $\tau_\mathcal{A}$.  More precisely,
continuity of translations on $G$ with respect to $\tau_\mathcal{A}$
is built into the definition of $\tau_\mathcal{A}$.  Similarly,
continuity of the group operations on $G$ at $e$ with respect to
$\tau_\mathcal{A}$ follows from the fact that the elements of
$\mathcal{A}$ are subgroups of $G$.

        If $G$ satisfies the first or $0$th separation condition
with respect to $\tau_\mathcal{A}$, then
\begin{equation}
\label{bigcap_{A in mathcal{A}} A = {e}}
        \bigcap_{A \in \mathcal{A}} A = \{e\},
\end{equation}
which means that for each $x \in G$ with $x \ne e$ there is an $A \in
\mathcal{A}$ such that $x \not\in A$.  Conversely, (\ref{bigcap_{A in
    mathcal{A}} A = {e}}) implies that $G$ is Hausdorff with respect
to $\tau_\mathcal{A}$.  As in Section \ref{translation-invariance,
  topology}, $G$ is automatically uniformly $0$-dimensional with
respect to the left and right uniformities associated to
$\tau_\mathcal{A}$.  If (\ref{bigcap_{A in mathcal{A}} A = {e}})
holds, then it follows that $G$ is uniformly totally separated with
respect to the left and right uniformities associated to
$\tau_\mathcal{A}$, as in Section \ref{topological dimension 0}.  One
can also check this more directly, using the fact that each subgroup
$A$ of $G$ is $A_L$ and $A_R$-separated from its complement, as in the
previous section.  This implies that $A$ is uniformly separated from
its complement with respect to the left and right uniformities
associated to $\tau_\mathcal{A}$ when $A \in \tau_\mathcal{A}$.  Of
course, any uniform space which is strongly totally separated is
totally separated as a topological space, and hence Hausdorff as a
topological space.

        Suppose that $G$ is already a topological group with respect
to some topology $\tau$, and let us take $\mathcal{A}$ to be the
collection of all open subgroups of $G$ with respect to $\tau$.  Thus
$\mathcal{A}$ is invariant under conjugations and closed under finite
intersections, as before.  In this situation,
\begin{equation}
\label{tau_mathcal{A} subseteq tau}
        \tau_\mathcal{A} \subseteq \tau,
\end{equation}
because $\mathcal{A}$ is contained in $\tau$ by hypothesis.  In
particular, this implies that open subgroups of $G$ with respect to
$\tau_\mathcal{A}$ are open with respect to $\tau$.  Conversely, every
open subgroup of $G$ with respect to $\tau$ is in $\mathcal{A}$ and
hence in $\tau_\mathcal{A}$, by construction.

        Suppose for the moment that $G$ is strongly totally separated
with respect to the left or right uniformity associated to $\tau$,
as in Section \ref{topological dimension 0}.  This implies that for
each $x \in G$ with $x \ne e$ there is a subset $E$ of $G$ such that
$e \in E$, $x \not\in E$, and $E$ is strongly separated from its
complement with respect to the left or right uniformity on $G$
associated to $\tau$, as appropriate.  As in the previous section,
it follows that there is an open subgroup $A$ of $G$ with respect
to $\tau$ such that $E$ is $A_L$ or $A_R$-separated from its complement,
as appropriate.  Note that $A \subseteq E$ under these conditions,
because $e \in E$, as before.  Thus $x \not\in A$, while $e \in A$
automatically, so that (\ref{bigcap_{A in mathcal{A}} A = {e}}) holds
with this choice of $\mathcal{A}$.

        Conversely, if (\ref{bigcap_{A in mathcal{A}} A = {e}}) holds
with this choice of $\mathcal{A}$, then we have already seen that $G$
is strongly totally separated with respect to the left and right
uniformities associated to $\tau_\mathcal{A}$.  This implies that $G$
is strongly totally separated with respect to the left and right
uniformities associated to $\tau$, because of (\ref{tau_mathcal{A}
  subseteq tau}).

\section{Commutative groups}
\label{commutative groups}
\setcounter{equation}{0}

        Let $A$ be a commutative group, with the group operations
expressed additively.  Thus left and right translations on $A$ are the
same.  If $d(x, y)$ is a $q$-semimetric on $A$ for some $q > 0$, then
$d(x, y)$ is invariant under translations on $A$ when
\begin{equation}
\label{d(x + a, y + a) = d(x, y)}
        d(x + a, y + a) = d(x, y)
\end{equation}
for every $a, x, y \in A$.  In this case, we have that
\begin{equation}
\label{d(-x, -y) = d(x, y)}
        d(-x, -y) = d(x, y)
\end{equation}
for every $x, y \in A$, using also the symmetry condition (\ref{d(x,
  y) = d(y, x)}).  If $d(x, y)$ is a translation-invariant
semi-ultrametric on $A$, then open and closed balls in $A$ with
respect to $d(x, y)$ centered at $0$ are subgroups of $A$, as before.

        If $A$ is a commutative topological group, then the corresponding
left and right uniformities on $A$ are the same.  One can check that
$A$ is a commutative topological group with respect to the topology
determined by any nonempty collection $\mathcal{M}$ of
translation-invariant $q$-semimetrics on $A$.  In this case, the
uniformity on $A$ just mentioned is the same as the one determined by
$\mathcal{M}$.

        Let $k$ be a field, and let $|\cdot|$ be a $q$-absolute value
function on $k$ for some $q > 0$, as in Section \ref{q-absolute value
functions}.  Thus
\begin{equation}
\label{d(x, y) = |x - y|, 2}
        d(x, y) = |x - y|
\end{equation}
defines a $q$-metric on $k$, which is obviously invariant under
translations on $k$ as a commutative group with respect to addition.
This leads to a topology on $k$ in the usual way, and it is easy to
see that $k$ is a topological group with respect to addition and this
topology.  If $|\cdot|$ is an ultrametric absolute value function on
$k$, then (\ref{d(x, y) = |x - y|, 2}) is an ultrametric on $k$, and
in particular $k$ is uniformly $0$-dimensional with respect to the
corresponding uniformity.  In this case, the open and closed balls in
$k$ of positive radius with respect to (\ref{d(x, y) = |x - y|, 2})
are open subgroups of $k$ with respect to addition.

        If $n$ is a positive integer, then let $n \cdot 1$ be the sum
of $n$ $1$'s in $k$, where $1$ refers to the multiplicative identity
element in $k$.  If $|\cdot|$ is an ultrametric absolute value function
on $k$, then it is easy to see that
\begin{equation}
\label{|n cdot 1| le 1}
        |n \cdot 1| \le 1
\end{equation}
for every $n \ge 1$, using (\ref{|1| = 1}).  A $q$-absolute value
function $|\cdot|$ on $k$ is said to be \emph{archimedian} if there
are positive integers $n$ such that $|n \cdot 1|$ is arbitrarily
large.  It suffices to have
\begin{equation}
\label{|n_0 cdot 1| > 1}
        |n_0 \cdot 1| > 1
\end{equation}
for some positive integer $n_0$, because
\begin{equation}
\label{|n_0^j cdot 1| = |(n_0 cdot 1)^j| = |n_0 cdot 1|^j}
        |n_0^j \cdot 1| = |(n_0 \cdot 1)^j| = |n_0 \cdot 1|^j
\end{equation}
for every positive integer $j$.  Otherwise, $|\cdot|$ is said to be
\emph{non-archimedian} on $k$ if there is a finite upper bound for $|n
\cdot 1|$ over all positive integers $n$, which implies that (\ref{|n
  cdot 1| le 1}) holds for every $n \ge 1$, by the preceding argument.

        Thus ultrametric absolute value functions are non-archimedian,
and it is well known that the converse also holds.  If $|\cdot|$ is an
archimedian $q$-absolute value function on $k$, then it follows in
particular that $k$ has characteristic $0$.  This leads to a natural
embedding of the field ${\bf Q}$ of rational numbers into $k$, and
$|\cdot|$ induces an archimedian $q$-absolute value function on ${\bf
  Q}$.  In this case, a famous theorem of Ostrowski implies that the
induced $q$-absolute value function on ${\bf Q}$ is the same as a
positive power of the standard absolute value function on ${\bf Q}$.
Another famous theorem of Ostrowski implies that $k$ can be identified
with a sub-field of the field ${\bf C}$ of complex numbers, using a
positive power of the standard absolute value function on ${\bf C}$.

        If $p$ is a prime number, then it is well known that the
\emph{$p$-adic absolute value} $|x|_p$ defines an ultrametric absolute
value function on ${\bf Q}$.  The \emph{$p$-adic numbers} are obtained
by completing ${\bf Q}$ with respect to the corresponding $p$-adic
metric.

        Let $V$ be a vector space over a field $k$, and let $N$ be a
$q$-seminorm on $V$ with respect to a $q$-absolute value function
$|\cdot|$ on $k$ for some $q > 0$.  Thus
\begin{equation}
\label{d(v, w) = N(v - w), 2}
        d(v, w) = N(v - w)
\end{equation}
defines a $q$-semimetric on $V$ which is invariant under translations
on $V$ as a commutative group with respect to addition.  If
$\mathcal{N}$ is a nonempty collection of $q$-seminorms on $V$, then
we get a collection $\mathcal{M}$ of translation-invariant
$q$-semimetrics on $V$ in this way.

\newpage


\printindex

\end{document}